\documentclass[12pt,english]{article}
\usepackage[T1]{fontenc}
\pdfoutput=1
\usepackage[latin9]{inputenc}
\usepackage[a4paper]{geometry}
\usepackage{color}
\usepackage{babel}
\usepackage{verbatim}
\usepackage{float}
\usepackage{algorithm2e}
\usepackage{amsmath}
\usepackage{amsthm}
\usepackage{amssymb}
\usepackage{graphicx}
\usepackage{csquotes}
\usepackage{tikz}
\usepackage{standalone}
\usepackage{bm}
\usepackage{pgfplots}
\pgfplotsset{compat=newest}
\pgfplotsset{plot coordinates/math parser=false}
\usepackage[unicode=true,pdfusetitle,
 bookmarks=true,bookmarksnumbered=false,bookmarksopen=false,
 breaklinks=false,pdfborder={0 0 0},pdfborderstyle={},backref=false,colorlinks=true]
 {hyperref}

\makeatletter
\theoremstyle{definition}
\newtheorem*{example*}{\protect\examplename}
\theoremstyle{plain}
\newtheorem{thm}{\protect\theoremname}
\theoremstyle{plain}
\newtheorem{lem}[thm]{\protect\lemmaname}
\theoremstyle{plain}

\makeatother

\providecommand{\examplename}{Example}
\providecommand{\lemmaname}{Lemma}
\providecommand{\theoremname}{Theorem}
\providecommand{\remarkname}{Remark}
\newcommand{\R}{\mathbb{R}}
\newcommand{\xx}{\mathbf{x}}
\newcommand{\y}{\mathbf{y}}

\usepackage[
backend=biber,
style=numeric,
sorting=nyt,
giveninits=true
]{biblatex}
\begin{document}
\title{Corrected trapezoidal rules for singular implicit boundary integrals}
\author{Federico Izzo\footnote{Corresponding author} \footnote{Department of Mathematics, KTH Royal Institute of Technology, Stockholm, Sweden (\href{mailto:izzo@kth.se}{izzo@kth.se})} \and Olof Runborg\footnote{Department of Mathematics, KTH Royal Institute of Technology, Stockholm, Sweden (\href{mailto:olofr@kth.se}{olofr@kth.se})}\and Richard Tsai\footnote{ Department of Mathematics and Oden Institute for Computational Engineering and Sciences, The University of Texas at Austin, Austin TX, USA (\href{mailto:ytsai@math.utexas.edu}{ytsai@math.utexas.edu})}}
\maketitle
\begin{abstract}
We present new higher-order quadratures for a family of boundary integral operators re-derived using the approach introduced in
[Kublik, Tanushev, and Tsai. \textit{J. Comp. Phys.} 247: 279-311, 2013].
In this formulation, a boundary integral over a smooth, closed hypersurface is transformed into an equivalent volume integral defined in a sufficiently thin tubular neighborhood of the surface. The volumetric formulation makes it possible to use the simple trapezoidal rule on uniform Cartesian grids and relieves the need to use parameterization for developing quadrature. 
Consequently, typical point singularities in a layer potential extend along the surface's normal lines. We propose new higher-order corrections to the trapezoidal rule on the grid nodes around the singularities. This correction is based on local decompositions of the singularity and is dependent on the angle of approach to the singularity relative to the surface's principal curvature directions. The proposed decomposition, combined with the volumetric formulation, leads to a special quadrature error cancellation.\\


\end{abstract}

\section{Introduction}

\label{sec:intro}

Boundary integral methods (BIMs) are employed in a wide range of applications for solving partial differential equations with conditions defined on boundaries of subregions and {at infinity}. In a BIM, one needs to solve a system of {boundary integral} equations (BIEs) involving singular integral operators acting on an unknown function defined on the boundaries and at infinity. 

Typical computational challenges for a boundary integral method involve developing high-order quadrature rules for the singular integrals and efficient dense matrix-vector computations for solving the resulting linear systems. Overcoming these challenges leads to highly efficient and accurate solutions for the associated partial differential equations.

We consider applications that require solving BIEs on a sequence of surfaces that are challenging to parametrize. These may include level set methods \cite{osher1988fronts,sethian1999level,osher2006level} and the closest point method \cite{ruuth2008simple,macdonald2008level, macdonald2010implicit} used to track evolving surfaces on a grid, particularly when the PDE solution is needed only at a small set far away from the surfaces.
In such situations, it is not immediately convenient to use any classical BIM.
The implicit boundary integral formulations are derived in \cite{kublik2013implicit}, aiming at these situations. 
We refer to this as the IBIM approach.
In \cite{tsaikublik16, kublik2018extrapolative}, further analysis related to the closest point projection is reported.
In \cite{chen2017implicit}, a similar formulation is derived to approximate the hypersingular integral equations arising from the Neumann problems of the Helmholtz equation. The method evaluates the limit of a family of surface integrals utilizing extrapolative averaging kernels.
An IBIM is applied to compute electrostatic potential from large molecules submerged in a solvent in \cite{zhong2018implicit}. That paper also demonstrates that IBIM, coupled with an ``off-the-shelf'' Fast Multipole Method, can easily be applied to solve the equations for very large molecules.
Partial differential equations arising from calculus of variation problems defined on closed surfaces can be solved with high order convergence rates  using similar strategies; see \cite{chu2018volumetric,hsu2019coupled, martin2020equivalent, martin2020equivalent-jsc}.

In the center of these formulations lie volume integrals with identical evaluations to the corresponding surface integrals. 
The volume integrals involve integration in thin tubular neighborhoods of the surfaces in the ambient space, and do not require surface parameterizations.
In principle, the integrals can be approximated on a wide range of meshing.
Among the existing work, these volume integrals are discretized on Cartesian grids using the trapezoidal rule and a lower order regularization of the layer singularities.

This paper presents higher-order accurate quadrature rules for the singular integrals arising from the non-parametric boundary integral formulation discussed above. 
In these formulations, the singularities in the integral operators concentrate along the surface normal lines; these lines generally do not lie on the grid. This feature is atypical in the more classical boundary integral formulations.
Our approach
is a generalization of the methods in \cite{marin2014corrected}: Regular trapezoidal rule-based summation is performed over the grid nodes lying in the regions, excluding small, grid-dependent neighborhoods around the singularities. This approach is called the punctured trapezoidal rule. We derive additional corrections corresponding to the skipped grid nodes and add them to the punctured trapezoidal rule. The resulting corrected trapezoidal rule is second-order accurate with respect to the uniform grid spacing of the underlying Cartesian grid.
We discover an additional benefit of the non-parametric approach --- there are cancellation of errors which leads to a improved order of accuracy in practice. \\

The structure of the paper is as follows: in Section 2 we present an overview of how to express solutions to Laplace and Helmholtz problems using boundary integral equations, and introduce the volumetric extension setting to express surface integrals via volume integrals. In Section 3 we present singularity regularization methods for volume integrals of Section 2. In Section 4 we present a general singularity correction framework for the volume integrals of Section 2, and in Section 5 we go into details about how to apply this methods to the Laplace singular kernels. Finally, Section 6 presents numerical results for the methods of Sections 3 and 5 applied to the evaluation of Laplace potentials.

\section{A short review of boundary integral equations}

Boundary integral methods can be used to solve 
the Laplace and homogeneous Helmholtz equations
in both bounded and unbounded domains. Given a bounded domain $D\subset\mathbb{R}^{3}$,
%$n=2,3$, 
the problems are
\[
\begin{cases}
-\Delta u=0 & \text{in }\Omega\\
u=f\ \text{ or }\ \frac{\partial u}{\partial \mathbf{n}}=g & \text{on }\partial\Omega
\end{cases}\quad,\quad\begin{cases}
\Delta u+\lambda^{2}u=0 & \text{in }\Omega\\
u=f\ \text{ or }\ \frac{\partial u}{\partial \mathbf{n}}=g & \text{on }\partial\Omega
\end{cases}\quad,\ \ 
\]
where $\Omega=D$ for the interior problem, and $\Omega=\mathbb{R}^{3}\setminus\overline{D}$ {complement of the closure of $D$}
for the exterior problem, and $\mathbf{n}$ is the outward pointing normal
to the surface $\partial\Omega=:\Gamma$. For the exterior Helmholtz
problem, in order to ensure uniqueness (see \S3, Thm 3.13 in \cite{colton2013integral}), the solution must also satisfy
the radiation (or Sommerfeld) condition

\begin{equation}
\lim_{\|\mathbf{x}\|\to\infty}\|\mathbf{x}\|\left(\frac{\partial u}{\partial r}(\mathbf{x})-\text{i}\lambda u(\mathbf{x})\right)=0\ \ ,\ \ r=\|\mathbf{x}\|\,.\label{eq:Helm_radiation_cond}
\end{equation}%^{\frac{n-1}{2}}
The Helmholtz equation with $\lambda=0$ becomes Laplace equation, and the
two problems are heavily related. In three dimensions, the fundamental solutions for Laplace and Helmholtz are respectively:
\begin{align*}
G_{0}(\mathbf{x},\mathbf{y})=\frac{1}{4\pi}\frac{1}{\|\mathbf{x}-\mathbf{y}\|}\ ,\ \ \ G_{\lambda}(\mathbf{x},\mathbf{y})=\frac{1}{4\pi}\frac{\exp(\text{i}\lambda\|\mathbf{x}-\mathbf{y}\|)}{\|\mathbf{x}-\mathbf{y}\|}\,.
\end{align*}

\subsection{Solutions as layer potentials}

A solution $u$ to $\Delta u+\lambda^2u=0$ in 
{$\Omega=\R^3\setminus \overline{D}$}, can be expressed, for $x\in \Omega$, as
\begin{align}
\text{(Single-Layer potential) }\ \ u(\mathbf{x})=&\mathcal{S}[\alpha](\mathbf{x}):=  \int_{\partial\Omega}G_{\lambda}(\mathbf{x},\mathbf{y})\alpha(\mathbf{y})\text{d}\sigma_{\mathbf{y}}\label{eq:SL},\\
\text{(Double-Layer potential) }\ \ u(\mathbf{x})=&\mathcal{K}[\beta](\mathbf{x}):=  \int_{\partial\Omega}\frac{\partial G_{\lambda}}{\partial \mathbf{n}_{{y}}}(\mathbf{x},\mathbf{y})\beta(\mathbf{y})\text{d}\sigma_{\mathbf{y}},\label{eq:DL}\\
\text{(Combined-Layer potential) }\ \ u(\mathbf{x})=&\mathcal{C}[\zeta](\mathbf{x}):=  \ \mathcal{S}[\zeta](\mathbf{x})+\text{i}\xi\mathcal{K}[\zeta](\mathbf{x}).\label{eq:CL}
\end{align}
The functions $G_{\lambda}$ and $\frac{\partial G_{\lambda}}{\partial \mathbf{n}_{y}}$ are called \emph{single-layer} {(SL)} and \emph{double-layer} {(DL)} \emph{kernels} respectively. The functions $\alpha$, $\beta$, $\zeta$, are called \emph{single-layer density}, \emph{double-layer density}, and \emph{combined-layer density}.
Their expressions are derived by the application of the Green's identities and the properties of the solutions to interior and exterior Helmholtz problems (see \S2.4, 2.5, 3.2 and 3.4 in \cite{colton2013integral}): for $\mathbf{y}\in\partial\Omega$
\begin{align*}
    \alpha(\mathbf{y})&:= \dfrac{\partial v_{SL}}{\partial \mathbf{n}}(\mathbf{y})-\dfrac{\partial u}{\partial \mathbf{n}}(\mathbf{y}) \,,\ \ \ \beta(\mathbf{y}):= u(\mathbf{y})-v_{DL}(\mathbf{y})\,, \\
    \zeta(\mathbf{y})&:=\dfrac{\partial v_{CL}}{\partial \mathbf{n}}(\mathbf{y})-\dfrac{\partial u}{\partial \mathbf{n}}(\mathbf{y})=\frac{1}{\text{i}\xi}(u(\mathbf{y})-v_{CL}(\mathbf{y}))\,,
\end{align*}
where $v_{SL},\,v_{DL},\,v_{CL}$ are solutions to $\Delta v+\lambda^2v=0$ in $\Omega^c={D}$ with boundary conditions:
\begin{align*}
    v_{SL}=u\,,\ \ \ v_{DL}=\dfrac{\partial u}{\partial\mathbf{n}}\,,\ \ \ v_{CL}+\text{i}\xi\dfrac{\partial v_{CL}}{\partial\mathbf{n}}=u+\text{i}\xi\dfrac{\partial u}{\partial\mathbf{n}}\,.
\end{align*}
The uniqueness of $v_{CL}$ requires $\xi\neq 0$. 
%The parameter i$\xi$ is imaginary in order to make the solution $v_{CL}$ unique, but the value of $\xi$ itself is arbitrary. 
In practice, the value for $\xi$ is usually tuned to improve the properties of the numerical methods.\\

When expressing the solution to the problem as a layer potential,
the density is unknown. Using the boundary conditions, we can find
%\CANC{boundary integral equations (}
BIEs
%\CANC{)}
to which the solution is the density.
In the Dirichlet problem with boundary conditions $u=f$ on $\partial\Omega$
the solution can be expressed as either a single-layer potential 
\begin{equation}
\int_{\partial\Omega}G_{\lambda}(\mathbf{x},\mathbf{y})\alpha(\mathbf{y})\text{d}\sigma_{\mathbf{y}}=f(\mathbf{x})\ \ ,\ \ \mathbf{x}\in\partial\Omega\label{eq:SL_dirich}
\end{equation}
or a double-layer potential:\begin{equation}
\int_{\partial\Omega}\frac{\partial G_{\lambda}}{\partial \mathbf{n}_{y}}(\mathbf{x},\mathbf{y})\beta(\mathbf{y})\text{d}\sigma_{\mathbf{y}}\mp\frac{1}{2}\beta(\mathbf{x})=f(\mathbf{x})\ \ ,\ \ \mathbf{x}\in\partial\Omega\label{eq:DL_dirich}
\end{equation}
with minus for the interior and plus for the exterior problem. 
Note that in general the double-layer formulation is preferable as it involves
the solution of an integral equation of the second kind: this leads to non-singular matrices when discretizing the integral operators with Nystr\"om methods.

For the Neumann problem with boundary conditions
$\frac{\partial u}{\partial \mathbf{n}}=g$ on $\partial\Omega$, the single-layer
formulation is preferable as it avoids the appearance of a hypersingular
kernel, and the boundary integral equations to solve is
\begin{equation}
-\int_{\partial\Omega}\frac{\partial G_{\lambda}}{\partial \mathbf{n}_{{x}}}(\mathbf{x},\mathbf{y})\alpha(\mathbf{y})\text{d}\sigma_{\mathbf{y}}\pm\frac{1}{2}\alpha(\mathbf{x})=g(\mathbf{x})\ ,\ \ \mathbf{x}\in\partial\Omega\label{eq:SL_neu}
\end{equation}
with plus for the interior and minus for the exterior problem. 

In addition to the single- and double-layer potentials, we also consider
the potential appearing in \eqref{eq:SL_neu}, called the double-layer conjugate potential:
\begin{equation}\label{eq:DLC}
\text{(Double-Layer Conjugate potential) }\ \ \mathcal{K}^{*}[\alpha](\mathbf{x}):=\int_{\partial\Omega}\frac{\partial G_{\lambda}}{\partial \mathbf{n}_{{x}}}(\mathbf{x},\mathbf{y})\alpha(\mathbf{y})\text{d}\sigma_{\mathbf{y}}\,.
\end{equation}
Its treatment is going to be analogous to the treatment of the
double-layer potential. The function $\frac{\partial G_{\lambda}}{\partial \mathbf{n}_{x}}$ is called \emph{double-layer conjugate} {(DLC)} \emph{kernel}.

It is important to note that these boundary integral equations (\ref{eq:SL_dirich}
-\ref{eq:SL_neu}) are valid for both 
Laplace and Helmholtz equation, and the only thing that changes in the formulae
is the parameter $\lambda$ (the wavenumber), and consequently the kernel $G_{\lambda}$. 

\subsection{Quadratures for singular integrals}
If the parametrization of a surface is known, the surface integration can be done straightforwardly by applying any preferred quadrature rule: given the parametrization $\mathbf{z}$ of $U\subset\partial\Omega$, 
$\mathbf{z}(\tau,\varsigma):=(x(\tau,\varsigma),y(\tau,\varsigma),z(\tau,\varsigma))$, $(\tau,\varsigma)\in W$,
\[
\int_{U}F(\mathbf{z})\text{d}\sigma_{\mathbf{z}}=\int_{W}F(\mathbf{z}(\tau,\varsigma))J(\tau,\varsigma)\text{d}\tau\text{d}\varsigma\approx\sum_{i,j}\omega_{ij}F(\mathbf{z}(\tau_{i},\varsigma_{j}))J(\tau_{i},\varsigma_{j})\,,
\] 
where $J(\tau,\varsigma)=|\mathbf{z}_\tau \times \mathbf{z}_\varsigma|$ is the surface area element, $\{(\tau_i,\varsigma_j)\}_{i,j}\subset\mathbb{R}^2$ are the nodes, and $\{\omega_{ij}\}_{i,j}\subset\mathbb{R}$ are the weights.

If the function $F$ is smooth, the freedom of choice of the position of the nodes and of the quadrature rule makes it easy to attain high accuracy. However, if the function $F$ is singular, for example in the origin, then using standard quadrature rules results in a great loss of accuracy. To remedy this loss, several classes of methods have been developed, for different kinds of singular integrands.\\

An important class of methods to handle these singular integrands, the methods of \textit{singularity subtraction}, approaches the problem by locally approximating the surface and evaluating the integral analytically, and then adding a correction term dependent on the surface approximation. If the kernels are similar to the ones for which the analytical results exist, singularity subtraction is applied and those same results are used \cite{farina2001evaluation, davis2007methods}.
  
Another class of methods, the methods of \textit{singularity cancellation}, use a change of variables to put themselves in a setting where it is possible to split the integral in a smooth one and a singular one which is only defined close to the singularity points, and can be computed to high accuracy using exact local parametrization and a suitable quadrature rule, e.g. trapezoidal rule in polar coordinates \cite{bruno2001fast,ying2006high}.
   
A third class of methods, the methods of \textit{singularity regularization}, relies on regularizing the kernel so that rules for smooth integrands can be applied, and then adding corrections to account for the different kernel based on analytical results, or on Richardson extrapolation \cite{beale2001method,haroldsen1998numerical,kublik2013implicit}.
  
A newer class of methods, called \emph{quadrature by expansion} (QBX), handles the problem by expanding and treating the kernel (and the corresponding layer potential) away from the surface target point, e.g. with Taylor or spherical harmonics, consequently working with smooth integrands, and then evaluate the results back on the surface. It relies on the smoothness of the expansion terms because the new target point is not on the surface, and on the convergence of the expansion in the surface target point \cite{qbxklockner2013quadrature}. 
  
Finally, the methods of \textit{singularity correction} aim to develop specialized quadrature rules to deal with families of singular integrands by modifying the weights of an existing quadrature rule, often trapezoidal rule, close to the singularity point \cite{kapur1997high,marin2014corrected,wumartinsson2020zeta}. Marin, Tornberg and Runborg \cite{marin2014corrected} developed corrections to the trapezoidal rule for singularities of the kind $|{x}|^{\gamma}$, $\gamma\in(-1,0)$, in one dimension and $\|\mathbf{x}\|^{-1}$ in two dimensions, proved convergence order, and found an analytic expression for the weights in one dimension. Wu and Martinsson expanded these results to $\log|x|$ in one dimension \cite{wumartinsson2020zeta} and found analytic expression for weights in one and two dimensions \cite{wumartinsson2020corrected}. \\
  
The majority of above-mentioned methods require explicit knowledge of the parametrization, and the possibility to choose the position of the nodes around the singularity; moreover, often the singularity point lies in one of the nodes. 
  
In our setting however, the position of the nodes on the surface is going to be determined by the projections of the nodes in the volume onto the surface, which cannot be assumed to have any particular structure (see Figure~\ref{fig:tiltedtorusplot}). Moreover, because of how the integrand is extended from the surface to the volume, instead of a single singularity point on the surface, we will have the singularity lying along a straight line in three dimensions.

We will approach this problem then by splitting the three-dimensional trapezoidal rule into the weighted sum of all the two-dimensional trapezoidal rules on each two-dimensional grid and correct each one separately.

\subsection{Volumetric extensions of the layer integrals}
\label{subsec:ImplicitParametrization}

Let $\Omega \subset\mathbb{R}^n$ be a bounded open set with $C^2$ boundaries, and $\partial\Omega=:\Gamma$. We shall refer to $\Gamma$ as the surface. Let $f$ be a function defined on $\Gamma$ (or $\mathbb{R}^n$).
In this Section we present an approach for extending 
a boundary integral
\begin{equation}\label{eq:surfint}
\int_{\Gamma}f(\mathbf{x})\text{d}\sigma_{\mathbf{x}},
\end{equation}
to a volumetric integral around the surface.
Instead of parameterizations, this approach relies on the Euclidean distance to the surface, and its derivatives. More precisely, we define the \textit{signed distance function}
\begin{equation}\label{eq:signed-distance}
    d_{\Gamma}(\mathbf{x}):=\begin{cases}
\quad\min_{\mathbf{y}\in\Gamma}\left\|\mathbf{x}-\mathbf{y}\right\|, & \text{ if } \mathbf{x}\in \Omega\\
-\min_{\mathbf{y}\in\Gamma}\left\|\mathbf{x}-\mathbf{y}\right\|, & \text{ if } \mathbf{x}\in \Omega^c 
\end{cases}
\end{equation}
and the \textit{closest point projection}  %\textbf{add more ref (CPM,etc)}
\begin{equation}\label{eq:closest-point-projection}
    P_{\Gamma}(\mathbf{x}):=\text{argmin}_{\mathbf{y}\in\Gamma}\left\|\mathbf{x}-\mathbf{y}\right\|\,.
\end{equation}
If there is more than one global minimum, we pick one randomly from the set. 
Let $\mathcal{C}_\Gamma$ denote the set of points in $\mathbb{R}^n$ which are equidistant to at least two distinct points on $\Gamma$.
The \emph{reach} $\tau$ is defined as $\inf_{\mathbf{x}\in \Gamma, \mathbf{y}\in\mathcal{C}_\Gamma} \|\mathbf{x}-\mathbf{y}\|$. Clearly, $\tau$ is restricted by the local geometry (the curvatures) and the global structure of $\Gamma$ (the Euclidean and geodesic distances between any two points on $\Gamma$). 

In this paper, we assume that $\Gamma$ is $C^2$ and has a non-zero reach. 
Let $T_{\varepsilon}$ denote the set of points of distance at most $\varepsilon$ from $\Gamma$:
\begin{equation}
    T_{\varepsilon}=\{\mathbf{x}\in\mathbb{R}^{n}:\ |d_{\Gamma}(\mathbf{x})|\leq \varepsilon\}.
\end{equation}
Then, for $\varepsilon<\tau$, $P_\Gamma$ is a diffeomorphism between the level sets of $d_\Gamma$ and
\[
P_\Gamma(\mathbf{x}) = \mathbf{x}-d_\Gamma(\mathbf{x})\nabla d_\Gamma(\mathbf{x}),~~~~\mathbf{x}\in T_\varepsilon.
\]

We define the \emph{extension} (or \emph{restriction}) of the integrand $f$ by
\begin{equation}\label{eq:extension-restriction-f}
\overline{f}(\mathbf{x}) := f(P_\Gamma \mathbf{x}),~~~\mathbf{x}\in\mathbb{R}^n.
\end{equation}

As in \cite{kublik2013implicit, tsaikublik16}, we can then rewrite the
surface integral (\ref{eq:surfint}), 
for any $\eta\in[-\varepsilon,+\varepsilon]$,
as % of a given function $f$ on $\Gamma$ as 
\begin{equation}\label{eq:surfint2}
\int_{\Gamma}f(\mathbf{x})\text{d}\sigma_x=\int_{\Gamma_\eta}\overline{f}(\mathbf{x}')J_{\eta}(\mathbf{x}')\text{d}\sigma_{\mathbf{x}'},
\end{equation}
where $J_{\eta}(\mathbf{x}')$ is the Jacobian of the transformation from  $\Gamma$ to the level set $\Gamma_{\eta}:=\{\mathbf{x}\in\mathbb{R}^{n}\ :\ d_{\Gamma}(\mathbf{x})=\eta\}.$ In $\mathbb{R}^3$, the Jacobian $J_\eta(\mathbf{x}')$ is a quadratic polynomial in $\eta$:
\begin{equation*}
    J_\eta(\mathbf{x}'):=1+2\eta \mathcal{H}(\mathbf{x}') + \eta^2 \mathcal{G}(\mathbf{x}') = \sigma_1\sigma_2(P^\prime_\Gamma \mathbf{x}'),
\end{equation*}
where $\mathcal{H}(\mathbf{x}')$ and $\mathcal{G}(\mathbf{x}')$ are respectively the mean and Gaussian curvatures of $\Gamma_\eta$ at $\mathbf{x}'$, and $\sigma_1\sigma_2(P^\prime_\Gamma \mathbf{x}')$ is the product of the first two singular values of the Jacobian matrix of $P_\Gamma$ evaluated at $\mathbf{x}'$. See \cite{tsaikublik16} for more detail.

\begin{figure}
    \begin{center}
    \includegraphics{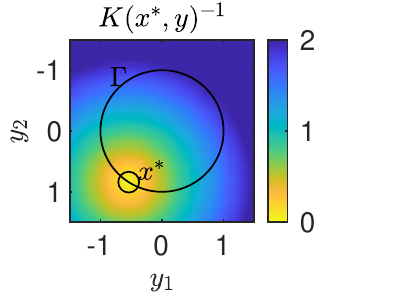}\includegraphics{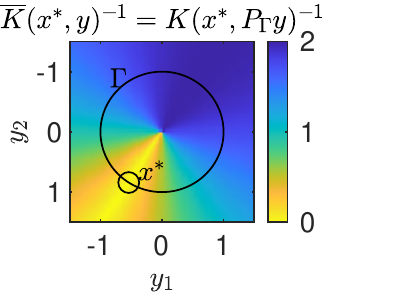}\includegraphics{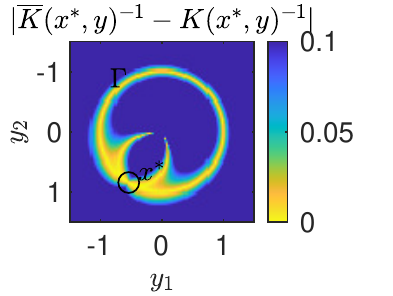}\vspace{0.2cm}\includegraphics[height=3cm]{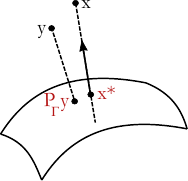}
    \end{center}
    \caption{Kernel restriction}
    \footnotesize{Visualization of the restriction of $K(\xx^*,\y)=\|\xx^*-\y\|^{-1}$ to the unit circle and of the singular properties of $\overline{K}$. 
    Since $K(\mathbf{x}^*, \mathbf{y})$ is singular at $\mathbf{y}^*$, $\overline{K}(\mathbf{x}^*, \mathbf{y})$ is singular along the line $P_\Gamma \mathbf{y}=\mathbf{x}^*.$ One observes that the gradient of $\overline{K}(\mathbf{x}^*, \mathbf{y})^{-1}$ (and thus $\overline{K}(\mathbf{x}^*, \mathbf{y})$) is orthogonal to the normal of the interface.}
    \label{fig:singularity_along_normals}
\end{figure}

To extend (\ref{eq:surfint}) to a volumetric integral, 
we now average the integral on the right hand side in 
(\ref{eq:surfint2})
over $\eta$ ranging from $-\varepsilon$ to $\varepsilon,$  using 
\begin{equation*}
    \delta_\varepsilon(\eta):= \frac{1}{\varepsilon}\phi\left(\frac{\eta}{\varepsilon}\right),
\end{equation*}
with $\phi\in C^\infty(\mathbb{R})$ supported in $[-1,1]$, and $\int_\R \phi(x)\text{d}x =1$. This means
\begin{align*}
\int_{\Gamma}f(\mathbf{x})\text{d}\sigma_{\mathbf{x}}&=\int_{-\varepsilon}^{+\varepsilon}\left\{ \delta_{\varepsilon}(\eta) \int_{\Gamma_\eta}f(P_\Gamma\mathbf{x})J_{\eta}(\mathbf{x})\text{d}\sigma_{\mathbf{x}}\right\}\text{d}\eta .
\end{align*}
Applying the coarea formula, we have
\begin{align*}
\int_{-\varepsilon}^{+\varepsilon}\left\{ \delta_{\varepsilon}(\eta) \int_{\Gamma_\eta}f(P_\Gamma\mathbf{x})J_{\eta}(\mathbf{x})\text{d}\sigma_{\mathbf{x}}\right\}\text{d}\eta  = \int_{T_\varepsilon}f(P_\Gamma\mathbf{x})J_{d_{\Gamma}(\mathbf{x})}(\mathbf{x})\delta_{\varepsilon}(d_{\Gamma}(\mathbf{x}))\text{d}\mathbf{x}.
\end{align*}
Thus from the surface integral, we derive a volume integral with the same evaluation:
\begin{equation}
    \int_{\Gamma}f(\mathbf{x})\text{d}\sigma_{\mathbf{x}}=\int_{\mathbb{R}^n}f(P_\Gamma\mathbf{x})\delta_{\Gamma,\varepsilon}(\mathbf{x})\text{d}\mathbf{x}\,,
\end{equation}
where 
\[
\delta_{\Gamma,\varepsilon}(\mathbf{x}):=J_{d_{\Gamma}(\mathbf{x})}(\mathbf{x})\delta_{\varepsilon}(d_{\Gamma}(\mathbf{x}))\ ,\ \ \mathbf{x}\in \R^n\,.
\]
Our primary focus is when $f(\mathbf{x})$ is replaced by a function $K(\mathbf{x},\mathbf{y})\zeta(\mathbf{y})$, 
{with $K:\R^n\times\R^n\to\R$ and} $K(\mathbf{x},\mathbf{y})$ 
 singular for $\mathbf{x}=\mathbf{y}$ (corresponding to the layer potentials reviewed in the previous Section):
\begin{equation}\label{eq:general_S}
\mathcal{J}[\zeta](\mathbf{x}):=\int_{\Gamma}K(\mathbf{x},\mathbf{y})\zeta(\mathbf{y})\text{d}\sigma_{\mathbf{y}}\ \ ,\ \ \ \mathbf{x}\in\R^n.
\end{equation}
%\RTcorr{
When a function $g:\Gamma\to\R$ is given, we may form an integral equation for the unknown density $\zeta$. For example, in the case of the double-layer potential \eqref{eq:DL_dirich}, the equation is
\begin{equation}\label{eq:BIE}
\mathcal{J}[\zeta](\mathbf{x}) \mp\frac{1}{2}\zeta(\mathbf{x}) = g(\mathbf{x}),~~~\mathbf{x}\in\Gamma.
\end{equation}
{Suppose that for any $\mathbf{x}$, we are interested in evaluating $K(\mathbf{x},\mathbf{y})$ at the point on $\Gamma$ that is closest to $\mathbf{y}$. This can be done by   }
\begin{equation}\label{eq:restriction_of_K}
    \overline{K} (\mathbf{x},\mathbf{y}):=K(\mathbf{x},P_\Gamma \mathbf{y})\ ,\ \ \mathbf{x},\mathbf{y}\in\R^n\,.%J_\Gamma (y) \delta_\varepsilon(d_\Gamma(y)).
\end{equation}
{Hence we refer to $\overline{K}(\mathbf{x},\mathbf{y})$ as the restriction of $K$.}
If $K(\mathbf{x},\mathbf{y})$ is singular for $\mathbf{x}=\mathbf{y}$, then $\overline{K}(P_\Gamma\mathbf{x},\mathbf{y})$ is singular on the set
\begin{equation}\nonumber
    \{ (\mathbf{x},\mathbf{y})\in \mathbb{R}^n\times\mathbb{R}^n: P_\Gamma \mathbf{x} = P_\Gamma \mathbf{y}\}, 
\end{equation}
i.e. for a fixed $\mathbf{x}^*\in \Gamma$, $\overline{K}(\mathbf{x}^*,\mathbf{y})$ is singular along the normal line passing through $\mathbf{x}^*$,
while $K(\mathbf{x}^*,\mathbf{y})$ is singular in a point.
In Figure~\ref{fig:singularity_along_normals} the singular behavior of $\overline{K}(\mathbf{x}^*,\mathbf{y})$ along the normal is illustrated.

In conclusion, instead of approximating \eqref{eq:general_S}, we approximate
\begin{equation}\label{eq:volume_potential_SL}
\overline{\mathcal{J}}[\rho](\mathbf{x}):= \int_{\mathbb{R}^n} \overline{K}(\mathbf{x},\mathbf{y})\rho(\mathbf{y}) \delta_{\Gamma, \varepsilon}(\mathbf{y}) \text{d}\mathbf{y},
\end{equation}
for functions $\rho$ that are integrable  in $T_\varepsilon$. 

Corresponding to \eqref{eq:BIE}, we have the equivalent implicit boundary integral equation
\begin{equation}\label{eq:IBIE}
\overline{\mathcal{J}}[\rho](P_\Gamma\mathbf{x}) \mp \frac{1}{2}\rho(\mathbf{x})= {g}(P_\Gamma\mathbf{x}),~~~\mathbf{x}\in T_\varepsilon. 
\end{equation}
{ The solution $\rho$ will coincide with the constant extension along the normals of $\zeta$. To see this, we write the two equations:
\begin{align*}
    \left(\overline{\mathcal{J}}\mp \frac{1}{2}I\right)[\zeta\circ P_\Gamma](\mathbf{x})=g(P_\Gamma \xx)\ \ ,\ \ \
    \left(\overline{\mathcal{J}}\mp \frac{1}{2}I\right)[\rho](\mathbf{x})=g(P_\Gamma \xx)\ \ ,\ \ \   \xx\in T_{\varepsilon}\,.
\end{align*}
The first equation corresponds to \eqref{eq:BIE} where the integral has been rewritten and the target point $\xx\in T_\varepsilon$ is projected onto $\Gamma$. The second is \eqref{eq:IBIE} where the equation is imposed for $\rho$ function defined in $T_\varepsilon$. If we take the difference of these two equations, we find
\[
\left(\overline{\mathcal{J}}\mp \frac{1}{2}I\right)[\rho-\zeta\circ P_\Gamma](\mathbf{x})=0\ \ ,\ \ \ \xx\in T_\varepsilon\,.
\]
The kernel of the operator on the left-hand side coincides with the kernel of the original operator, so whenever the solution is unique, $\rho(\xx)=\zeta(P_\Gamma \xx)$ for any $\xx\in T_\varepsilon$.}\\ 

In this paper, we will concentrate on developing numerical quadratures for the extended singular integral operator $\overline{J}[\rho](\mathbf{x})$ for $\mathbf{x}\in\Gamma$ (equivalently, $\overline{J}[\rho](P_\Gamma\mathbf{x})$ for $\mathbf{x}\in T_\varepsilon)$. 
The quadrature rules will be constructed based on the trapezoidal rule for the grid nodes $T^h_\varepsilon := T_\varepsilon \cap h\mathbb{Z}^3$, which corresponds to the portion of the uniform Cartesian grid $h\mathbb{Z}^3$ within $T_\varepsilon$.
Since the integrand in \eqref{eq:volume_potential_SL} is singular for $\mathbf{x}\in\Gamma$, the trapezoidal rule should be corrected near 
$\mathbf{x}$ for faster convergence. 
{Correction will be defined by summing the judiciously derived weights over a set of grid nodes denoted by $N_h(\mathbf{x}).$ The sum will be denoted by $\mathcal{R}_h(\mathbf{x})$. }
Ultimately, the quadrature for $\overline{J}[\rho](P_\Gamma\mathbf{x})$ 
will involve the regular Riemann sum of the integrand in $T^h_\varepsilon\setminus N_h(\mathbf{x})$, and  the correction $\mathcal{R}_h(\mathbf{x})$ in $N_h(\mathbf{x})$:

\begin{equation}\label{eq:corrected_trapezoidal_general_form}
  \overline{\mathcal{J}}[\rho](\mathbf{x}) \approx 
   \sum_{\mathbf{y}_m\in  T_\varepsilon^h \setminus N_h(\mathbf{x})}\overline{K}(\mathbf{x},\mathbf{y}_{m})\rho(\mathbf{y}_{m})\delta_{\Gamma, \varepsilon}(\mathbf{y}_m)h^3+ \mathcal{R}_h(\mathbf{x})\hspace{0.2cm}.
\end{equation}
The contribution of this paper is a high order, trapezoidal rule-based, quadrature rule for $\mathcal{J}$ via $\overline{\mathcal{J}}$.

Figure~\ref{fig:tiltedtorusplot} demonstrates a typical configurations the points $\mathbf{y}_m$ in the summation for a torus. 

In the following two Sections, we will see how to build the correction term $\mathcal{R}_h(\mathbf{x})$ using two different approaches: a function regularization independent of trapezoidal rule (Section \ref{subsec:Laplace_regularizations}), and the corrected trapezoidal rule (Section \ref{sec:corr_TR}). Each approach will determine the {set} $N_h(\mathbf{x})$ differently.

\section{Correction via regularization of singularity}

\label{subsec:Laplace_regularizations}

In this Section we present an approach that locally regularizes a singular kernel before the discretization.  In this approach, a special Lipschitz continuous function, $\Psi$, is used to replace the kernel in a neighborhood around the kernel's singularity. 
The integral with the regularized kernel is then extended following \eqref{eq:volume_potential_SL}.  
Again, the resulting implicit boundary integral can be discretized on different meshings. {When the trapezoidal rule is applied to $\overline{J}$ involving the locally regularized $\overline{K}$, we find an expression of the kind \eqref{eq:corrected_trapezoidal_general_form} where the term $\mathcal{R}_h(\xx)$ involves a local sum of the integrand with $\overline{K}$ replaced by $\overline{\Psi}$.}

\subsection{A localized regularization approach}
We consider regularization of $K(\mathbf{x},\mathbf{y})$ 
constructed in the following fashion:
\begin{equation}\label{eq:kernel_regularization_general}
    K_{r_0}^{reg}(\mathbf{x},\mathbf{y}) =
    \begin{cases}
      K(\mathbf{x},\mathbf{y}), & \|\mathbf{x}-\mathbf{y}\|\geq r_0, \\
      \Psi_{\Gamma, r_0}(\mathbf{x},\mathbf{y}), &
      \|\mathbf{x}-\mathbf{y}\|< r_0,
    \end{cases}
\end{equation}
where $\Psi_{\Gamma,r_0}$ thus is a function which substitutes $K$ close to the singularity point.
We choose this as a simple function (constant, or linear in $\|\mathbf{x}-\mathbf{y}\|$) which
approximates $K(\mathbf{x},\cdot\,)$ weakly for $C^1$ functions on the $r_0$ neighborhood of $\mathbf{x}$, such that
\begin{equation}\label{weak-approx}
    \int_{\Gamma\cap B_{r_0}(\mathbf{x})} \Psi_{\Gamma, r_0}(\mathbf{x},\mathbf{y}) \rho(\mathbf{y}) \text{d}\sigma_\mathbf{y} \approx \int_{\Gamma\cap B_{r_0}(\mathbf{x})} K(\mathbf{x},\mathbf{y}) \rho(\mathbf{y}) \text{d}\sigma_\mathbf{y},~~~\rho\in C^1(\Gamma).
\end{equation}
Here, $B_{r_0}(\mathbf{x})$ is the ball with radius $r_0$, centered at $\mathbf{x}$.

Applying the trapezoidal rule to the integral \eqref{eq:volume_potential_SL} with the regularized kernel~\eqref{eq:kernel_regularization_general}, we get a correction to the trapezoidal rule of the form \eqref{eq:corrected_trapezoidal_general_form}:

\begin{align}\label{eq:Kreg_TR}
    \overline{\mathcal{J}}[\rho](\mathbf{x}) \approx & \ h^3\sum_{\mathbf{y}_m\in T_\varepsilon^h}\overline{K}_{r_0(h)}^{reg}(\mathbf{x},\mathbf{y}_m)\rho(\mathbf{y}_m)\delta_{\Gamma, \varepsilon}(\mathbf{y}_m)  \\
    = & \ \ \, \ h^3\sum_{\mathbf{y}_m\in T_\varepsilon^h\setminus N_h(\mathbf{x})}\overline{K}(\mathbf{x},\mathbf{y}_m)\rho(\mathbf{y}_m)\delta_{\Gamma, \varepsilon}(\mathbf{y}_m)  \nonumber \\
    & + h^3\sum_{\mathbf{y}_m\in T_\varepsilon^h\cap N_h(\mathbf{x})}\overline{\Psi}_{\Gamma,r_0(h)}( \mathbf{x},\mathbf{y}_m)\rho(\mathbf{y}_m)\delta_{\Gamma, \varepsilon}(\mathbf{y}_m) \nonumber
\end{align}
where: $\overline{K}$, $\overline{K}^{reg}_{r_0}$, and $\overline{\Psi}_{\Gamma,r_0}$ are the restrictions \eqref{eq:restriction_of_K} of $K$, $K_{r_0}^{reg}$, and $\Psi_{\Gamma,r_0}$ respectively. The formula \eqref{eq:corrected_trapezoidal_general_form} holds with:
\begin{align*}
N_h(\mathbf{x}) &= \{ \mathbf{y}\in T^h_\varepsilon\,:\, \|P_\Gamma \mathbf{y}-P_\Gamma \mathbf{x}\|<r_0(h) \}\,, \\
\mathcal{R}_h(\mathbf{x}) &= h^3\sum_{\mathbf{y}_m\in T_\varepsilon^h\cap N_h(\mathbf{x})}\overline{\Psi}_{\Gamma,r_0(h)}( \mathbf{x}, \mathbf{y}_m)\rho(\mathbf{y}_m)\delta_{\Gamma, \varepsilon}(\mathbf{y}_m)\,.
\end{align*} 
 
In order to determine $\Psi_{\Gamma,r_0}$ we want it to satisfy (\ref{weak-approx}) for $\rho\equiv 1$ with an error at most $O(r_0^2)$. However, given the lack of an explicit parametrization of the surface, we  
approximate $\Gamma$ in the integrals in \eqref{weak-approx} by a suitable paraboloid, $\tilde\Gamma_{\mathbf{x}}$, defined from the principal curvatures of $\Gamma$ at $\mathbf{x}$ (as shown in \cite{kublik2013implicit}). The domain $\Gamma\cap B_{r_0}(\mathbf{x})$
is furthermore replaced by a neighborhood $\mathcal{M}(\mathbf{x},r_0) \approx \tilde\Gamma_{\mathbf{x}}\cap B_{r_0}(\mathbf{x})$.
Eventually, we seek $\Psi_{\Gamma,r_o}$ satisfying
\begin{equation}\label{eq:IBIM_reg_integral}
    %\frac{1}{|\mathcal{M}(\mathbf{x},r_0)|}
    \int_{\mathcal{M}(\mathbf{x},r_0)}K(\mathbf{x},\mathbf{y})\text{d}\sigma_{\mathbf{y}} = \int_{\mathcal{M}(\mathbf{x},r_0)}\Psi_{\Gamma,r_0}(\mathbf{x},\mathbf{y})\text{d}\sigma_{\mathbf{y}}+\mathcal{O}\left({r_0^p}\right), \qquad p\geq 2.
\end{equation}
If $v$ is a Lipschitz continuous function on $\Gamma$, we can write
\begin{align*}
    \nonumber
    \int_{\mathcal{M}(\mathbf{x},r_0)}K(\mathbf{x},\mathbf{y})v(\mathbf{y})\text{d}\sigma_{\mathbf{y}} &=
    v(\mathbf{x})\int_{\mathcal{M}(\mathbf{x},r_0)}K(\mathbf{x},\mathbf{y})\text{d}\sigma_{\mathbf{y}}+\int_{\mathcal{M}(\mathbf{x},r_0)}(v(\mathbf{y})-v(\mathbf{x}))K(\mathbf{x},\mathbf{y})\text{d}\sigma_{\mathbf{y}}  \\
    &= v(\mathbf{x})\int_{\mathcal{M}(\mathbf{x},r_0)}K(\mathbf{x},\mathbf{y})\text{d}\sigma_{\mathbf{y}}+\mathcal{O}(r_0^2)  \\
    &= v(\mathbf{x})\int_{\mathcal{M}(\mathbf{x},r_0)}\Psi_{\Gamma,r_0}(\mathbf{x},\mathbf{y})\text{d}\sigma_{\mathbf{y}}+\mathcal{O}(r_0^p)+\mathcal{O}(r_0^2)\,,
\end{align*} % |\mathcal{M}(\mathbf{x},r_0)
by using $|\mathcal{M}(\mathbf{x},r_0)|\sim r_0^2$ and $K(\mathbf{x},\mathbf{y})\sim\|\mathbf{x}-\mathbf{y}\|^{-1}$ for the kernels we are interested in.
This approach is used in \cite{kublik2013implicit} and \cite{chen2017implicit}. \\

The rest of this Section will be now dedicated to showing results of this approach for the Laplace and Helmholtz double-layer kernels.

\subsection{Application to the Laplace and Helmholtz double-layer kernels}

Consider the double-layer kernel $K(\mathbf{x},\mathbf{y})=\dfrac{\partial G_0}{\partial \mathbf{n}_y}(\mathbf{x},\mathbf{y})$ for Laplace. In~\cite{kublik2013implicit} the function $\Psi_{\Gamma, r_0}(\mathbf{x},\mathbf{y})$ is built as a constant 
function, $\Psi_{\Gamma, r_0}(\mathbf{x},\mathbf{y})\equiv C_{\Gamma,r_{0}}$,
\[
\int_{\mathcal{M}(\mathbf{x},{r_{0}})}\frac{\partial G_{0}}{\partial \mathbf{n}_{y}}(\mathbf{x},\mathbf{y})\text{d}\sigma_{\mathbf{y}}\approx \int_{\mathcal{M}(\mathbf{x},{r_{0}})}C_{\Gamma,r_{0}}\,\text{d}\sigma_{\mathbf{y}}\,.
\]
The constant $C_{\Gamma,r_{0}}$ represents the average of the integrand
on the set, and $\Psi_{\Gamma,r_0}(\mathbf{x},\mathbf{y})=C_{\Gamma,r_0}$ regularizes the double-layer kernel:
\begin{equation}\label{eq:KregIBIM-constant}
K_{r_0,C}^{reg}(\mathbf{x},\mathbf{y}):  =
    \begin{cases}
      \dfrac{\partial G_{0}}{\partial\mathbf{n}_{y}}(\mathbf{x},\mathbf{y}), & \|\mathbf{x}-\mathbf{y}\|\geq r_0, \\
      C_{\Gamma,r_0}, &
      \|\mathbf{x}-\mathbf{y}\|< r_0.
    \end{cases}
\end{equation}
The expression for $C_{\Gamma,r_0}$ found in this setting, dependent on the principal curvatures and $r_0$, is:
\begin{align*}\label{eq:IBIM_C_Cr0}
C_{\Gamma,r_{0}}:= & \frac{\kappa_{1}+\kappa_{2}}{8\pi r_{0}}-\frac{\kappa_{1}+\kappa_{2}}{512\pi}\left(13\kappa_{1}^{2}-2\kappa_{1}\kappa_{2}+13\kappa_{2}^{2}\right)r_{0}\\
 & +(\kappa_{1}+\kappa_{2})\left\{ \frac{(\kappa_{1}^{2}+\kappa_{2}^{2})\left(5\kappa_{1}^{2}-2\kappa_{1}\kappa_{2}+5\kappa_{2}^{2}\right)}{4096\pi}+\frac{\kappa_{1}^{4}+2\kappa_{1}^{2}\kappa_{2}^{2}+\kappa_{2}^{4}}{512\pi}\right\} r_{0}^{3}\nonumber \,.
\end{align*}
The calculations and details about the setting together with the exact definition of $\mathcal{M}(\mathbf{x},r_0)$
can be found in the~\ref{sub:AppendixRegularization}. 

When treating the Helmholtz double-layer kernel, we can apply this constant regularization approach to the additional term which differentiates it from the Laplace double-layer kernel: the gradient of the Helmholtz fundamental solution in three dimensions is 
\[
\nabla G_{\lambda}(\mathbf{z})=\frac{1}{4\pi}\frac{\exp(\text{i}\lambda\|\mathbf{z}\|)}{\|\mathbf{z}\|^{3}}\left(\text{i}\lambda\|\mathbf{z}\|-1\right)\mathbf{z}\,;
\]
hence, the double-layer kernel for Helmholtz takes the form:
\begin{align}
\frac{\partial G_{\lambda}}{\partial \mathbf{n}_{y}}(\mathbf{x},\mathbf{y}) & =\frac{1}{4\pi}\frac{(\mathbf{x}-\mathbf{y})^T \mathbf{n}_{y}}{\|\mathbf{x}-\mathbf{y}\|^{3}}\exp(\text{i}\lambda\|\mathbf{x}-\mathbf{y}\|)[1-\text{i}\lambda\|\mathbf{x}-\mathbf{y}\|]\nonumber \\
 & =\exp(\text{i}\lambda\|\mathbf{x}-\mathbf{y}\|)\left\{ \frac{\partial G_{0}}{\partial \mathbf{n}_{y}}(\mathbf{x},\mathbf{y})-\frac{\text{i}\lambda}{4\pi}\frac{(\mathbf{x}-\mathbf{y})^T \mathbf{n}_{y}}{\|\mathbf{x}-\mathbf{y}\|^{2}}\right\}. \label{eq:DL_Helmholtz}
\end{align}
In the expression above
the factor $\exp(\text{i}\lambda\|\mathbf{x}-\mathbf{y}\|)$
is a Lipschitz continuous function in $\mathbf{y}$, and we know how to deal numerically with the
Laplace double-layer kernel $\frac{\partial G_0}{\partial \mathbf{n}_y}$ (using regularizations \eqref{eq:KregIBIM-constant} or \eqref{eq:KregIBIM-linear}, or the corrected trapezoidal rule which will be the focus of Section~\ref{sec:corr_TR}), so we focus only on the secondary kernel
\begin{equation*}
\dfrac{(\mathbf{x}-\mathbf{y})^T \mathbf{n}_{y}}{\|\mathbf{x}-\mathbf{y}\|^{2}}\label{eq:Gtilde}
\end{equation*}
which, if $\kappa_1\neq\kappa_2$, is undetermined in $\mathbf{x}=\mathbf{y}$, as the value depends on the direction of approach. 
The maximum and minimum limit values are the ones found traveling
along the principal directions, equal to $\frac{\kappa_{1}}{2}$ and
$\frac{\kappa_{2}}{2}$ respectively. 

Then we wish to find a constant function
$\Psi_{\Gamma, r_0}(\mathbf{x},\mathbf{y})\equiv \tilde{C}_{\Gamma,r_{0}}$
such that the integral around the singularity point is approximated well,
\[
\int_{\mathcal{M}(\mathbf{x},{r_{0}})}\dfrac{(\mathbf{x}-\mathbf{y})^T \mathbf{n}_{y}}{\|\mathbf{x}-\mathbf{y}\|^{2}}\text{d}\sigma_{\mathbf{y}}=\int_{\mathcal{M}(\mathbf{x},{r_{0}})}\tilde{C}_{\Gamma,r_{0}}\,\text{d}\sigma_{\mathbf{y}}+\mathcal{O}(r_{0}^{4}).  
\]
This requirement gives %where
\begin{align}
\tilde{C}_{\Gamma,r_{0}}:= & \frac{\kappa_{1}+\kappa_{2}}{4}-\frac{\kappa_{1}+\kappa_{2}}{256}\left(-13\kappa_{1}^{2}+2\kappa_{1}\kappa_{2}-13\kappa_{2}^{2}\right)r_{0}^{2}\, .\label{eq:Gtilde_IBIM_Cr0}
\end{align}
It is interesting to notice that the first term in $\tilde{C}_{\Gamma,r_{0}}$
is the average of the maximum and minimum limits of the integrand. The function $\frac{(\mathbf{x}-\mathbf{y})\cdot \mathbf{n}_{y}}{\|\mathbf{x}-\mathbf{y}\|^{2}}$, for
$\mathbf{x},\mathbf{y}\in\Gamma$, can then be regularized using (\ref{eq:Gtilde_IBIM_Cr0}):
\begin{equation}\label{eq:GtildeIBIM-const}
\tilde{K}_{C,r_0}^{reg}(\mathbf{x},\mathbf{y}):=
    \begin{cases}
      \dfrac{(\mathbf{x}-\mathbf{y})^T \mathbf{n}_{y}}{\|\mathbf{x}-\mathbf{y}\|^{2}}\,, & \|\mathbf{x}-\mathbf{y}\|\geq r_0, \\[0.4cm]
      \tilde{C}_{\Gamma,r_{0}}, &
      \|\mathbf{x}-\mathbf{y}\|< r_0\,.
    \end{cases} 
\end{equation}

\begin{figure}
    \begin{center}
    \includegraphics{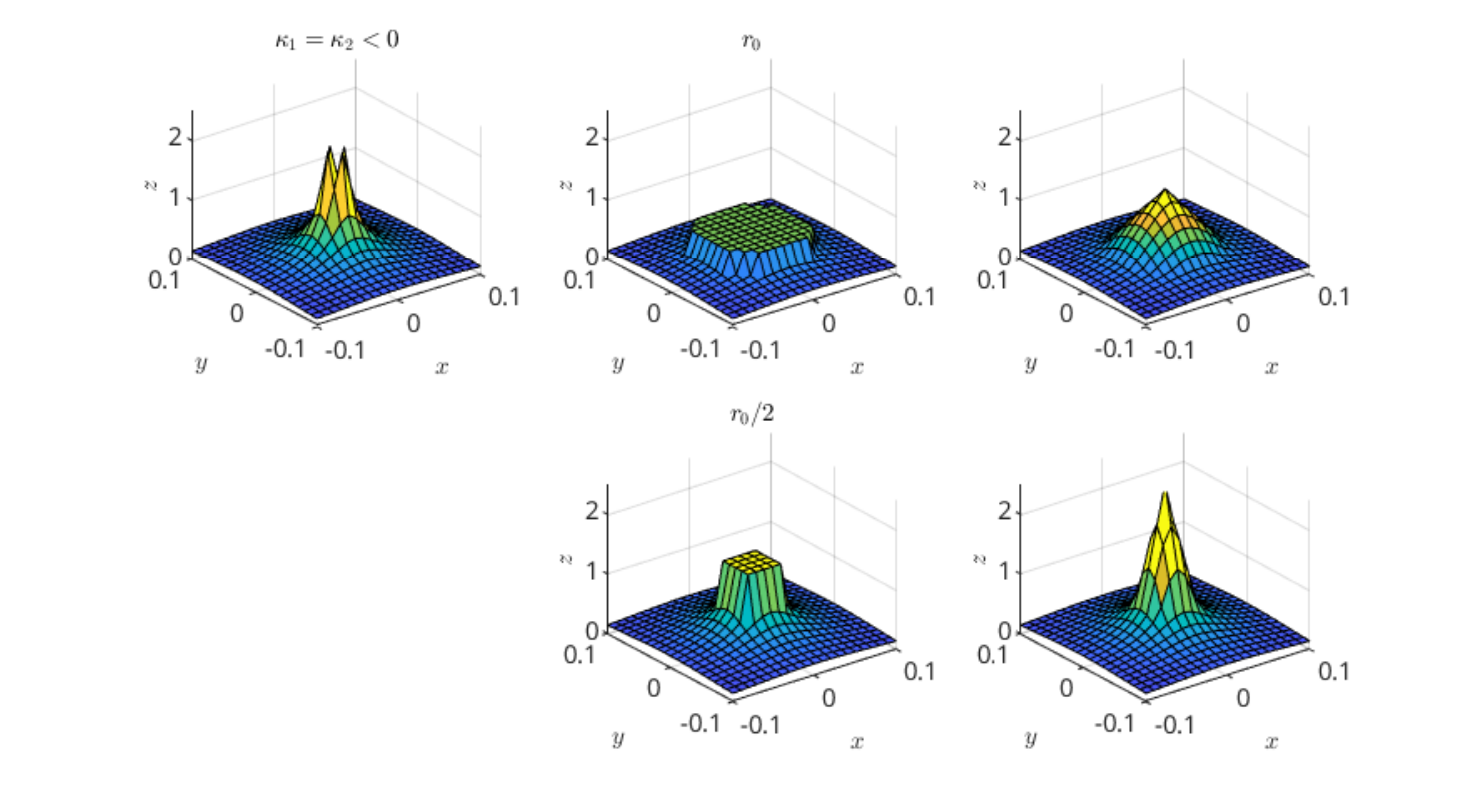}\\
    \includegraphics{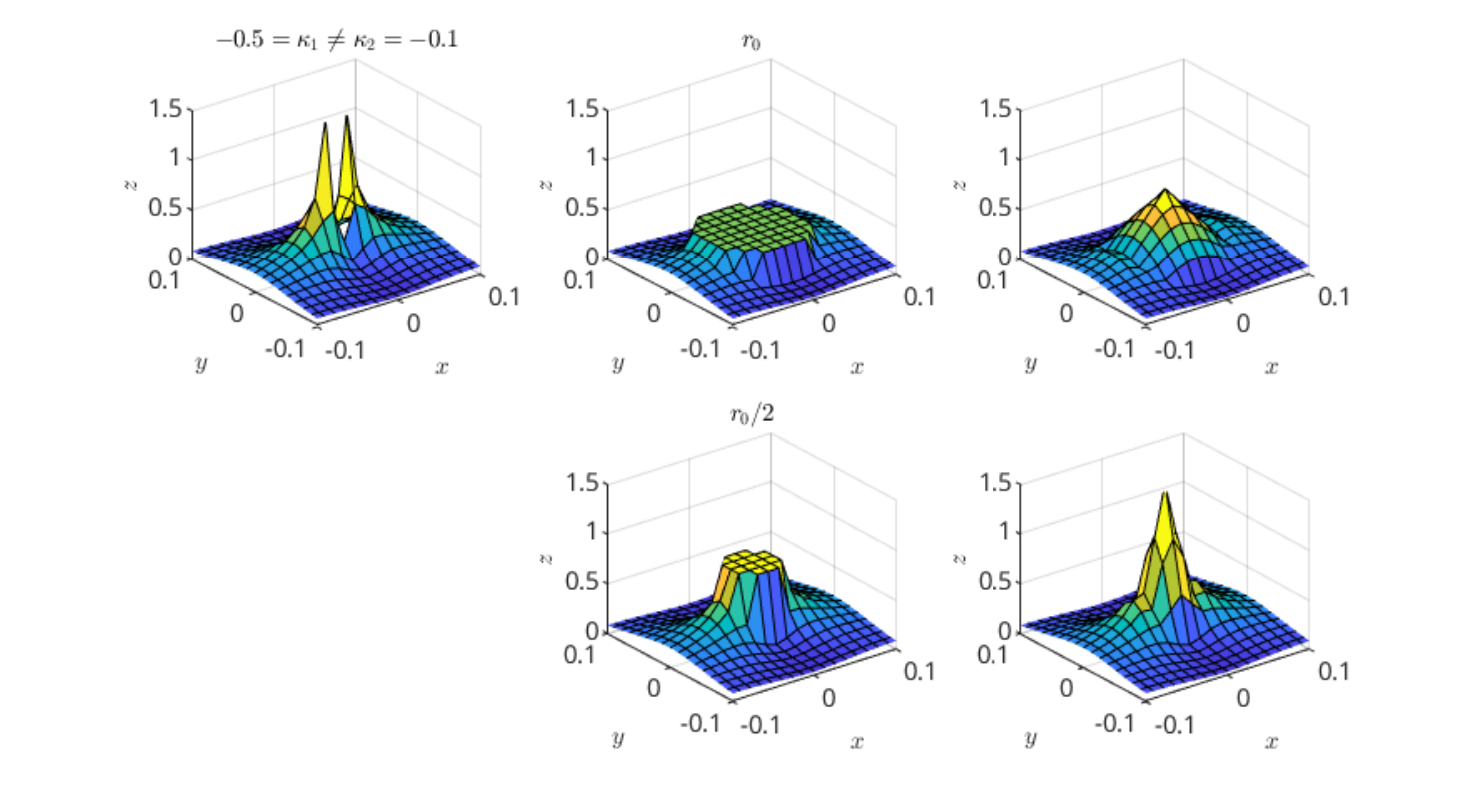}
    \end{center}
    \caption{Singular behavior - 1}
    \footnotesize{Curvatures with same sign: same (first and second row) and different  values (third and fourth row). Left column: double-layer kernel with target point in the origin, without any regularization. Center column: constant regularization. Right column: linear regularization. The second and fourth rows show the constant (center) and cappuccio (right) regularizations with $r_0$ halved compared to the first and third rows respectively.}
    \label{fig:regularizations1}
\end{figure}

\subsubsection*{New regularization with linear function (\emph{cappuccio})}
\label{sub:linear_regularization}
Here, we consider regularizing with a class of function that are linear with respect to the distance to the singularity. We construct 
the hood-like (\emph{cappuccio} in Italian) function
$\Psi^L_{\Gamma, r_0}(\mathbf{x},\mathbf{y}):=a_0\,\|\mathbf{x}-\mathbf{y}\|/r_0+a_{1}$:
\begin{figure}
    \begin{center}
    \includegraphics{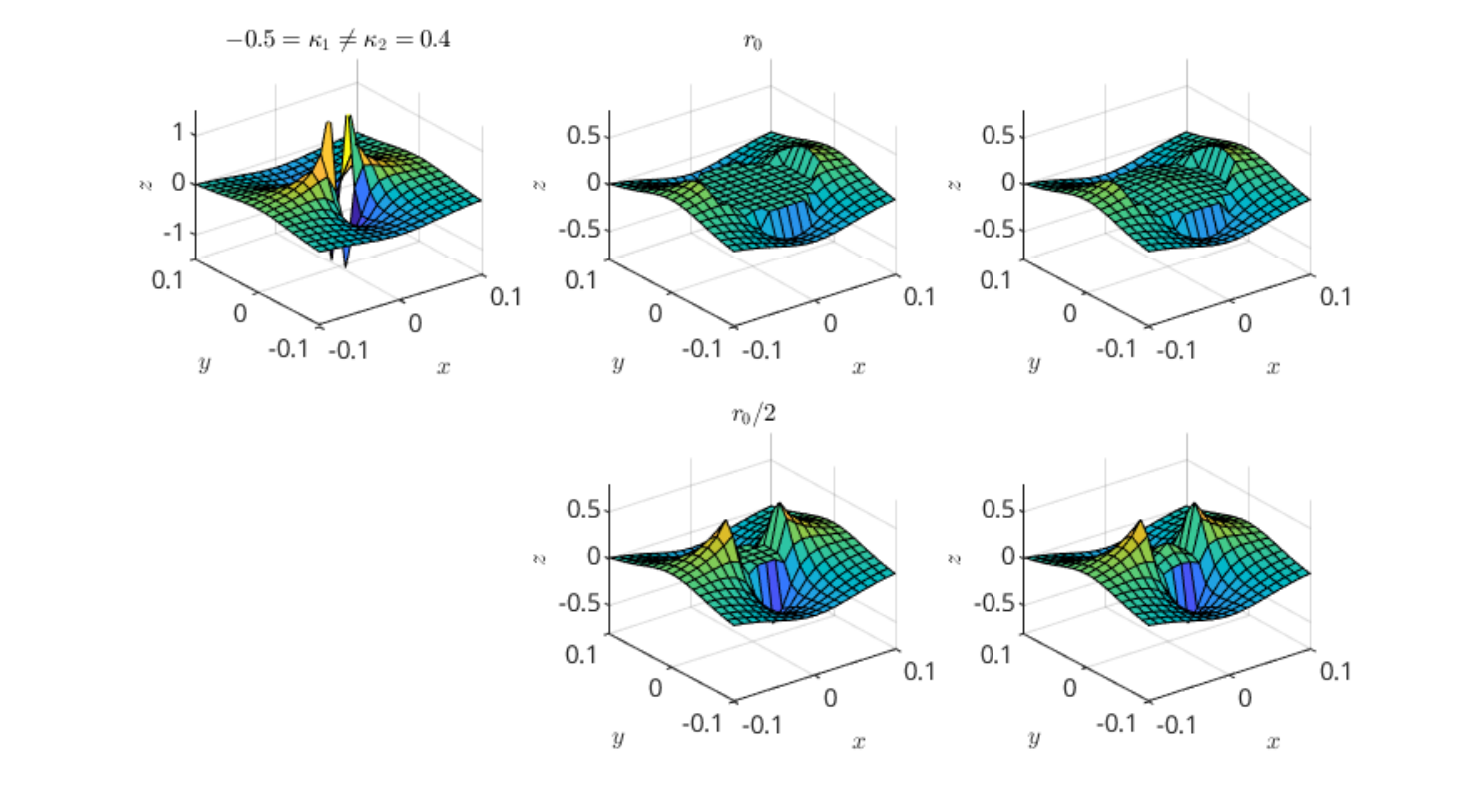} %-eps-converted-to
    \end{center}
    \caption{Singular behavior - 2}
    \footnotesize{Curvatures with different sign. Left column: double-layer kernel with target point in the origin, without any regularization. Center column: constant regularization. Right column: linear regularization. The second row shows the constant (center) and cappuccio (right) regularizations with $r_0$ halved compared to the first row. }
    \label{fig:regularizations2}
\end{figure}
\[
\int_{\mathcal{M}(\mathbf{x},r_0)}\frac{\partial G_{0}}{\partial \mathbf{n}_{y}}(\mathbf{x},\mathbf{y})\text{d}\sigma_{\mathbf{y}}\approx\int_{\mathcal{M}(\mathbf{x},r_0)}\Psi^L_{\Gamma,r_{0}}(\mathbf{x},\mathbf{y})\,\text{d}\sigma_{\mathbf{y}}\,.
\]
%\FI{
This condition imposes one constraint; the second constraint we impose is that 
\[
\Psi_{\Gamma,r_0}^L(\mathbf{x},\mathbf{y})\Big|_{\|\mathbf{x}-\mathbf{y}\|
=r_0}=\frac{1}{|\partial\mathcal{M}(\mathbf{x},r_0)|}\int_{\partial\mathcal{M}(\mathbf{x},r_0)}\frac{\partial G_{0}}{\partial \mathbf{n}_{y}}(\mathbf{x},\mathbf{y})\text{d}\sigma_{\mathbf{y}}\,,
\]
which means that $\Psi^L_{\Gamma,r_0}$ takes as outermost value the average of the kernel on the boundary of $\mathcal{M}(\mathbf{x},r_0)$.

The motivation is that the discontinuity in the regularized kernels can be significantly smaller, for small $r_0$, than the ones regularized by constants.  Consequently, the quadrature errors can be smaller. See Figures~\ref{fig:regularizations1}  and \ref{fig:regularizations2} for a comparison. 
%\CANC{In particular, Figure~\ref{fig:regularizations2} suggests ...}
{In particular, the third columns in Figures~\ref{fig:regularizations1} and~\ref{fig:regularizations2} show that the existing and proposed regularizations will lead to discontinuity between the regularized and original functions if the direction of approach to the target point is not taken into account. Consequently we believe} that a possible future work is to develop a {continuous} regularization {by including dependence} 
%\CANC{which depends}
on the principal directions $\bm{\tau}_1,\bm{\tau}_2$ in 
addition to the curvatures of the surface, e.g. $\Psi_{\Gamma,r_0}=\Psi_{\Gamma,r_0}(\mathbf{x}-\mathbf{y}, \kappa_1, \bm{\tau}_1, \kappa_2, \bm{\tau}_2)$.

We get the following expression for the regularized kernel,
\begin{equation}\label{eq:KregIBIM-linear}
    K_{r_0,L}^{reg}(\mathbf{x},\mathbf{y}):  =
    \begin{cases}
      \dfrac{\partial G_{0}}{\partial \mathbf{n}_{y}}(\mathbf{x},\mathbf{y}), & \|\mathbf{x}-\mathbf{y}\|\geq r_0\,, \\[0.4cm]
      a_0\,\dfrac{\|\mathbf{x}-\mathbf{y}\|}{r_0}+a_{1}, &
      \|\mathbf{x}-\mathbf{y}\|< r_0\,,
    \end{cases}
\end{equation}
where
\begin{align*}
a_0 = & -\frac{3(\kappa_1+\kappa_2)}{16\pi r_0}+\frac{3(\kappa_1+\kappa_2)}{5120\pi}(21\kappa_1^2-2\kappa_1\kappa_2+21\kappa_2^2)r_0\,, \\
a_1 = & \frac{\kappa_1+\kappa_2}{4\pi r_0}-\frac{3(\kappa_1+\kappa_2)}{2560\pi}(23\kappa_1^2-6\kappa_1\kappa_2+23\kappa_2^2)r_0\,.
\end{align*}
Note that $\Psi^L_{\Gamma,r_0}$ scales as
$1/r_0$ for small $r_0$.

In Section~\ref{sec:numerical_examples}, we shall present some numerical convergence studies of the approaches mentioned in this Section.

\section{The corrected trapezoidal rules}

\label{sec:corr_TR}
We have shown in Section \ref{subsec:ImplicitParametrization} that the singular integrals of interest can be characterized by their singular behavior along lines in $\mathbb{R}^3$.  We view the trapezoidal rule 
on a three-dimensional uniform Cartesian grid
as the sum over the trapezoidal rules applied to the two-dimensional uniform grids. On each two-dimensional grid, the case is reduced to correction of trapezoidal rule for functions that are singular only at a single point. However this point is typically not lying on any grid nodes. 

We will first present the trapezoidal rule and the existing methods for correcting it to achieve higher order convergence rates. We will then present our generalization of these works.

\subsection{The punctured trapezoidal rules}

Let $f:\mathbb{R}^n\rightarrow\mathbb{R}$ be a 
compactly supported smooth function. 
We are interested in approximating its integral $\int_{\R^n}f(\mathbf{x})\text{d}\mathbf{x}$ by utilizing values of $f$ on the uniform grid $h\mathbb{Z}^n$. By the compact support, 
the trapezoidal rule applied to $f$ becomes the following simple Riemann sum:
\begin{equation}
  T_{h}[f] := h^n\sum_{\mathbf{y}\in h\mathbb{Z}^n} f(\mathbf{y})\, . \label{eq:trapezoidalrule_nD}
\end{equation}
When $f$ is compactly supported
the order of accuracy of such summations depends on the regularity of $f$: if $f\in C^p$, the error is  $\mathcal{O}(h^p)$ (see \S 25.4.3 in \cite{abramowitzstegun} and \S 5.1 in \cite{isaacsonkeller}); the trapezoidal rule enjoys spectral accuracy if $f\in C^\infty$. 

When $f$ is continuous in $\R^n\setminus \{\mathbf{x}_0\}$ and singular at $\mathbf{x}_0$, where $\int_{\R^n}f(\mathbf{x})\text{d}\mathbf{x}$ exists as a Cauchy principal value, 
it is natural to modify the trapezoidal rule by skipping the summation over the grid nodes within certain distance to $\mathbf{x}_0$:
\begin{equation}
T_{h}^{0}[f]:=h^n\sum_{\mathbf{y}\in h\mathbb{Z}^n \setminus N_h(\mathbf{x}_0)}f(\mathbf{y}) \label{eq:puncturedtrapez}
\end{equation}
where $N_h(\mathbf{x}_0)$ determines which grid nodes we remove. We will call \eqref{eq:puncturedtrapez} the \textit{punctured trapezoidal rule} when $N_h$ includes only a single grid node; in other words,  
\begin{equation}
\label{eq:puncturedTRcutoffregion}
    N_h(\mathbf{x})=\{\mathbf{y}\in\R^n : \|\mathbf{x}-\mathbf{y}\|_\infty\leq h/2\}\,.
\end{equation} 
The punctured trapezoidal rule converges, but with lower order rates at best, even though $f$ may be $C^\infty$ in the punctured domains.  For example, in one dimension for $f(x)=\log|x|$, the order of convergence is sublinear $\mathcal{O}(h^p)$, $p<1$, and in two dimensions for $f(\mathbf{x})=\|\mathbf{x}\|^{-1}$ the order is 1. The large decrease in order is exactly the property we would like to address with the correction technique. 

The idea is to add a correction term to \eqref{eq:puncturedtrapez}, which 
makes up for the integral over $N_h$.
In the following, we describe an approach for defining such corrections in detail. 

\subsection{Corrections for the trapezoidal rule}
\label{sub:corr_1D}

From this point forward, we will assume the function $f$ can be factored into the following form
\begin{equation}\label{eq:fsv}
    f(\mathbf{x})=s(\mathbf{x}-\mathbf{x}_0)v(\mathbf{x}),~~~\mathbf{x}\in\mathbb{R}^n\setminus\{\mathbf{x}_0\}
    %:\mathbb{R}^n\rightarrow\mathbb{R}
\end{equation}
where $s$ represents an integrable function, singular in
the origin, and $v$ represents a smooth compactly supported function in $\mathbb{R}^n$.
In this Section we discuss a general approach to developing high order quadratures for the integration of such type of functions. In Section~\ref{sec:corr_laplace_kernels} we will provide specific choices of $s$ for use with single- and double-layer kernels arising from the Laplace or Helmholtz operator.

The trapezoidal rule is a sum of
the function values on the grid,
where all values have the same
weight $h^n$.
Improving the order of accuracy of the trapezoidal rule by modifying the weights
close to the singularity point has been an approach studied and applied successfully 
with different kinds of singular behaviors and in different dimensions. See for example \cite{kapur1997high, marin2014corrected, wumartinsson2020zeta}. 

The following is a brief presentation of the one- and two-dimensional corrections found in \cite{marin2014corrected}, where $\mathbf{x}_0$ is always assumed to be the origin. 

The starting point is the punctured trapezoidal rule in one dimension.
When $s(x)=|x|^\gamma$ for $-1<\gamma<0$,
an error expansion of the following type can be derived,
$$
\int s(x)v(x)\text{d}x=T_{h}^{0}[s\, v]+h^{1+\gamma} \omega \,v(0)  + O(h^{3+\gamma}).
$$
The goal is to find the constant $\omega$, which is \emph{independent
of $v$} (but depends on $\gamma$), and use it to correct the rule as
$$
  Q_h^0[f]:=T_{h}^{0}[f]+h^{1+\gamma}\omega v(0).
$$
While $T_{h}^{0}$ is of order ${1+\gamma}$, 
the method $Q_h^0$ is of order ${3+\gamma}$.
Note that $Q_h^0$ only modifies the original trapezoidal rule in
one point; the value in the singular point is replaced by
the value of the smooth part $v(0)$, weighted by $\omega$ and a 
suitable power of $h$. 
In general $\omega$ is a
functional of the singular function $s$
and we write $\omega=\omega[s]$.

In order to find $\omega[s]$ we define $\bar{\omega}(h)$
as the actual error of $T_{h}^{0}$ for
a fixed $h$ and a smooth test function $g$ with $g(0)=1$,
scaled by $h^{1+\gamma}$. More precisely,
\begin{equation}
\int s(x)g(x)\text{d}x=T_{h}^{0}[s\, g]+h^{1+\gamma}\bar{\omega}(h).
\label{eq:whdef}
\end{equation}
From the error expansion above, since $g(0)=1$, we see that
$$
\bar{\omega}(h) = \frac1{h^{1+\gamma}}\left(\int s(x)g(x)\text{d}x-T_{h}^{0}[s\, g]
\right)
 = \omega[s] + O(h^2).
$$
Hence, the weight $\bar{\omega}(h)$ converges to $\omega[s]\neq 0$ for $h\to0^+$, independent of $g$. 

This is a crucial property, which makes it possible 
to compute, store and reuse $\omega$ for the integration of any 
{integrand} of the kind \eqref{eq:fsv}. In order to do that
one needs to be able to accurately compute the integral $\int s(x)g(x)dx$
containing the test function. As this 
computation is only needed for
one function $g$, it can be done either by analytical means
or adaptive high order numerical integration.
If the test function $g$ is chosen more flat at the singularity
point, such that $0=g'(0)=g''(0)=\dots$ one can show that
the convergence
will be faster, which makes the numerical computations easier.
In fact, $g$ would be ideally the constant function $g({x})=1$  but in order to avoid dealing with the boundary conditions of trapezoidal rule and keep the expression of $T_h^0$ equal to the Riemann sum with the exclusion of a single node, $g$ is taken compactly supported.

Higher order corrections are also possible,
where more terms in the error expansion are cancelled. 
The weights
must then be modified in more points close to the singularity.
The condition (\ref{eq:whdef})
can be interpreted as requiring that $T^0_h$, corrected with
the weight $\bar\omega(h)$,
integrate $s(x) g(x)$ exactly. When multiple weights are used,
the weights are similarly defined by requiring that 
the modified method integrates not only $s(x) g(x)$ exactly
but also $s(x) g(x) x$, $s(x) g(x) x^2$, \ldots. A set of $h$-dependent
weights are then obtained, which converge as $h\to 0^+$.

One can also apply the same idea to other singularities.
In \cite{wumartinsson2020zeta} this was done for
$s(x)=-\log|x|$. Then the factor $h^{1+\gamma}$
must be replaced by an expression
$a(h)=h(2-\log h)$ and a second order method is obtained
\[
  \int s(x)v(x)\text{d}x=Q_h^0[s\, v] +\mathcal{O}(h^{2}),
  \qquad Q_h^0[s\, v]=T_h^0[s\, v]+ a(h)\omega[s]\, v(0).
\]
In two dimensions, similar to one dimension, for functions \eqref{eq:fsv} with $s(\mathbf{x})=\|\mathbf{x}\|^{-1}$, the corrected trapezoidal rule is defined as
\begin{equation*}
  Q^0_h[s\, v] := T_h^0[s\, v] + h\,\omega[s]\, v(\mathbf{\mathbf{0}})
\end{equation*}
where $\omega[s]$ is calculated as the limit:% \displaystyle
\begin{equation}
  \omega[s] := \lim_{\delta\to0^+} 
  \frac{1}{\delta}\left( \int_{\mathbb{R}^2}s(\mathbf{x})g(\mathbf{x})\text{d}\mathbf{x}-T_\delta^0[s\, g]\right),
\end{equation}
for a test function $g$ with $g(\mathbf{0})=1$.
In \cite{marin2014corrected} it was proven that 
the corrected method 
for $s(\mathbf{x})=\|\mathbf{x}\|^{-1}$ is 
third order accurate,
\begin{equation}
\label{thirdordererror}
  \int_\R s(\mathbf{x})v(\mathbf{x})\text{d}\mathbf{x}=Q_h^0[s\, v] +\mathcal{O}(h^{3})\,.
\end{equation}

\subsubsection*{Corrections for singularity unaligned to the grid}

\begin{figure}
    \begin{center}
        \includegraphics{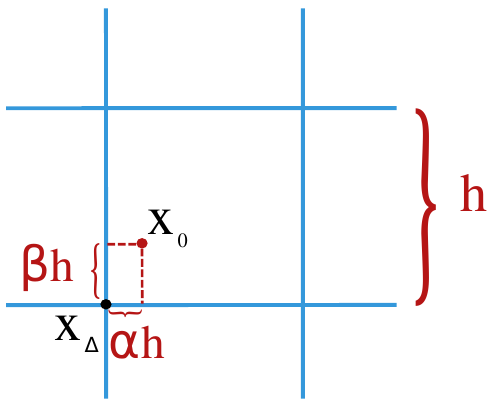}
    \end{center}
    \caption{Singularity unaligned to the grid}
    \footnotesize{Position of the singularity point $\mathbf{x}_0$ relative to the closest grid node $\mathbf{x}_\Delta$; the parameters $\alpha,\beta$ are used to characterize its position relative to the grid.}
    \label{fig:single_correction_grid_plot}
\end{figure}

To prepare for the proposed quadrature rules for implicit boundary integrals, we first generalize the approach presented in Section~\ref{sub:corr_1D} to the case when
the singularity does not lie on a grid node. 
We consider two dimensions and retain the assumption that $f$ 
can be factorized as $s(\mathbf{x}-\mathbf{x}_0)v(\mathbf{x})$ where $s$
has a singularity 
and $v$ is smooth and compactly supported. 
However, the singularity is now in a
point $\mathbf{x}_0$ which may not be part of the grid.
We let $\mathbf{x}_\Delta$
be the grid node closest to $\mathbf{x}_0$,
$$
 \mathbf{x}_\Delta = \arg\min_{\mathbf{x}\in h\mathbb{Z}^2}
 \left\|\mathbf{x}-\mathbf{x}_0\right\|,
$$
or one of the closest in the case that it may not be unique,
such that
$$
  (\alpha,\beta)=\frac{\mathbf{x}_0-\mathbf{x}_\Delta}{h},~~~\text{for some $\alpha,\beta\in[-1/2,1/2$)\,,}
$$
as shown in Figure~\ref{fig:single_correction_grid_plot}.
When $(\alpha,\beta)\neq \mathbf{0}$ the usual trapezoidal rule
is well-defined also for
the singular function and the same type of error expansion
holds as for the punctured trapezoidal rule in the previous
Section. However, the error constant is not uniform
and blows up as $(\alpha,\beta)\to \mathbf{0}$. We therefore use
the punctured trapezoidal rule also for unaligned grids
as the base method for correction.

The singular functions considered
in this paper are of the form
$f(\mathbf{x})=s(\mathbf{x}-\mathbf{x}_0)v(\mathbf{x})$
where $|s(\mathbf{x})|\sim ||\mathbf{x}||^{-1}$.
For those functions the same
scaling in $h$
as $s(\mathbf{x})= ||\mathbf{x}||^{-1}$
is appropriate and
we define the \emph{single-correction trapezoidal rule} for 
unaligned grids in two dimensions as
\begin{equation}
  \bar{Q}^{2D}_h[f] := 
  T_h^{0}[s(\,\cdot -\mathbf{x}_0) \,v(\,\cdot\,)] 
  + h\,\omega[s;\alpha,\beta]v(\mathbf{x}_\Delta).
  \label{eq:single_correction_quadrature}
\end{equation}
The weight is given as the limit of the sequence:% \displaystyle
\begin{equation}
  \omega[s;\alpha,\beta] := \lim_{\delta\to0^+} \omega_\delta[s;\alpha,\beta].
  \label{eq:single_correction_weight}
\end{equation}
where as before
$\omega_\delta[s;\alpha,\beta]$
is defined using a smooth compactly supported test function $g$ with $g(\mathbf{0})=1$,
\begin{align}
  \omega_\delta[s;\alpha,\beta] :=& 
  \frac{ \int_{\mathbb{R}^2}{s}(\mathbf{x}-\mathbf{x}_0){g}(\mathbf{x}-\mathbf{x}_0)\text{d}\mathbf{x}-T_\delta^{0}\Big[s(\,\cdot -\mathbf{x}_0)\, g(\,\cdot -\mathbf{x}_0)\Big]}{\delta\, g(\mathbf{x}_\Delta-\mathbf{x}_0)}\nonumber\\
  =& \frac{ \int_{\mathbb{R}^2}{s}(\mathbf{x}){g}(\mathbf{x})\text{d}\mathbf{x}-T_\delta^{0}\Big[s(\,\cdot -(\alpha,\beta)\delta)\, g(\,\cdot -(\alpha,\beta)\delta)\Big]}{\delta\, g(-(\alpha,\beta)\delta)}.
  \label{eq:single_correction_weight-before-taking-limit}
\end{align}
In the last step we shifted
the exact integral by ${\bf x}_0$
and the trapezoidal rule by ${\bf x}_\Delta$, the closest node of the grid $\delta\mathbb{Z}^2$,
to show that the weight,
in addition to $s$,
only depends on the 
difference ${\bf x}_\Delta-{\bf x}_0$,
i.e. on $\alpha$ and $\beta$,
not on ${\bf x}_0$ itself.

The expression converges quickly when the stepsize $\delta$ is halved, and the accurate computation of $\omega$ is possible without 
needing specialized quadratures for singular integrands. In this paper, $\omega[s;\alpha,\beta]$ will be computed offline and tabulated for a suitable set of $(\alpha,\beta)$, and for relevant functions $s$; for values outside of the tabulation, we will use interpolation. 
In the next Section, we will discuss a few specific cases involving layer kernels, and we shall then present more details about the approximation of $\omega$ via tabulation and interpolation.

When $\alpha=\beta=0$ we get the same weights as in the aligned case. In particular, for $s(\mathbf{x})=\|\mathbf{x}\|^{-1}$ the limit
\eqref{eq:single_correction_weight} will find the same weight as the one found in \cite{marin2014corrected}. For the more general kernels considered in this paper and with unaligned grids{, in numerical experiments} we observe an error expansion 
of the type
\begin{equation}\label{eq:error_2D}
\int_{\R^2} s(\mathbf{x}-\mathbf{x}_0)v(\mathbf{x})\text{d}\mathbf{x}=
\bar{Q}_h^{2D}[s(\,\cdot -\mathbf{x}_0) \,v(\,\cdot\,)]  + F_1(\alpha,\beta)h^2 + F_2(\alpha,\beta)h^3+\mathcal{O}(h^4),
\end{equation}
where $F_1$ is a smooth function of $\alpha$, $\beta$.
Moreover, $F_1(0,0)=0\neq F_2(0,0)$,
which means that we have the same third order error \eqref{thirdordererror}
as for $s(\mathbf{x})=\|\mathbf{x}\|^{-1}$
when the grid is aligned. 
The properties of $F_1$ will be further explored in Section~\ref{sec:error-cancellation}. Proving these error results rigorously is in program for future research.

\subsection{Corrected trapezoidal rules for implicit boundary integrals}
\label{sub:corr_3D}

We describe our approach in developing corrected trapezoidal rules for the family of integrals defined in \eqref{eq:volume_potential_SL}:
\begin{equation}\label{eq:model-IBIM-integral}% volume_potential_SL
    \int_{\mathbb{R}^3} \overline{K}(\mathbf{x}^*,\mathbf{y})\rho(\mathbf{y}) \delta_{\Gamma,\varepsilon}(\mathbf{y})\text{d}\mathbf{y},~~~\mathbf{x}^*\in\Gamma.
\end{equation}
Without loss of generality, we consider, as the target point, $\mathbf{x}^*=(x^*,y^*,z^*)\in\Gamma$, where
the surface normal at $\mathbf{x}^*$ is $\mathbf{n}=(n_1,n_2, 1)$: from this point forward we will only consider this case, and if the normal direction points instead more towards the $\vec{x}$ or $\vec{y}$ directions, we can apply a change of coordinates and proceed with the same reasoning.
To simplify notation we let $f$ be the integrand
\begin{equation}\label{eq:kernelfunction}
f(\mathbf{y}):=\overline{K}(\mathbf{x}^*, \mathbf{y})\rho(\mathbf{y})\delta_{\Gamma,\varepsilon}(\mathbf{y})\,.
\end{equation}
Notice that $f$ in \eqref{eq:kernelfunction} is compactly supported in $T_\varepsilon$, if $\Gamma$ is a compact point set.
Furthermore, since the restricted kernel $\overline{K}(\mathbf{x}^*, \mathbf{y})=\overline{K}(\mathbf{x}^* + t\,\mathbf{n}, \mathbf{y})$ for all $t\in\mathbb{R}$,
the integrand $f$ is singular along this line.
 
\begin{figure}
    \begin{center}
        \includegraphics[height=4.5cm]{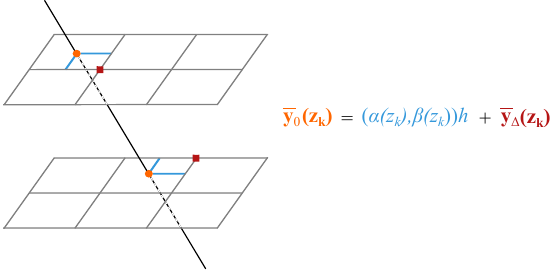}
    \end{center}
    \caption{Intersections of the line in three dimensions}
    \footnotesize{
    Intersection of the line passing through $\mathbf{x}^*$ with direction $\mathbf{n}$ with the planes $\{z=z_k\}$: on every plane, the intersection will be  $\mathbf{y}_{0}(z_k)$ (orange circle), and the closest grid node will be $\mathbf{y}_{\Delta}(z_k)$ (red square). The parameters which characterize the position of $\mathbf{y}_{0}(z_k)$ with respect to the grid $h\mathbb{Z}^2$ are $(\alpha(z_k),\beta(z_k))$ such that $\mathbf{\bar y}_{0}(z_k)=(\alpha(z_k),\beta(z_k))h+\mathbf{\bar y}_{\Delta}(z_k)$.}
    \label{fig:visualization_normal_surface_grid}
\end{figure}

The plan is then to construct a quadrature for \eqref{eq:model-IBIM-integral} ``plane-by-plane'' on the grid $h\mathbb{Z}^3$. First the standard trapezoidal rule
is used in the
$z$-direction,
\begin{align*}
\int_{\mathbb{R}^{3}}f(x,y,z)\text{d}x\text{d}y\text{d}z &= \int_{\mathbb{R}} \left\{ \int_{\mathbb{R}^2}f(x,y,z)\text{d}x\text{d}y\right\}\text{d}z 
\approx 
h\sum_{k}\int_{\mathbb{R}^{2}}f(x,y,kh)\text{d}x\text{d}y. \label{eq:plane-by-plane-sum}
\end{align*}
Then, the corrected trapezoidal rule is used to compute
the integrals on each plane,
$$
\int_{\mathbb{R}^2}f(x,y,kh)\text{d}x\text{d}y
\approx \bar{Q}^{2D}_h[f(\,\cdot\,,\,\cdot\,,kh)].
$$
See Figure~\ref{fig:visualization_normal_surface_grid} for an illustration. 
As the $z$-component of $\mathbf{n}$ is 1, $f$
is singular in one point only when restricted to the planes.
We can therefore use the quadrature rules described above. 

We let 
$\mathbf{\bar y}=(x,y)$ denote a point in the $xy$-plane
and introduce the projection onto this plane
\begin{equation}\label{eq:z-projection}
    \mathbf{\bar y}=\pi\mathbf{y}=\pi(x,y,z)=(x,y)\,.
\end{equation}
For a fixed $z$, the singular point 
$\mathbf{\bar y}_0$ of $f(\,\cdot\,,\,\cdot\,,z)$
is then  given by
\begin{equation}\label{eq:y0(z)_expression}
\mathbf{\bar y}_0(z) =\pi
\mathbf{y}_0(z), \qquad \mathbf{y}_0(z)=\mathbf{x}^*+(z-{z}^*)\mathbf{n}
=(\mathbf{\bar y}_0(z),z)\, .
\end{equation}

The corresponding closest grid node is denoted by $\mathbf{\bar y}_\Delta(z)$
and
its shift parameters
$\alpha=\alpha(z)$,
$\beta=\beta(z)$.
The factorization of $f(\,\cdot\,,\,\cdot\,,z)$
will be of the kind
\begin{equation}\label{eq:f_s_v_3D_splitting}
f(\mathbf{\bar y},z)=
s(\mathbf{\bar y}-\mathbf{\bar y}_0(z),z) v(\mathbf{\bar y},z)\, ,
\end{equation}
where $s$ is smooth 
in the second argument,
which ensures that the
partial integral $\int f\text{d}x\text{d}y$
is smooth in $z$,
justifying the use of the
standard trapezoidal rule
in this variable.
The functions $s$ and $v$ correspond to factorizations specific to the kernel $K$ and the geometry of $\Gamma$.  They will be discussed in detail in the next Section. 

With this notation we can now give the precise form of the
corrected method
$$
\int_{\mathbb{R}^2}f(x,y,z)\text{d}x\text{d}y
\approx {T}^{0}_h[
s(\,\cdot-\mathbf{\bar y}_0(z),z)
v(\,\cdot\,,z)]
+h\omega[s(\,\cdot\,,z);\alpha(z),\beta(z)]v(\mathbf{\bar y}_\Delta(z),z).
$$
After the discretization, $z_k=kh$, and
we can write the full method as
\begin{comment}
    \begin{align}
        Q_h[f] &=h\sum_{k\in \mathbb{Z}} \left\{
        {T}^{0}_h[
        s(\,\cdot-\mathbf{\bar y}_{0,k},z_k)
        v(\pi\cdot\,,z_k)]
        +h\,\omega[s(\,\cdot,z_k);\alpha_k,\beta_k]v(\mathbf{\bar y}_{\Delta,k},z_k)\right\}
        \nonumber \\
        &=h^3\sum_{k\in \mathbb{Z}}\ \sum_{\mathbf{\bar y}\in (
        h\mathbb{Z}^2\setminus N_h(\mathbf{y}_{0,k}))}
        s(\mathbf{\bar y}-\mathbf{\bar y}_{0,k},z_k)
        v(\mathbf{\bar y},z_k)
        +h^2\sum_{k\in \mathbb{Z}}\omega[s(\,\cdot,z_k);\alpha_k,\beta_k]v(\mathbf{{y}}_{\Delta,k})
        \nonumber \\
        &=h^3\sum_{\mathbf{x}\in
        (h\mathbb{Z}^3\setminus \mathcal{N}_h(\mathbf{x}^*))}
        f(\mathbf{x})
        +h^2\sum_{k\in \mathbb{Z}}\omega[s(\,\cdot,z_k);\alpha_k,\beta_k]v({\mathbf{{y}}}_{\Delta,k}),\label{eq:Q1_3D}
    \end{align}
\end{comment}

\begin{align}
Q_h[f] =& h\sum_{k\in \mathbb{Z}} \left\{
{T}^{0}_h\left[
s(\,\cdot-\mathbf{\bar y}_{0}(z_k),z_k)
v(\pi\,\cdot\,,z_k)\right]
+h\,\omega\left[s(\,\cdot\,,z_k);\alpha(z_k),\beta(z_k)\right]v(\mathbf{\bar y}_{\Delta}(z_k),z_k)\right\}
\nonumber \\
=&h^3\sum_{k\in \mathbb{Z}}\ \sum_{\mathbf{\bar y}\in (
h\mathbb{Z}^2\setminus N_h(\mathbf{y}_{0}(z_k)))}
s(\mathbf{\bar y}-\mathbf{\bar y}_{0}(z_k),z_k)
v(\mathbf{\bar y},z_k)
\nonumber\\&+h^2\sum_{k\in \mathbb{Z}}\omega[s(\,\cdot\,,z_k);\alpha(z_k),\beta(z_k)]v(\mathbf{{y}}_{\Delta}(z_k))
\nonumber \\
=& h^3\sum_{\mathbf{x}\in
(h\mathbb{Z}^3\setminus \mathcal{N}_h(\mathbf{x}^*))}
f(\mathbf{x})
+h^2\sum_{k\in \mathbb{Z}}\omega[s(\,\cdot\,,z_k);\alpha(z_k),\beta(z_k)]v({\mathbf{{y}}}_{\Delta}(z_k)),\label{eq:Q1_3D}
\end{align}
where
\begin{equation}\label{eq:puncturedTRcutoffregion3D}
    \mathcal{N}_h(\mathbf{x}^*):= \bigcup_{k} N_h(\mathbf{y}_{0}(z_k)) = \bigcup_{k} \mathbf{y}_{\Delta}(z_k)\,.
\end{equation}
Then we can see that this method is of the form \eqref{eq:corrected_trapezoidal_general_form} with $\mathcal{N}_h(\mathbf{x})$ as in \eqref{eq:puncturedTRcutoffregion3D} and
\begin{equation}\label{correction_term_correctedTR}
    \mathcal{R}_h(\mathbf{x})=h^2\sum_{k\in\mathbb{Z}} \omega\left[s(\,\cdot\,,z_k);\alpha(z_k),\beta(z_k)\right]v(\mathbf{y}_{\Delta}(z_k))\,.
\end{equation}

\section{Factorization of the Laplace kernels and the resulting quadratures} 

\label{sec:corr_laplace_kernels}

In the previous Section \ref{sub:corr_3D} we have seen the corrected trapezoidal rule \eqref{eq:Q1_3D} for the family of integrals of the kind \eqref{eq:model-IBIM-integral}. We have however not gone into detail about the form of the singular functions and the consequent splitting \eqref{eq:f_s_v_3D_splitting} involved, as they depend on the functions \eqref{eq:kernelfunction} and the corresponding kernels 
\[
\overline{K}(\mathbf{x},\mathbf{y})=K(\mathbf{x},P_\Gamma \mathbf{y})
\]
where $K$ is one of the 
{Laplace layer kernels}: 
\begin{equation}
\begin{array}{rlrr}
\text{(SL)}:\ \ G_{0}(\mathbf{x},\mathbf{y}) =& \dfrac{1}{4\pi}\dfrac{1}{\|\mathbf{x}-\mathbf{y}\|} \ ,\\[0.4cm] \text{(DL)}:\ \ \dfrac{\partial G_{0}}{\partial \mathbf{n}_{y}}(\mathbf{x},\mathbf{y})= & \dfrac{1}{4\pi}\dfrac{(\mathbf{x}-\mathbf{y})^T \mathbf{n}_{y}}{\|\mathbf{x}-\mathbf{y}\|^{3}}\,,\\[0.4cm]
\text{(DLC)}: \ \ \dfrac{\partial G_{0}}{\partial \mathbf{n}_{x}}(\mathbf{x},\mathbf{y})= & -\dfrac{1}{4\pi}\dfrac{(\mathbf{x}-\mathbf{y})^T \mathbf{n}_{x}}{\|\mathbf{x}-\mathbf{y}\|^{3}}\,.
\end{array}\label{eq:IBIM_kernels}
\end{equation}
In this section, we will derive the proposed quadrature rules for the above kernels.
In the same way the surface integral of the single-layer potential 
\eqref{eq:SL} in the non-parametric setting becomes the volume 
integral \eqref{eq:volume_potential_SL}, the double-layer 
\eqref{eq:DL} and double-layer conjugate \eqref{eq:DLC} potentials are 
extended to volume integrals in the tubular neighborhood $T_\varepsilon\subset\R^3$. 
Between the DL and DLC kernels we will only consider the DLC kernel, because the singularity behavior is identical.\\

{ In Section~\ref{sub:corr_allkernels} we present our approach to define the factorization $f=s~v$ in \eqref{eq:f_s_v_3D_splitting}. $s$ describes $K$ close to the singularity point and it will be written as the product of $S~\ell$, where $S$ takes a very simple form that is easy to work with.  
After introducing the function $S$, we will explain how to find the corresponding function $\ell$ in Section~\ref{sub:section511-finding-ell}, and then show the full expression of \eqref{eq:Q1_3D} in Section~\ref{sub:section512-factorization-splitting}.
Finally, in Section \ref{sub:section52-weights}, we will describe how the weights $\omega[s;\alpha,\beta]$ are defined from $S\,\ell$ and computed. A brief outline of the process can be found in Table \ref{tab:section5table}. }

\begin{table}
    \caption{Key ingredients in Section 5.}
    \label{tab:section5table}
    \begin{align*}
        \hline\\
         \int_{\R^3}&\underbrace{\overline{K}(\xx^*,\y){\rho(\y)\delta_{\Gamma,\varepsilon}(\y)}}_{=f(\y)}\text{d}\y\,\approx \,  Q_h[f]\, =\, \underbrace{h^3\sum_{\y\in h\mathbb{Z}^3\setminus \mathcal{N}_h(\xx^*)}f(\y)}_{=T_h^0[f]}+\underbrace{h^2\sum_{k\in\mathbb{Z}}\omega_k \, v(\y_{\Delta}(kh))}_{\mathcal{R}_h(\xx^*)}  \nonumber\\[0.4cm]
        \hline\\
        \text{In}&\text{{ Section \ref{sub:corr_allkernels}: }} \\
        &\left\{ 
         \begin{array}{l}
            \begin{array}{l}
                f(\y)=f(\bar\y,z) = \underbrace{S(\y - \y_0(z),\mathbf{n})\,\ell\left(\dfrac{\bar\y - \bar\y_0(z)}{\|\bar\y - \bar\y_0(z)\|};z\right)}_{=s(\bar\y-\bar\y_0(z),z)\,\text{ from \eqref{eq:f_s_v_3D_splitting}}}\,v(\y)\,, \\[0.5cm]
                \text{where} \\[0.2cm]
                S(\mathbf{r},\mathbf{n}) := \|\mathbf{n}\|\left\|\mathbf{r}\times\mathbf{n}\right\|^{-1}\ ,\ \ \mathbf{r},\mathbf{n}\in\R^3\,,
            \end{array}
            \\[0.1cm]
            \left.
                \ell(\mathbf{q};z) := \lim_{t\to 0^+} \dfrac{\overline{K}(\xx^*,t(q_1,q_2,0)+\y_0(z))}{S(t(q_1,q_2,0),\mathbf{n})}\ ,\ \ \mathbf{q}\in\mathbb{S}^1\,,
            \right\} \ \text{\footnotesize{ in Section \ref{sub:section511-finding-ell}.}}\\[0.6cm] 
            \left.\begin{array}{l}
                \text{Additive splitting:}\\
                \text{\footnotesize{ (in Section \ref{sub:section512-factorization-splitting})}}\\
                f(\y) = S(\y - \y_0(z),\mathbf{n})\, \hat v (\y) + S(\y - \y_0(z),\mathbf{n})\,\ell\left(\dfrac{\bar\y - \bar\y_0(z)}{\|\bar\y - \bar\y_0(z)\|};z\right) V(\y)\,, \\[0.3cm] % s(\bar\y-\bar\y_0(z),z)
                \text{where}\\[0.2cm]
                V(\y):=\rho(\y)\,\delta_{\Gamma,\varepsilon}(\y)\,, \\[0.5cm]
                \hat v(\bar\y,z) :=  
               \begin{cases}
                    \left( \dfrac{\overline{K}(\xx^*,\y)}{ S(\y - \y_0(z),\mathbf{n})}-\ell\left(\dfrac{{\bar \y}-{\bar \y}_0(z)}{\|{\bar \y}-{\bar \y}_0(z)\|};z\right) \right)\,V({\y})\,, & \ \y\neq \y_0(z)\,,\\[0.5cm]
                    \,0\,, & \ \y=\y_0(z)\,.
                \end{cases} \\[0.9cm]
            \end{array}\right\}  
        \end{array}\right. \nonumber\\[0.4cm]
        \hline\\
         \text{In}&\text{{ Section \ref{sub:section52-weights}: }} \\  
         &\left\{
         \begin{array}{l}
         \begin{array}{l}
         \omega_k := \omega(kh)\,,  \\[0.3cm]
         \omega(z) := \omega[ s(\bar\y,z);\,\alpha(z),\beta(z) ]  \\[0.2cm]
         \hspace{1cm} \approx c_0\,\omega\left[\dfrac{1}{\|\bar\y\|};\,\alpha,\beta\right] + \sum_{j=1}^N\Bigg\{c_j\, \omega\left[\dfrac{\cos(2j\,\psi(\bar\y))}{\|\bar\y\|};\,\alpha(z),\beta(z)\right] \\ 
         \hspace{1.5cm} + d_j\, \omega\left[\dfrac{\sin(2j\,\psi(\bar\y))}{\|\bar\y\|};\,\alpha(z),\beta(z)\right]\Bigg\}\,,
         \end{array}\\[1.5cm]
         \begin{array}{l}
         \text{where $c_0$, $c_j,d_j$ come from }\\[0.2cm]
         s( (\bar\y,z),\,\mathbf{n} )\,\|\bar\y\|\, \approx\, c_0+\sum_{j=1}^N \left\{ c_j \cos(2j\,\psi(\bar\y))+ d_j \sin(2j\,\psi(\bar\y)) \right\}\,,\\[0.2cm]
         \text{and } \bar\y = \|\bar\y\| \left(\,\cos(\psi(\bar\y)),\,\sin(\psi(\bar\y))\right)\,.
         \end{array}
         \end{array}
         \right. \nonumber\\[0.4cm]
         \hline
    \end{align*}
\end{table}

\subsection{Correction formula for the three kernels}
\label{sub:corr_allkernels}
Let $\mathbf{x}^*\in\Gamma$ be a target point with normal to the surface $\mathbf{n}=(n_1,n_2,n_3)$, $n_3\neq 0$, $\|\mathbf{n}\|=1$. 
{ We want to apply the three-dimensional second order correction formula \eqref{eq:Q1_3D} to the layer potential $\int_{\R^3}\overline{K}(\xx^*,\y)\rho(\y)\delta_{\Gamma,\varepsilon}(\y)\text{d}\y$ with one of the three layer kernels \eqref{eq:IBIM_kernels}, for example the single-layer kernel $K(\xx,\y)=G_0(\xx,\y)=(4\pi\|\xx-\y\|)^{-1}$. }
The starting point is the following singular function:
\begin{equation}\label{eq:s-the-main-singular-factor}
    S:(\mathbf{r},\mathbf{n})\in \left(\R^3\setminus\{\mathbf{0}\}\right)\times\R^3\mapsto \left(\frac{\left\| \mathbf{r}\times \mathbf{n} \right\|}{\|\mathbf{n}\|}\right)^{-1}\in\R\,,
\end{equation}
which represents the reciprocal of the distance from a point $\mathbf{r}$ to the line with direction $\mathbf{n}$ passing through the origin: $\{t\,\mathbf{n}\,:\,t\in\R\}$. 

\begin{figure}
    \begin{center}
    \includegraphics{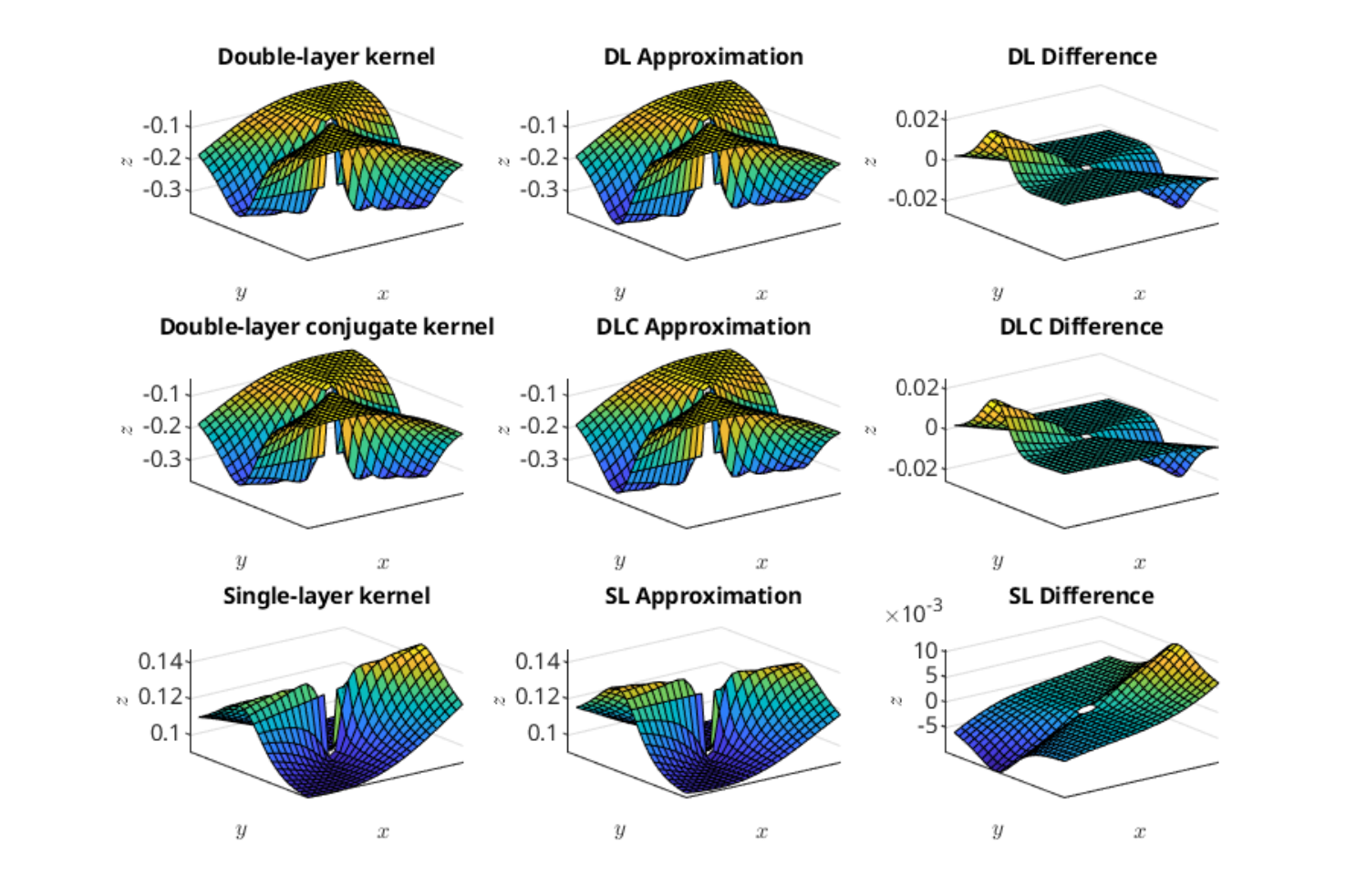}%-eps-converted-to
    \end{center}
    \caption{Discontinuous behavior}
    \footnotesize{Left column: kernel multiplied by the singular function \eqref{eq:s-the-main-singular-factor}. Center column: function $\ell$ {analytically} found by taking the limit \eqref{eq:ell_limit}. Right column: plot of the function used in the additive splitting $f/S-\ell$.}
    \label{fig:singularities}
\end{figure}

This choice is equivalent to approximating the distance on the denominator of the three kernels \eqref{eq:IBIM_kernels} as 
$\|\mathbf{x}^*-P_\Gamma\mathbf{y}\|\approx \|\mathbf{x}^*-\mathbf{y}_{TM(\mathbf{x}^*)}\|$, where $\mathbf{y}_{TM(\mathbf{x}^*)}$ is the projection of $\mathbf{y}$ onto the tangent plane to $\Gamma$ at $\mathbf{x}^*$.  
The reciprocal of the distance $\|\mathbf{x}^*-\mathbf{y}_{TM(\mathbf{x}^*)}\|$ is conveniently given by $S(\mathbf{x}^*-\mathbf{y},\mathbf{n})$.

With $S$, {we write} formula \eqref{eq:f_s_v_3D_splitting} {as}
\begin{equation}\label{eq:first_splitting_S_vtilde}
f(\mathbf{\bar y},z) =  S((\mathbf{\bar y},z) - \mathbf{y}_{0}(z), \mathbf{n}) \tilde{v}(\mathbf{\bar y},z). 
\end{equation}
The function $\tilde{v}(\,\cdot\,,z)$ defined via $f$ and $S$ in \eqref{eq:first_splitting_S_vtilde}, while bounded, is discontinuous at $\mathbf{\bar y}_0(z)$ for all three kernels. {In Figure~\ref{fig:singularities} the left column shows the behaviour of $f/S$ for $f$ written using the three Laplace layer kernels.} 

The next step is to isolate the discontinuous behavior in $\tilde{v}(\mathbf{\bar y},z)$. 
We observe {Figure~\ref{fig:singularities}} that $\tilde v$ has different limits at $\mathbf{\bar y}_0(z)$, depending on the approaching angle. So we will derive the function $\ell$ which has the same discontinuity. More precisely,
\begin{equation}\label{eq:ell_limit}
        {\ell}(\mathbf{q};z):=\lim_{t\rightarrow 0} \frac{\overline{K}(\mathbf{x}^*,(tq_1,tq_2,0)+{\mathbf{y}}_{0}(z))}{{S}((tq_1,tq_2,0),\mathbf{n})},~~~  \mathbf{q}=(q_1,q_2)\ ,\ \ \|\mathbf{q}\|=1\,.
\end{equation}

We write $\tilde v$ as
\begin{align}\label{eq:second_splitting_ell_v}
\tilde{v}(\mathbf{\bar y},z)=&\,\ell\left(\frac{\mathbf{\bar y}-\mathbf{\bar y}_0(z)}{\|\mathbf{\bar y}-\mathbf{\bar y}_0(z)\|};z\right) v(\mathbf{\bar y},z)%\ ,\ \ \text{where } \mathbf{\hat y}:=\,\frac{\mathbf{\bar y}}{\|\mathbf{\bar y}\|}\,.
\end{align}

This defines the new function $v$ which is smooth; hence $f=S \cdot\ell\cdot v$, and the weights for the corrected trapezoidal rule will be therefore derived for {the singular function
\begin{equation*}
    s(\y)=s(\bar\y,z)=S(\y,\mathbf{n})\,\ell\left(\dfrac{\bar\y}{\|\bar\y\|};z\right)\,.
\end{equation*}}

\subsubsection{Derivation of the formula for the new factor} \label{sub:section511-finding-ell}
Due to the closest point projection in $\overline{K}$, the limit function $\ell$ depends on the intrinsic geometry of $\Gamma$ (the principal curvatures) as well as the orientation and distance of $\mathbf{y}_0(z)$ to $\Gamma$ (the signed distance $\eta=\eta(z)$).
For different values of $z$ but fixed direction $\mathbf{q}$, $P_\Gamma$ may map the lines $(t\mathbf{q},0)+\mathbf{y}_0(z)$ to different curves on $\Gamma$. The speed at which these curves, $P_\Gamma((t\mathbf{q},0)+\mathbf{y}_0(z))$, pass through the target point $\mathbf{x}^*$ may vary, depending on the curvature of the curve. 

To calculate explicitly the limit that defines $\ell$, we will replace the projection $P_\Gamma$ by $P_{\tilde\Gamma}$, a high order local approximation of the closest point projection to the osculating paraboloid at $\mathbf{x}^*$. In the following, we will present the derivation of an explicit formula for $\ell$, based on
\begin{equation}\label{eq:ell_limit_using_paraboloid}
        {\ell}(\mathbf{q};z)=\lim_{t\rightarrow 0} \frac{{K}(\mathbf{x}^*,P_{\tilde\Gamma}((tq_1,tq_2,0)+{\mathbf{y}}_{0}(z)))}{{S}((tq_1,tq_2,0),\mathbf{n})},
\end{equation}
and building up from the simplest case.
In the derivations, we let $\bm{\tau}_1,\bm{\tau}_2,\mathbf{n}$ be the orthonormal basis of $\R^3$ composed of the principal directions $\bm{\tau}_i$, with corresponding principal curvature $\kappa_i$, ordered such that $\bm{\tau}_1\times\bm{\tau}_2=\mathbf{n}$. 
We write a point
in this basis as
\[
(a,b,c)_{T}:=a\bm{\tau}_1+b\bm{\tau}_2+c\,\mathbf{n}\,.%+\mathbf{x}^*\,.
\]
We first assume the $z$ planes are parallel to the tangent plane $TM(\mathbf{x}^*)$, i.e. $\mathbf{n}=\mathbf{e}_3$. Thus we want to find the limit 
\begin{equation}\label{eq:ell_limit_tangent_plane}
        \lim_{t\rightarrow 0} \frac{{K}(\mathbf{x}^*,P_{\tilde\Gamma}( 
        (tp_1,tp_2,\eta)_T+\mathbf{x}^*
        ))}{{S}((tp_1,tp_2,\eta)_T+\mathbf{x}^*,\mathbf{n})} ,\ ~~~\mathbf{p}=(p_1,p_2)\,,\ \ \sqrt{p_1^2+p_2^2}=1\,.
\end{equation}

\begin{figure}[ht]
    \begin{center}
    \includegraphics[scale=0.48]{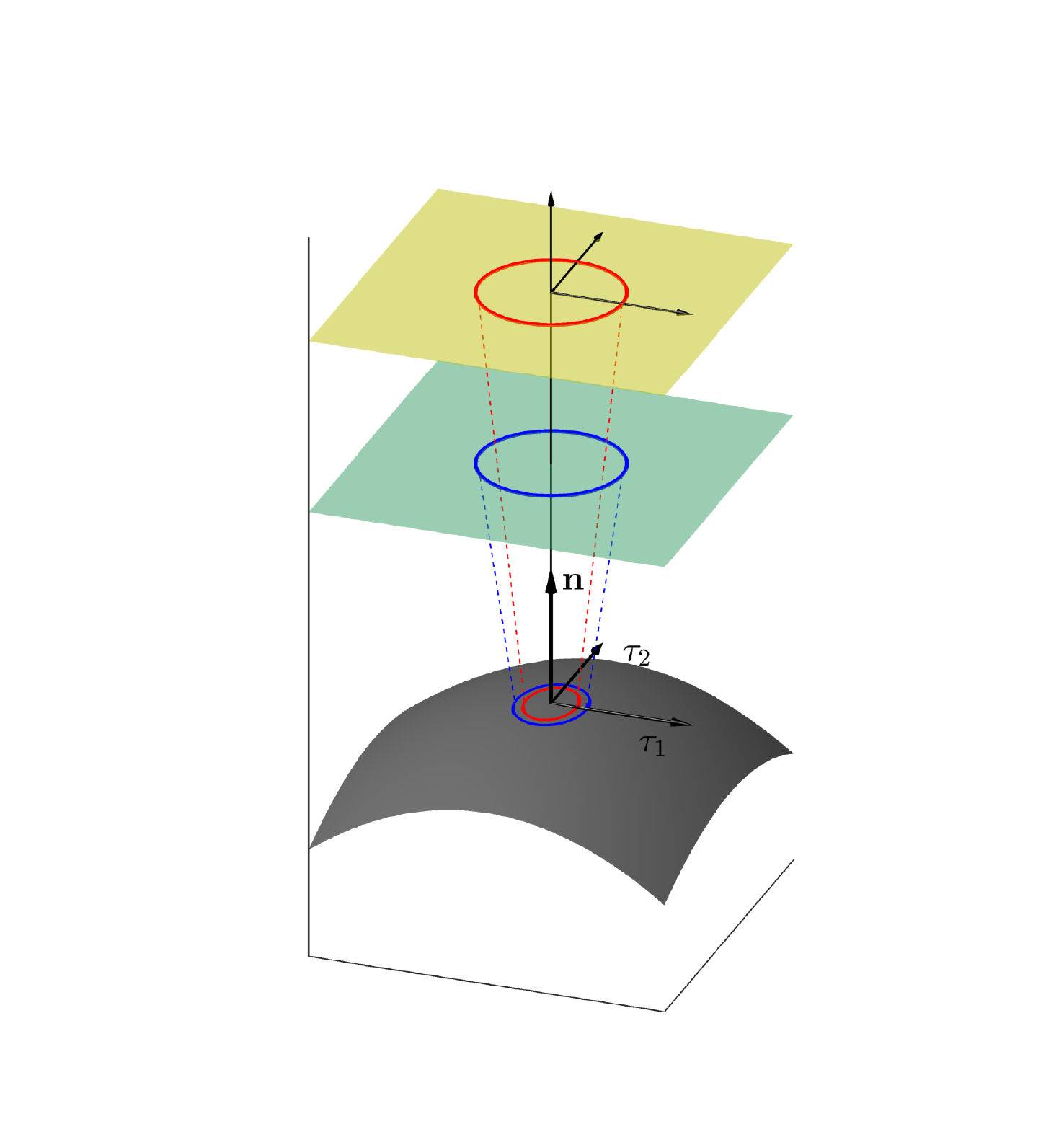}    \includegraphics[scale=0.48]{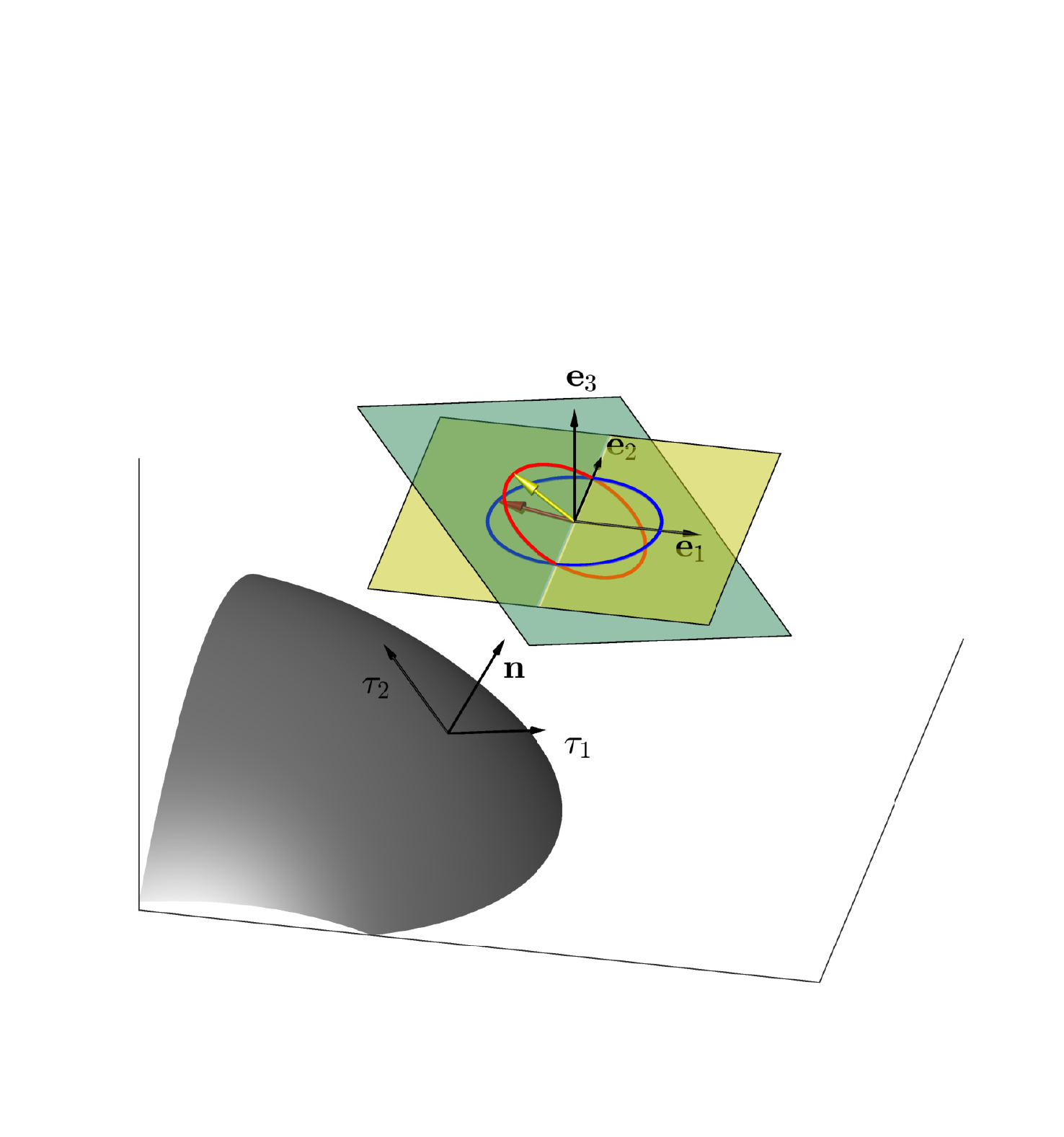}
    \end{center}
    \caption{Limit computation setting}
    \footnotesize{The surface is approximated around the target point with a paraboloid defined by the surface's principal curvatures and directions. 
    Points on the circles on each plane are mapped to the closest points on the paraboloid $\tilde\Gamma$ (instead of $\Gamma$) for calculation of the limit defined in \eqref{eq:ell_limit_using_paraboloid}. 
    Right plot: $\mathbf{n}\neq\vec{z}$; a circle (blue) drawn on the $z$ plane (yellow plane) around the singular point $\mathbf{y}_0(z)$ becomes an ellipse (red) when projected on the plane $\Pi_z$ (blue plane) parallel to the tangent plane of $\Gamma$ in $\mathbf{x}^*$. The angle $\psi$ between $\mathbf{e}_1$ and a given direction (red direction) on the plane $z$ will correspond to the angle $\bar \psi$ between $\bm{\tau}_1$ and the projected direction (yellow direction) on the plane $\Pi_z$.}
    \label{fig:ell_calculation_setting_plots}
\end{figure}

For convenience, we translate the problem so that $\mathbf{x}^*$ is in the origin, and work only in the $\bm{\tau}_1,\bm{\tau}_2,\mathbf{n}$ basis. See a representation of this setting in the left plot of Figure~\ref{fig:ell_calculation_setting_plots}.
The paraboloid $\tilde \Gamma$ will then be 
\[
\tilde \Gamma = \left\{ \left(x,y,\frac{1}{2}(\kappa_1x^2+\kappa_2y^2)\right)_T\,:\,x,y\in\R \right\}\,.
\]
In a sufficiently small neighborhood of the origin, a point $(x,y,\eta)_T$  and its 
closest point $$P_{\tilde\Gamma}((x,y,\eta)_T)=(\bar x, \bar y, \left(\kappa_1\bar x^2+\kappa_2\bar y^2\right)/2)_T$$
satisfy
\[
    (x,y,\eta)_T-\left(\bar x, \bar y, \frac{1}{2}\left(\kappa_1\bar x^2+\kappa_2\bar y^2\right)\right)_T=\left(\eta-\frac{1}{2}\left(\kappa_1\bar x^2+\kappa_2\bar y^2\right)\right) (-\kappa_1\bar x,-\kappa_2 \bar y,1)_T\,;
\]
i.e. the vector pointing to the closest point on $\tilde\Gamma$ should be normal to the surface, with magnitude equal to the distance to the surface. 
Along $(tp_1,tp_2,\eta)_T$,  
we have: 
\begin{equation}\label{eq:xbar_ybar_projection_asymptotics}
    \begin{cases}
    tp_1 = \bar p_1-\kappa_1\bar p_1[\eta-(\kappa_1\bar p_1^2+\kappa_2\bar p_2^2)/2] \\[0.3cm]
    tp_2 = \bar p_2-\kappa_2\bar p_2[\eta-(\kappa_1\bar p_1^2+\kappa_2\bar p_2^2)/2]
    \end{cases}\ \Rightarrow \ \
    \begin{cases}
    \bar p_1 = \dfrac{t}{1-\kappa_1\eta}~p_1+\mathcal{O}(t^3)\,,\\[0.3cm]
    \bar p_2 = \dfrac{t}{1-\kappa_2\eta}~p_2+\mathcal{O}(t^3)\,.
    \end{cases}
\end{equation}

For example if $\mathbf{p}=(1,0)$ the limit is taken along the $\bm{\tau}_1$ direction, the projected point will travel along the curve corresponding to the first principal direction, and the limit value for the double-layer conjugate kernel will be:
\begin{align*}
\lim_{t\rightarrow 0} \left.\frac{\overline{K}(\mathbf{x}^*,(tp_1,tp_2,\eta)_T)}{{S}((tp_1,tp_2,\eta)_T,\mathbf{n})}\right|_{\mathbf{p}=(1,0)}& = \lim_{t\to 0^+} -\frac{1}{4\pi}\frac{(\mathbf{x}^*-P_{\tilde\Gamma}((t,0,\eta)_T))\cdot\mathbf{n}}{\|\mathbf{x}^*-P_{\tilde\Gamma}((t,0,\eta)_T)\|^3}\|(t,0,\eta)_T-(0,0,\eta)_T\|\\
&=\lim_{t\to0^+} -\frac{1}{4\pi}\dfrac{ -\frac{1}{2}\kappa_1\frac{t^2}{(1-\kappa_1\eta)^2}+\mathcal{O}(t^4) }{\left[\frac{t^2}{(1-\kappa_1\eta)^2}+\mathcal{O}(t^4)+\frac{1}{4}\frac{\kappa_1^2t^4}{(1-\kappa_1\eta)^4}+\mathcal{O}(t^6)\right]^{3/2}}t \\ &=\frac{1}{8\pi}\kappa_1(1-\kappa_1\eta)\,.
\end{align*}
In general, we have
\begin{align}
\nonumber \lim_{t\rightarrow 0} &\frac{\overline{K}(\mathbf{x}^*,(tp_1,tp_2,\eta)_T)}{{S}((tp_1,tp_2,\eta)_T,\mathbf{n})}
=\lim_{t\to 0^+} -\frac{1}{4\pi}\frac{(\mathbf{x}^*-P_{\tilde\Gamma}((tp_1,tp_2,\eta)_T))^T\mathbf{n}}{\|\mathbf{x}^*-P_{\tilde\Gamma}((tp_1,tp_2,\eta)_T)\|^3}\left\|(tp_1,tp_2,\eta)_T-(0,0,\eta)_T\right\|\nonumber\\
\nonumber &=\lim_{t\to0^+} -\frac{1}{4\pi}\dfrac{ -\frac{1}{2}\kappa_1\frac{t^2p_1^2}{(1-\kappa_1\eta)^2}-\frac{1}{2}\kappa_2\frac{t^2p_2^2}{(1-\kappa_2\eta)^2}+\mathcal{O}(t^4) }{\left\{\frac{t^2p_1^2}{(1-\kappa_1\eta)^2}+\frac{t^2p_2^2}{(1-\kappa_2\eta)^2}+\mathcal{O}(t^4)+\frac{1}{4}\left[\frac{\kappa_1^2p_1^2}{(1-\kappa_1\eta)^4}+\frac{\kappa_2^2p_2^2}{(1-\kappa_2\eta)^4}+\mathcal{O}(t^2)\right]^2t^4\right\}^{3/2}}t\\
&=\frac{1}{8\pi}\dfrac{\kappa_1\dfrac{p_1^2}{(1-\kappa_1\eta)^2}+\kappa_2\dfrac{p_2^2}{(1-\kappa_2\eta)^2}}{\left[ \dfrac{p_1^2}{(1-\kappa_1\eta)^2}+\dfrac{p_2^2}{(1-\kappa_2\eta)^2} \right]^{3/2}}\,. \label{eq:ell_limit_first_case}
\end{align}

We now consider the more general case in which $\mathbf{n}\neq\mathbf{e_3}$. Again we consider the target point to be in the origin: $\mathbf{x}^*=\mathbf{0}$. We define the plane $\Pi_z:=\{(x,y,\eta(z))_T\,:\,x,y\in\R\}$ parallel to the tangent plane $TM(\mathbf{x}^*)$ at distance $\eta(z):=d_\Gamma(\mathbf{y}_0(z))$. Fixed $z\in\R$, the projection
\[
P_{\Pi_z}:\mathbf{x}=(x_1,x_2,x_3)\in\R^3\mapsto (I-\mathbf{n}\otimes\mathbf{n})\mathbf{x}+\frac{z}{n_3}\mathbf{n}\in\Pi_z 
\]
takes a point $\mathbf{x}$ to the intersection of the line $\{\mathbf{x}+t\mathbf{n}\,:\,t\in\R\}$ and the plane $\Pi_z$.

Let $\mathbf{q}=(q_1,q_2)\in\mathbb{S}^1$, and let {$(tq_1, tq_2,0)+\mathbf{y}_0(z)$} be a point on the $z$ plane. To find the limit, we first consider the projection of the line $(tq_1,tq_2,0)+\mathbf{y}_0(z)$ onto the plane $\Pi_z$, and then apply the previous formula; in the expression of 
\eqref{eq:ell_limit_first_case} 
we consequently have 
\[
(t\,p_1(q_1,q_2),t\,p_2(q_1,q_2),\eta(z))_T=P_{\Pi_z}((tq_1,tq_2,0)+\mathbf{y}_0(z))
\]
instead of $(tq_1,tq_2,0)+\mathbf{y}_0(z)$. This change does not affect the limit expression because of the property
\begin{align*}
    P_{\tilde \Gamma}\left( (tx,ty,0)+\mathbf{y}_0(z) \right) &= P_{\tilde \Gamma}\left( (tx,ty,\eta(z))_T \right) +\mathcal{O}(t^2)
\end{align*}
for small values of $t$, which expresses how the orientation of the plane $z$ with respect to the basis $\bm{\tau}_1,\bm{\tau}_2,\mathbf{n}$ does not affect significantly the projection of points close to the singularity point $\mathbf{y}_0(z)$.

Let  
\[
(\tilde p_1,\tilde p_2,0)_T=P_{\Pi_0} (q_1,q_2,0)
\]
be the projection of $(q_1,q_2,0)$ onto the tangent plane
$TM(\mathbf{x}^*)=\Pi_{0}$; then
\begin{align}
    \tilde p_1&=\tilde p_1(q_1,q_2)=\dfrac{q_1\cos\theta_0-\dfrac{q_2}{c}(ab\cos\theta_0+(b^2+c^2)\sin\theta_0)}{\sqrt{1+\frac{(a\cos\theta_0+b\sin\theta_0)^2}{c^2}}}\,,\label{eq:p_transformation_1}\\
    \tilde p_2&=\tilde p_2(q_1,q_2)=\dfrac{q_1\sin\theta_0+\dfrac{q_2}{c}(ab\sin\theta_0+(a^2+c^2)\cos\theta_0)}{\sqrt{1+\frac{(a\cos\theta_0+b\sin\theta_0)^2}{c^2}}}\,,\label{eq:p_transformation_2}
\end{align}
and \eqref{eq:ell_limit_first_case} will be valid with
\begin{equation}\label{eq:transformation_normalization}
\mathbf{p}=(p_1,p_2)=(\tilde p_1,\tilde p_2)/\sqrt{\tilde p_1^2+\tilde p_2^2}\,.
\end{equation}

The expressions relating $\tilde p_1, \tilde p_2$ to $p_1,p_2$ can be visualized in the right plot of Figure~\ref{fig:ell_calculation_setting_plots}. The circle  $\{t(q_1,q_2,0)+\mathbf{y}_0(z): q_1^2+q_2^2=1\}$ on the plane $z$, projected onto the plane $\Pi_z$ 
will become an ellipse in general. 

The parameters $a,b,c,\theta_0$ relate $\bm{\tau}_1,\bm{\tau}_2$ to the standard $\R^3$ basis $\{\mathbf{e}_i\}_{i=1}^3$.
They are given by
$$
(a,b,c)=
(\sin\theta\cos\xi,\sin\theta\sin\xi,\cos\theta)\,,
$$ 
where $\theta$ is such that $\cos\theta={n}_3$, i.e. it is the second spherical coordinate of $\mathbf{n}$; $\xi$ is such that $\tan\xi=\mathbf{e}_3^T\bm{\tau}_2/\mathbf{e}_3^T\bm{\tau}_1$, hence it is the angle between $\bm{\tau}_1$ and the projection of $\mathbf{e}_3$ on the plane $\bm{\tau}_1,\bm{\tau}_2$; and $\theta_0$ is such that $\tan\theta_0=\mathbf{e}_1^T\bm{\tau}_2/\mathbf{e}_1^T\bm{\tau}_1$, meaning it is the angle between $\bm{\tau}_1$ and the projection of $\mathbf{e}_1$ on the plane $\bm{\tau}_1,\bm{\tau}_2$. 

{Given a unit vector $(q_1, q_2)$ on the $z$ plane, formulae (\ref{eq:p_transformation_1}-\ref{eq:p_transformation_2}) map it to the unit vector
$p_1\bm{\tau}_1+ p_2\bm{\tau}_2$. The unit direction $(p_1,p_2)$ on the plane $\Pi_z$, $(p_1,p_2,0)_T$, is the corresponding direction in which the limit (in the definition \eqref{eq:ell_limit_using_paraboloid} of $\ell(\mathbf{q}; z)$) will be evaluated. Recall that the formula for the limit for any given direction is already derived in \eqref{eq:ell_limit_first_case}. For convenience, we write the unit vector $(q_1, q_2)=(\cos\psi,\sin\psi)$, and correspondingly we will treat $p_1, p_2$, defined in (\ref{eq:p_transformation_1}-\ref{eq:transformation_normalization}), as functions of $\psi$.}

Thus we write the formulae for $\ell$ for the double-layer conjugate, double-layer, and single-layer kernels as: 
\begin{align}\label{eq:ell_DL_DLC}
    \text{(DL-DLC): }\,\ell(\psi,z)&=\frac{1}{8\pi}\dfrac{\kappa_1\dfrac{p_1^2(\psi)}{(1-\kappa_1\eta(z))^2}+\kappa_2\dfrac{p_2^2(\psi)}{(1-\kappa_2\eta(z))^2}}{\ \left\{ \dfrac{p_1^2(\psi)}{(1-\kappa_1\eta(z))^2}+\dfrac{p_2^2(\psi)}{(1-\kappa_2\eta(z))^2} \right\}^{3/2}}\,,
    \\[0.4cm]
    \text{(SL): }\,\ell(\psi,z)&=\frac{1}{4\pi}{\ \left\{ \dfrac{p_1^2(\psi)}{(1-\kappa_1\eta(z))^2}+\dfrac{p_2^2(\psi)}{(1-\kappa_2\eta(z))^2} \right\}^{-\frac{1}{2}}}\,.\label{eq:ell_SL}
\end{align}
Again,  $\eta(z)=d_\Gamma(\mathbf{y}_0(z))$ is the signed distance of $\mathbf{y}_0(z)$ to the surface. \\

{The center and right column in Figure~\ref{fig:singularities} illustrate how the formulae found approximate the behaviour $f/S$. The center column plots the function $\ell$ in (\ref{eq:ell_DL_DLC}-\ref{eq:ell_SL}) found for the three Laplace layer kernels, while the right column shows the difference between the function $\ell$ found and the values $f/S$.}

\subsubsection{The quadrature formulae} \label{sub:section512-factorization-splitting}
With $\ell$ defined above, we will work with the following singular function-smooth function factorization:
\begin{align}\label{eq:splitting_allkernels}
    &f(\mathbf{\bar y},z) = s(\mathbf{\bar y}-\mathbf{\bar y}_0(z),z)v(\mathbf{\bar y},z) \vspace{0.5cm}\\[0.2cm]
    \text{with }\,& s(\mathbf{\bar y},z) = 
    S((\mathbf{\bar y},z),\mathbf{n})\,\ell\left( \mathbf{\bar y}/\|\mathbf{\bar y}\|; z \right)\nonumber %\\
\end{align}
where the weight for the corrected trapezoidal rule applied to $f$ is computed with $s=S\cdot\ell$. The function $s$ completely captures the asymptotic behavior of the given kernel $\overline{K}$ in $\mathbf{\bar y}_0(z)$
and we can apply \eqref{eq:single_correction_quadrature} to $f(\mathbf{\bar y},z)$: 
\begin{equation}\label{eq:factorize_quad_non_additive}
    {\bar Q}^{2D}_h[f(\cdot,z)]=
    T_h^{0}[s(\,\cdot-\mathbf{\bar y}_0(z),z)\, v(\,\cdot\,)] + h \, \omega[s(\,\cdot\,,z);\alpha(z),\beta(z)] \dfrac{f(\mathbf{\bar y}_{\Delta},z)}{s(\mathbf{\bar y}_{\Delta}-\mathbf{\bar y}_0(z),z)}\,, 
\end{equation}
if $(\alpha(z),\beta(z))\neq(0,0)$, otherwise
\begin{equation*}
    {\bar Q}^{2D}_h[f(\,\cdot\,,z)]=
    T_h^{0}[s(\,\cdot-\mathbf{\bar y}_0(z),z)\, v(\,\cdot\,)]+h\,\omega[s(\,\cdot\,,z);0,0]\,v(\mathbf{\bar y}_{0},z)\,.
\end{equation*}
As long as $\ell$ is non-zero, 
the function $v$ is well-defined
through \eqref{eq:splitting_allkernels}, away from $\mathbf{\bar y}_0(z)$
and by continuity
at $\mathbf{\bar y}_0(z)$.
As can be seen from
\eqref{eq:ell_DL_DLC} and
\eqref{eq:ell_SL},
this is always the case for the
single-layer kernel
but for the double-layer
kernels in general only if $\kappa_1$
and $\kappa_2$ have the same sign.
Even though $\ell$ is zero only
at isolated points,
we cannot apply formula~\eqref{eq:factorize_quad_non_additive} as it is numerically problematic. We need a different approach
which works as follows. 

We use the discontinuity subtraction from $\overline{K}$: first,
to shorten the formulae below, we define
the smooth function
\[
V(\mathbf{\bar y},z) := \rho(\mathbf{\bar y},z)\delta_{\Gamma,\varepsilon}(\mathbf{\bar y},z)
\]
so that $f(\mathbf{\bar y},z)=\overline{K}(\mathbf{x}^*,(\mathbf{\bar y},z))V(\mathbf{\bar y},z)$. 
We then replace the
splitting in 
\eqref{eq:splitting_allkernels}
by
$$
   f(\mathbf{\bar y},z) = S((\mathbf{\bar y},z)-(\mathbf{\bar y}_0(z),z),\mathbf{n})\hat{v}(\mathbf{\bar y},z)+
   s(\mathbf{\bar y}-\mathbf{\bar y}_0(z),z)V(\mathbf{\bar y},z),
$$
where
$$
   \hat{v}(\mathbf{\bar y},z)=
\left(   \frac{\overline{K}(\mathbf{x}^*,(\mathbf{\bar y},z))}{
   S((\mathbf{\bar y},z)-(\mathbf{\bar y}_0(z),z),\mathbf{n})}-\ell\left(\frac{\mathbf{\bar y}-\mathbf{\bar y}_0(z)}{\|\mathbf{\bar y}-\mathbf{\bar y}_0(z)\|};z\right) \right)
   \,V(\mathbf{\bar y},z)\,.
$$
Then $\hat{v}$ is well-defined
everywhere, bounded and continuous around
$\mathbf{\bar y}_{0}$, by construction of $\ell$ via the limit \eqref{eq:ell_limit}.
%; moreover the terms are both bounded close to the singularity point.
We therefore rewrite  \eqref{eq:factorize_quad_non_additive} as
\begin{align} 
{\bar Q}^{2D}_h&[f(\,\cdot\,,z)] 
= {\bar Q}^{2D}_h[\overline{K}(\mathbf{x}^*,(\,\cdot\,,z))V(\,\cdot\,,z)] \nonumber\\
=& {\bar Q}^{2D}_h\Bigg[ \left(\frac{\overline{K}(\mathbf{x}^*,(\,\cdot\,,z))}{S((\,\cdot\,,z)-\mathbf{ y}_0(z),\mathbf{n})}-\ell\left(\frac{\,\cdot-\mathbf{\bar y}_0(z)}{\|\,\cdot-\mathbf{\bar y}_0(z)\|};z\right) \right)S((\,\cdot\,,z)-\mathbf{ y}_0(z),\mathbf{n})\,V(\,\cdot,z) \nonumber\\
&\hspace{1.0cm} +S((\,\cdot\,,z)-\mathbf{ y}_0(z),\mathbf{n})\,\ell\left(\frac{\,\cdot-\mathbf{\bar y}_0(z)}{\|\,\cdot-\mathbf{\bar y}_0(z)\|};z\right)V(\,\cdot,z) \Bigg]\nonumber \\ 
=& T_h^{0}\left[\overline{K}(\mathbf{x}^*,(\,\cdot\,,z))V(\,\cdot\,,z)\right]\nonumber\\
&\hspace{0.5cm} +h\Bigg\{ \omega[{S((\,\cdot\,,z)-\mathbf{ y}_0(z),\mathbf{n})};\alpha(z),\beta(z)]\nonumber\\
&\hspace{1.5cm}\cdot \left( \frac{K(\mathbf{x}^*,(\mathbf{\bar y}_{\Delta},z))}{S((\mathbf{\bar y}_\Delta,z)-(\mathbf{\bar y}_0(z),z),\mathbf{n})}-\ell\left(\frac{{\mathbf{\bar y}}_{\Delta}-\mathbf{\bar y}_0(z)}{{\|\mathbf{\bar y}}_{\Delta}-\mathbf{\bar y}_0(z)\|};z\right) \right)  \nonumber\\ 
& \hspace{1.5cm} + \omega[s(\,\cdot\,,z);\alpha(z),\beta(z)] \Bigg\}V(\mathbf{\bar y}_{\Delta},z)\,,\label{eq:additivesplit}
\end{align}
where inside the braces only the second term remains if $(\alpha(z),\beta(z))=(0,0)$. 

Our corrected trapezoidal rule \eqref{eq:additivesplit} for the implicit boundary integral takes the form:
\begin{align}
    Q_h[f]=& h^3\sum_{\mathbf{y}\in \left(h\mathbb{Z}^3\setminus N_h(\mathbf{x}^*)\right)}f(\mathbf{y}) + h^2\sum_{k\in\mathbb{Z}}V(\mathbf{y}_{\Delta}(z_k)) \mathcal{R}_{h,k}(\mathbf{y}_{\Delta}(z_k)), \label{eq:Q1_3D_additivesplit}
\end{align}
where
\begin{align}
    \label{eq:additive_split_correction}
    &\mathcal{R}_{h,k}(\mathbf{y}_{\Delta}(z_k)) :=\\ \nonumber
    &= 
    \begin{cases}
    \omega[s(\cdot,z_k);\alpha(z_k),\beta(z_k)] + \omega[S((\,\cdot\,,z_k)-\mathbf{y}_0(z_k),\mathbf{n});\alpha(z_k),\beta(z_k)]\cdot & \\[0.3cm]
    \hspace{0.6cm} \cdot\left( \dfrac{ \overline{K}(\mathbf{x}^*,\mathbf{y}_{\Delta}(z_k)) }{ S((\mathbf{\bar y}_{\Delta}(z_k),z_k)-\mathbf{y}_0(z_k),\mathbf{n}) } -\ell\left(  \dfrac{\mathbf{\bar y}_{\Delta}(z_k)-\mathbf{\bar y}_{0}(z_k)}{\|\mathbf{\bar y}_{\Delta}(z_k)-\mathbf{\bar y}_{0}(z_k)\|};z_k\right) \right) & \\
    &  \text{ if } \alpha(z_k),\beta(z_k)\neq 0, \\ & \\
    \omega[s(\cdot,z_k);0,0], & \text{ otherwise.}
    \end{cases}
\end{align}

\begin{comment}
    \begin{align}
        \label{eq:additive_split_correction}
        &\mathcal{R}_{h,k}(\mathbf{y}_{\Delta}(z_k)) :=\\ \nonumber
        &= 
        \begin{cases}
        \omega[S((\,\cdot\,,z_k)-\mathbf{y}_0(z_k),\mathbf{n});\alpha(z_k),\beta(z_k)]
        \left( \dfrac{ \overline{K}(\mathbf{x}^*,\mathbf{y}_{\Delta}(z_k)) }{ S((\mathbf{\bar y}_{\Delta}(z_k),z_k)-\mathbf{y}_0(z_k),\mathbf{n}) } -\ell\left(  \dfrac{\mathbf{\bar y}_{\Delta}(z_k)-\mathbf{\bar y}_{0}(z_k)}{\|\mathbf{\bar y}_{\Delta}(z_k)-\mathbf{\bar y}_{0}(z_k)\|};z_k\right) \right)
        &  \\
        \ +\,\omega[s(\cdot,z_k);\alpha(z_k),\beta(z_k)]\,,
        & \text{ if } \alpha(z_k),\beta(z_k)\neq 0, \\ & \\
        \omega[s(\cdot,z_k);0,0], & \text{ otherwise.}
        \end{cases}
    \end{align}
\end{comment}
The general correction form \eqref{eq:corrected_trapezoidal_general_form} of this method is then valid for $\mathcal{N}_h(x)$ as in \eqref{eq:puncturedTRcutoffregion3D} and 
$$
\mathcal{R}_h(\mathbf{x})=h^2\sum_{k\in\mathbb{Z}}V(\mathbf{y}_{\Delta}(z_k)) \mathcal{R}_{h,k}(\mathbf{y}_{\Delta}(z_k))\,.
$$

\subsection{Approximation and tabulation of the weights} \label{sub:section52-weights}

We approximate the singular functions using a Fourier interpolation, then tabulate the weights for the simpler singular terms of the expansion, and compose the general weights for any behavior needed. 

Given a $\pi$-periodic function $\ell(\,\cdot\,;z)$, such as \eqref{eq:ell_DL_DLC} or \eqref{eq:ell_SL}, we wish to compute the weight $\omega[s(\,\cdot\,,z);\alpha,\beta]$, where $s$ comes from the factorization \eqref{eq:splitting_allkernels}:
\[
s(\mathbf{\bar y},z)=S((\mathbf{\bar y},z),\mathbf{n})\ell(\mathbf{\bar y}/\|\mathbf{\bar y}\|;z)\,.
\]
Both the factors in this expression can be seen as functions of the angle of approach, $\psi$, to the singular point $\mathbf{0}$:
\[
S((\mathbf{\bar y},z),\mathbf{n})\,\ell\left( {\mathbf{\bar y}}/{\|\mathbf{\bar y}\|}; z \right) = \frac{S_{\mathbf{n}}(\psi(\mathbf{\bar y}))}{\|\mathbf{\bar y}\|}\ell(\psi(\mathbf{\bar y});z)\ \ \text{ where } \ \ 
\mathbf{\bar y}=\|\mathbf{\bar y}\|(\cos(\psi(\mathbf{\bar y})),\sin(\psi(\mathbf{\bar y})))\,.
\]
Given the dependence on $z$ is only present in $\ell$ through $\eta(z)$, we write $\omega[s;\alpha,\beta]$ instead of $\omega[s(\cdot,z);\alpha,\beta]$.

We use Fourier interpolation to approximate this function with a trigonometric polynomial:
\begin{equation}\label{eq:s_ell_Fourier}
S_{\mathbf{n}}(\psi){\ell}(\psi;z)\approx c_{0}+\sum_{j=1}^{N}\left[c_{j}\cos(2j\psi)+d_{j}\sin(2j\psi)\right] \,.
\end{equation}
Since the weight $\omega[s;\alpha,\beta]$ is a linear functional of the singular function, 
\begin{align}\label{eq:om1,4_lgen}
\omega[s;\alpha,\beta] \approx&\ c_{0}\,\omega\left[\frac{1}{\|\mathbf{\bar y}\|};\alpha,\beta\right]+\\
\nonumber & +\sum_{j=1}^{N}\left\{c_{j}\,\omega\left[\frac{\cos(2j\psi(\mathbf{\bar y}))}{\|\mathbf{\bar y}\|};\alpha,\beta\right]+d_{j}\,\omega\left[\frac{\sin(2j\psi(\mathbf{\bar y}))}{\|\mathbf{\bar y}\|};\alpha,\beta\right]\right\}
\,.
\end{align}

Therefore, we can precompute the weights for the basis functions and for certain values of $\alpha$ and $\beta$:
\begin{equation}
\left\{\omega\left[\frac{\cos(2j\psi(\mathbf{\bar y}))}{\|\mathbf{\bar y}\|};\alpha,\beta\right]\right\}_{j=0}^N \ \text{ and } \ \left\{\omega\left[\frac{\sin(2j\psi(\mathbf{\bar y}))}{\|\mathbf{\bar y}\|};\alpha,\beta\right]\right\}_{j=1}^{N}\,.
\end{equation}
For values outside of the precomputed tables, we interpolate.

In the evaluation of  $\omega\left[\frac{\cos(2j\psi(\mathbf{\bar y}))}{\|\mathbf{\bar y}\|};\alpha,\beta\right]$ using formulae \eqref{eq:single_correction_weight-before-taking-limit} and \eqref{eq:single_correction_weight}, we need the value of the integral:\\ $\iint_{\R^2} \frac{\cos(2j\psi(\mathbf{\bar y}))}{\|\mathbf{\bar y}\|} g(\mathbf{\bar y})\text{d}\mathbf{\bar y}$ for some test function $g$. 
We use $g(\mathbf{\bar y})=\exp(-\|\mathbf{\bar y}\|^8)$  which has derivatives $\partial^{\mathbf{k}}g(\mathbf{0})=0$, $\mathbf{k}\in\mathbb{N}^2,\ 0<|\mathbf{k}|<8$. Then
\begin{align*}
\iint_{\R^2} \frac{\cos(2j\psi(\mathbf{\bar y}))}{\|\mathbf{\bar y}\|} g(\mathbf{\bar y})\text{d}\mathbf{\bar y} &=\int_0^{2\pi}\text{d}\psi \int_0^\infty\text{d}r\left\{ \frac{{r\,\exp(-r^8)}\cos(2j\psi)}{r} \right\}\\
&=\left(\int_0^{2\pi} \cos(2j\psi) \text{d}\psi\right)\left( \int_0^\infty \exp(-r^8) \text{d}r\right)\\
&=2\pi \delta_{0j}\int_0^\infty \exp(-r^8) \text{d}r\approx 2\pi \delta_{0j}\int_0^R \exp(-r^8) \text{d}r\,,
\end{align*}
where $\delta_{0j}$ is the Kronecker delta; $R$ is set to 1.9 as the integrand is essentially zero at double precision. 
 The integral $\int_0^R \exp(-r^8) \text{d}r$ is computed once with high precision using common integration libraries and reused for all instances as a constant. \\

It is impractical to tabulate precomputations of the coefficients $\{c_j\}_{j=0}^N$, $\{d_j\}_{j=1}^N$ as they depend on too many variables.
Instead, we compute $c_j$ and $d_j$ on the fly, by solving the square linear system 
\begin{equation*}
    \begin{aligned}
    c_{0}+\sum_{j=1}^{N}\left[c_{j}\cos(2j\psi_i)+d_{j}\sin(2j\psi_i)\right] =& S_{\mathbf{n}}(\psi_i)\ell(\psi_i;z),~~\psi_i=\frac{i\pi}{2N+1},~i=0,\cdots, 2N.
    \end{aligned}
\end{equation*}
Because of the smoothness of the $\pi$-periodic functions we deal with, we can use a relatively small $N$ in order to accurately approximate the weights.

In the next Section, we present numerical convergence studies using weights computed with $N=22$. The small linear system can be inverted efficiently, e.g. using FFT, with a negligible computational time. 

In the convergence studies, an array of weights
of dimensions (45,101,101) has been precomputed, with $101$ values
for $\alpha$ and $\beta$ in $\left[-\frac{1}{2},\frac{1}{2}\right]$
and 45 for the Fourier series with $N=22$. 
Biquintic interpolation is used to approximate the weights for given $(\alpha,\beta)$ outside of the precomputed values.

In the case $\ell\equiv 1$,  $\{c_j\}_{j=0}^N$ and  $\{d_j\}_{j=1}^N$ depend only on $\theta$ and $\phi$. In that case, they can also be precomputed, stored, and used in an interpolation process when needed.

\section{Numerical Examples}

\label{sec:numerical_examples}

We demonstrate the convergence and accuracy of the proposed quadrature rules by evaluating the double layer potential with constant density on the surface $\Gamma\subset\mathbb{R}^3$. {We demonstrate the numerical errors computed by the proposed corrected trapezoidal rule for approximating
\[
I=\int_\Gamma \frac{\partial G_0}{\partial\mathbf{n}_y}\left(\mathbf{x}^*,\mathbf{y}\right)\text{d}\sigma_{\mathbf{y}}\,.
\]
The value of $I$ is known explicitly to be $-1/2$ for any $\mathbf{x}^*\in\Gamma$. So, we report 
\begin{equation}\label{eq:constant_density_error}
E_1(h):=\left| Q_h\left[\frac{\partial G_0}{\partial\mathbf{n}_y}\left(\mathbf{x}^*,\,\cdot\,\right)\right]+\frac{1}{2}\right|\,,     
\end{equation}
for several randomly chosen $\mathbf{x}^*\in\Gamma$. We will compare the results for the four different quadrature rules, including the two new quadrature rules $Q^L_{IBIM}$ defined by the regularization \eqref{eq:KregIBIM-linear} and $Q_h$ defined in \eqref{eq:Q1_3D_additivesplit}.}

%Mainly, we present the errors in computing the following identity:
%$$\int_{\Gamma}\frac{\partial G_{0}}{\partial\mathbf{n}_y}(x,y)\text{d}\sigma_y=-\frac{1}{2}\ \  ,\ \ \ x\in\Gamma\,.$$
The integral is first extended to the tubular neighborhood of  $T_\varepsilon$, as in \eqref{eq:general_S}-\eqref{eq:restriction_of_K}, using the compactly supported $C^{\infty}$ averaging function 
\begin{equation}\label{eq:averagingfunction}
    \phi(x)=
    \begin{cases}
        a\,\exp\left( \dfrac{2}{x^2-1} \right), & \text{ if } |x|<1, \\
        0, & \text{ otherwise};
    \end{cases}
\end{equation}
here $a\approx {\texttt{7.51393}}$ normalizes the integral $\int_{\mathbb{R}} \phi(x) \text{d}x$ to 1.

\begin{figure}
    \begin{center}
    \includegraphics{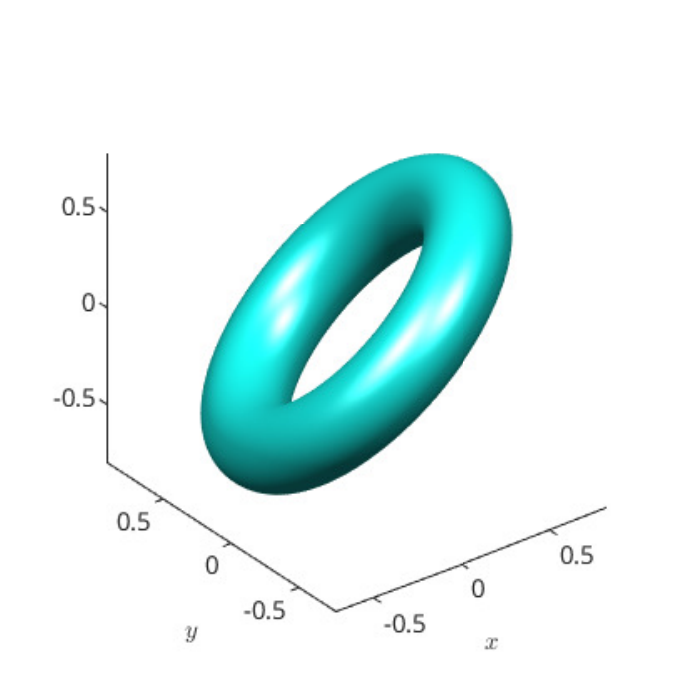}\includegraphics{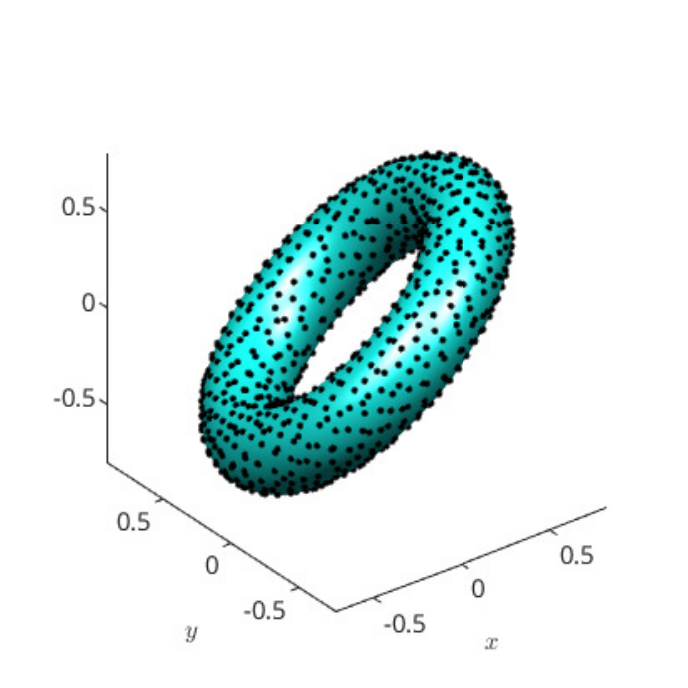}
    \end{center}
    \caption{Torus}
    \footnotesize{Left: the torus used in the tests. Right: the torus and the projections of the Cartesian grid nodes inside the tubular neighborhood $T_\varepsilon$. The projected nodes serve as the quadrature nodes.}
    \label{fig:tiltedtorusplot}
\end{figure}

The surfaces chosen for the tests are a sphere and a torus, centered in a random point in 3D and rotated with random angles along the $x$-, $y$- and $z$-axes. 

The sphere is characterized by center 
\[ \mathbf{C}=(\texttt{0.5475547095598521}, \texttt{0.6864792402110276}, \texttt{0.3502726366462485})\cdot \texttt{10}^{\texttt{-1}} \]
and radius $R=0.7$. The torus is described by the following parametrization
\begin{equation}\label{eq:torus}
    \mathcal{T}(\theta,\phi)=Q\left(
    \begin{matrix}
    (R_2\cos\theta+R_1)\cos\phi \\
    (R_2\cos\theta+R_1)\sin\phi \\
    R_2\sin\theta
    \end{matrix}
    \right)+\mathbf{C}
\end{equation}
{where $R_1=0.7$, $R_2=0.2$, $\mathbf{C}$ is the same as the sphere, and $Q=Q_z(c)Q_y(b)Q_x(a)$ is the composition of the three rotation matrices. The terms $Q_x(a)$, $Q_y(a)$, and $Q_z(a)$ are the matrices corresponding to a rotation by an angle $a$ around the $x$, $y$, and $z$ axes respectively.}

The parameters used for the rotations were: 
\begin{align*}
a&=\texttt{0.2440241225550843}\\
b&=\texttt{0.7454097947651017}\\
c&=\texttt{0.2219760487439292}\cdot \texttt{10}^{\texttt{1}}
\end{align*}

Of course to test our algorithms, we retain no information about the parameterizations. The test sphere and torus are represented only by $d_\Gamma$ and $P_\Gamma$ on the given grid.
Figure~\ref{fig:tiltedtorusplot} shows the torus that we use and the points used in the quadrature rule in a configuration. 
The Jacobian $J_\Gamma$ is approximated using a fourth-order centered differencing of $P_\Gamma$ on the grid, see \cite{tsaikublik16}.

\begin{figure}[ht]
    \begin{center}
    \includegraphics{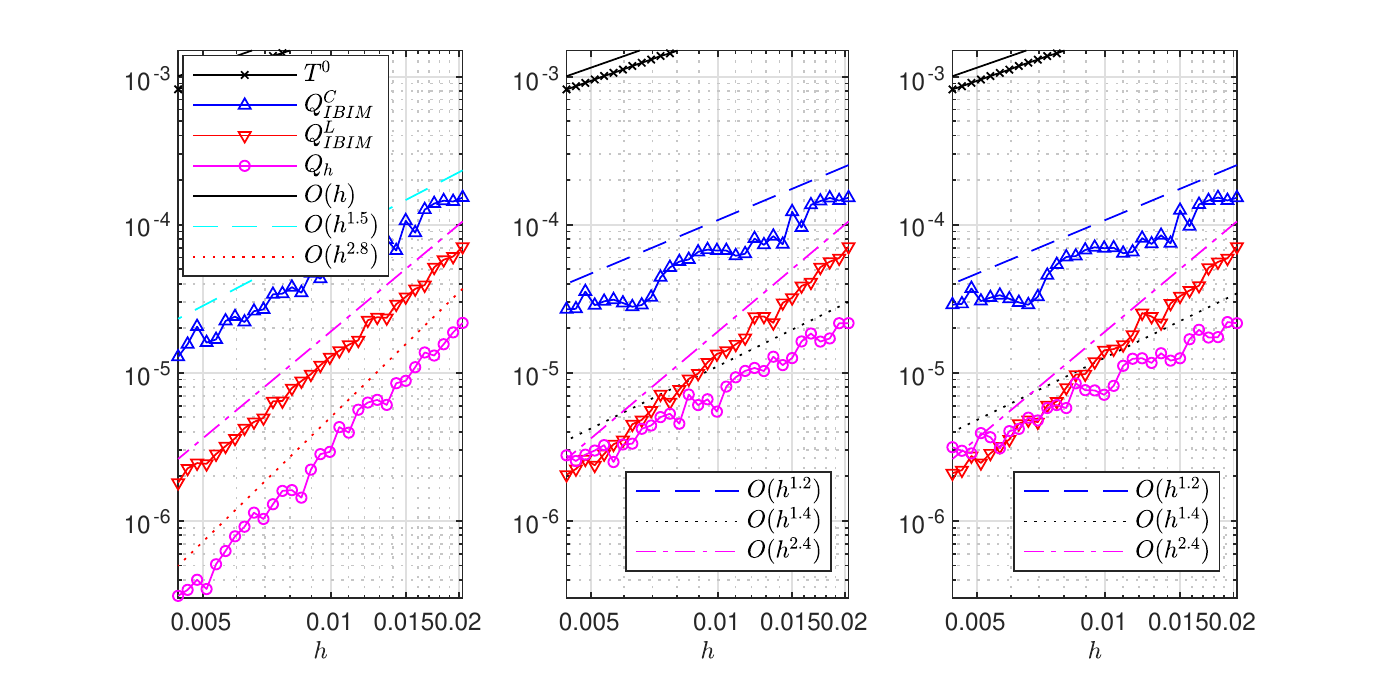}
    \end{center}
    \caption{Sphere tests}
    \footnotesize{Averaged {$E_1(h)$} errors { \eqref{eq:constant_density_error}} on 50 random target points on a sphere. Errors for the punctured trapezoidal rule \eqref{eq:puncturedtrapez} (black crosses) and the three considered methods in the evaluation of the double-layer potential: $Q^C$ constant regularization \eqref{eq:KregIBIM-constant} (blue upward triangles); $Q^L$ \textit{cappuccio} regularization \eqref{eq:KregIBIM-linear} (red downward triangles); $Q_h$ corrected trapezoidal rule \eqref{eq:Q1_3D_additivesplit} (magenta circles). The three plots reflect three different settings for the tubular neighborhood width: left plot $\varepsilon = 0.1$; center plot $\varepsilon\sim h^{0.7}$; right plot $\varepsilon\sim h^{0.8}$.}
    \label{fig:sphere_convergence}
\end{figure}

\subsection{Convergence studies}
We compare the numerical orders of convergence for the quadrature rules discussed in this paper. The quadratures are defined on the grid nodes in $T_\varepsilon$. The parameter $\varepsilon$, which describes the width of the tubular neighborhood, comes into play through the function $\delta_{\Gamma,\varepsilon}$, and the factor $h/\varepsilon$ determines the number of grid nodes in a cross section of the tubular neighborhood $T^h_\varepsilon$. It consequently determines how well the integrand is resolved. The errors for the proposed quadratures are formally $\mathcal{O}\left( \left(\frac{h}{\varepsilon}\right)^p \right)$,  $p\ge 2$, as $h\rightarrow 0$.  Thus, for fixed $\varepsilon= \mathcal{O}(1)$, we see the order of convergence resembling $p$.
If we choose $\varepsilon\sim h^{1/q}$, we will formally have the errors scale as $\mathcal{O}(h^{p(1-1/q)})$. If $\varepsilon\sim h$, then the method will not converge formally, but in the range of the grid resolution considered in practice, the method may yield results with acceptable accuracy. 

\begin{figure}[ht]
    \begin{center}
    \includegraphics{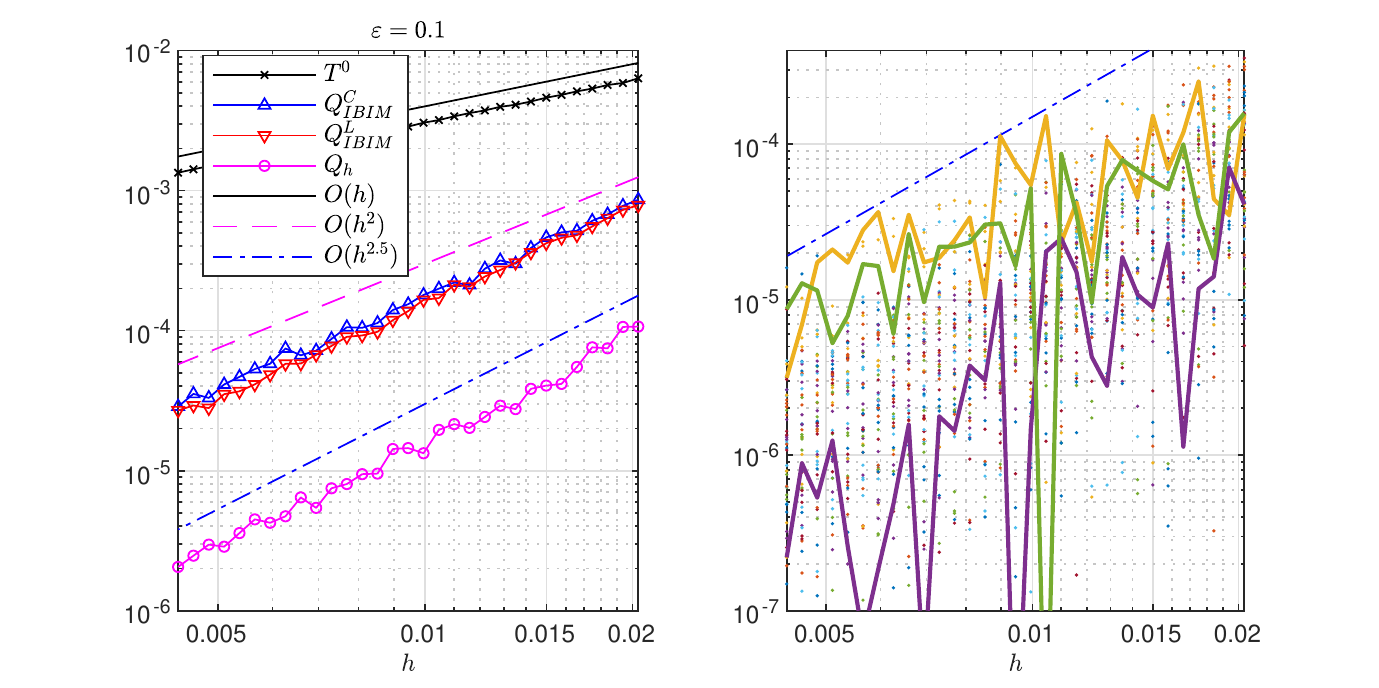}
    \end{center}
    \caption{Torus tests - 1}
    \footnotesize{Tilted torus with 50 random target points on the surface. Tubular neighborhood width constant with respect to $h$: $\varepsilon=0.1$. Left figure: mean error {$E_1(h)$} plotted for the four different methods considered (punctured trapezoidal rule \eqref{eq:puncturedtrapez} (black crosses), constant regularization \eqref{eq:KregIBIM-constant} (blue upward triangles), linear regularization \eqref{eq:KregIBIM-linear} (red downward triangles), and corrected trapezoidal rule \eqref{eq:Q1_3D_additivesplit} (magenta circles). 
    Right figure: distribution of the {$E_1(h)$} errors for the 50 target points for the corrected trapezoidal rule. The yellow, green and purple convergence lines correspond to three of the randomly generated target points: they correspond respectively to the parameters $(\theta,\phi)=($\texttt{0.674795533436653},\texttt{1.5287503395568336}$)$, $($\texttt{5.5902567180364535},\texttt{3.0915183172680867}$)$, $($\texttt{3.0463292511788698},\texttt{5.738447188350594}$)$.}
    \label{fig:tiltedtorusDL}
\end{figure}

\begin{figure}[ht]
    \begin{center}
    \includegraphics{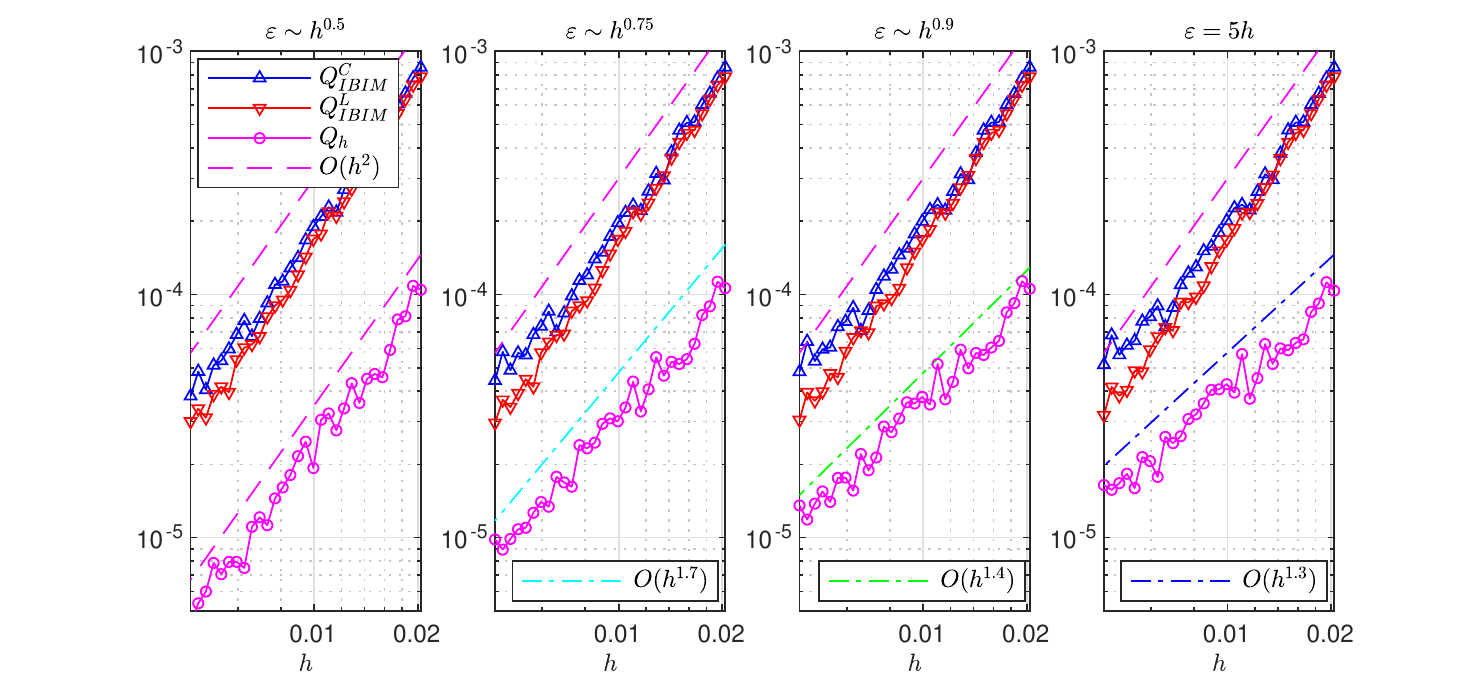}%\\
    \end{center}
    \caption{Torus tests - 2}
    \footnotesize{Tilted torus with 50 random target points; tubular neighborhood width $\varepsilon$ dependent on $h$. Mean error for the three considered methods in the evaluation of the double-layer potential. Top plots: $Q^C_{IBIM}$ constant regularization (blue upward triangles); $Q^L_{IBIM}$ linear regularization, \textit{cappuccio} (red downward triangles); $Q_h$ corrected trapezoidal rule (magenta circles). From left to right, the four plots represent: $\varepsilon \sim h^{0.5}$, $\varepsilon \sim h^{0.75}$, $\varepsilon \sim h^{0.9}$, $\varepsilon = 5h$.}
    \label{fig:torus_convergence}
\end{figure}

In Figures \ref{fig:sphere_convergence}, \ref{fig:tiltedtorusDL}, and \ref{fig:torus_convergence} the errors are shown as function of the Cartesian grid's spacing $h$; the first figure shows the errors for the sphere, while the others show them for the tilted torus. 
We show also how the errors scale for $\varepsilon\sim h^0$ (left plot in Figure~\ref{fig:sphere_convergence}, and Figure \ref{fig:tiltedtorusDL}), and $\varepsilon\sim h^\alpha$ for different $\alpha$s (other plots in Figure \ref{fig:sphere_convergence}, and Figure \ref{fig:torus_convergence}).

In Figure~\ref{fig:tiltedtorusDL} we show the error curves for three target points; affected by their relative positions to the grid and the surface, the errors at some target point is larger than at others. Furthermore, the errors at each target point oscillate as one varies $h$. 

We now show that the accuracy of the weights used in the tests is sufficient for the discretization used. In the previous tests, the number of terms in the Fourier expansion was $2N+1$ with $N=22$, and each term was tabulated in $\alpha,\beta$ with $N_{\alpha,\beta}=101$ values each. We repeated the same test as in Figure~\ref{fig:tiltedtorusDL} for the smallest $h\approx$\texttt{0.00437}, first decreasing $N$ to $11$, and then decreasing $N_{\alpha,\beta}$ to $51$. The corresponding results are in the following table:

\begin{center}
\begin{tabular}{ c r c  } 
    \hline
    $N$ & $N_{\alpha,\beta}$ & avg. error \\ \hline
    22 & 101 & \texttt{2.05289}$\cdot\texttt{10}^{\texttt{-6}}$ \\ 
    22 & 51 & \texttt{2.05290}$\cdot\texttt{10}^{\texttt{-6}}$ \\
    11 & 101 & \texttt{2.05277}$\cdot\texttt{10}^{\texttt{-6}}$ \\ \hline
\end{tabular}
\end{center}
The table suggests that for this range of parameters, the error from the correction of trapezoidal rule is dominating. 

For $\varepsilon\sim\mathcal{O}(1)$, the formal order of convergence for the proposed quadrature is $\mathcal{O}(h^2)$, but we have observed a rate of $\mathcal{O}(h^{2.5})$. In the next subsection, we present a property of our quadrature that we believe leads to the increase in accuracy.

\subsection{Order increase from error cancellation}\label{sec:error-cancellation}

\begin{figure}[ht]
    \begin{center}
    \includegraphics{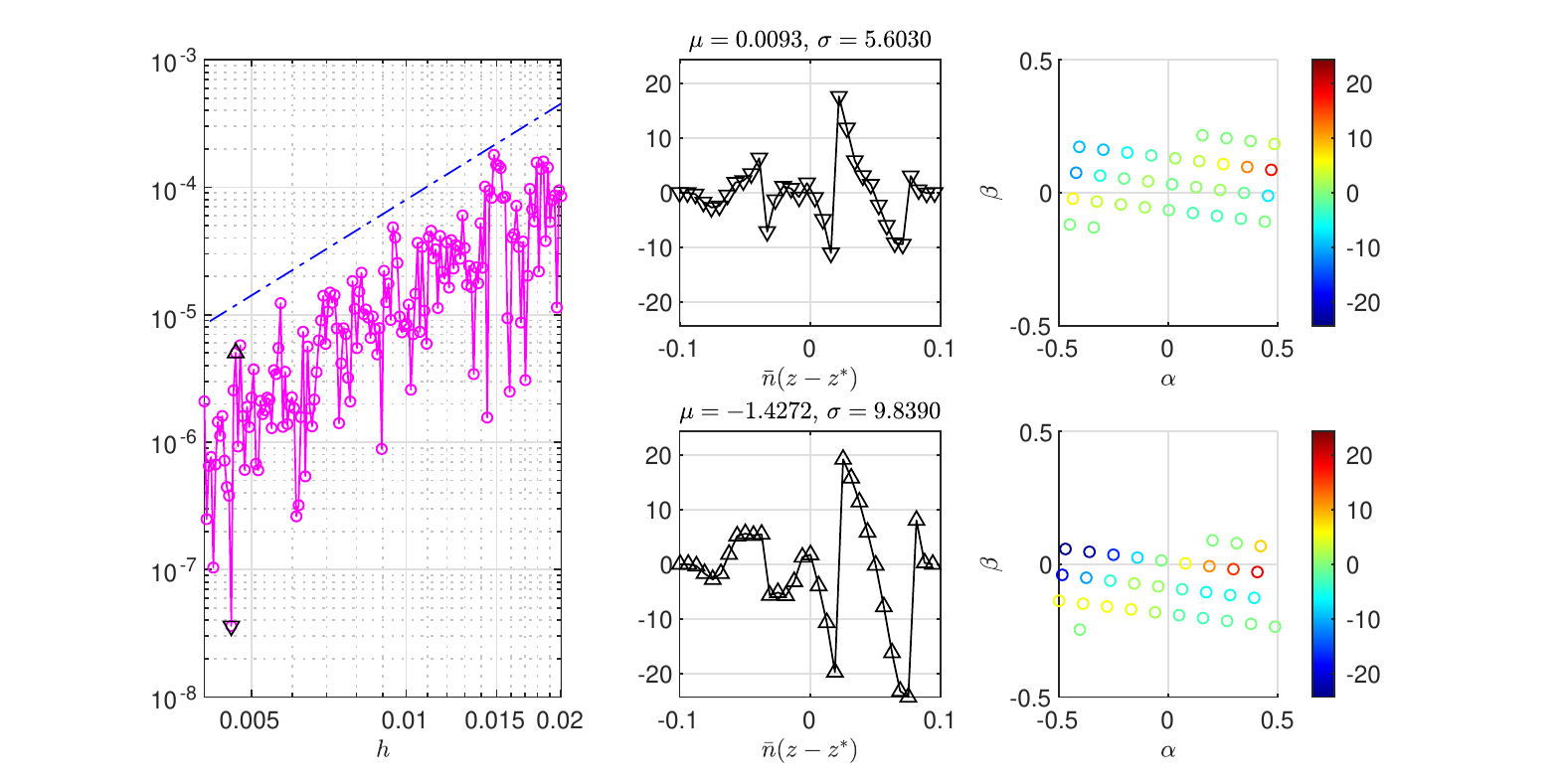}
    \end{center}
    \caption{Error behavior and nonsmooth convergence}
    \footnotesize{Single target point on a ``flat'' torus. In the left half of the figure, it is shown the convergence behavior of the error for the double-layer potential, with two discretizations highlighted (upward and downward triangles). In the right half of the figure, top plots show the error behavior for the discretization $h\approx 4.57\cdot 10^{-3}$ (downward triangle in the left figure) on the corrected planes as function of $\eta(z)=\bar n(z-z^*)$ with $\bar n=\sqrt{n_1^2+n_2^2+1}$, which corresponds to the error $3.61\cdot 10^{-8}$. Bottom plots show the error behavior for the discretization $h\approx 4.66\cdot 10^{-3}$ (upward triangle in the left figure) on the corrected planes as function of $\eta(z)$
    which corresponds to the error $5.03\cdot 10^{-6}$. The left figures show the signed error as function of $\eta(z)$; the right plot shows the distribution of the $(\alpha(h,k),\beta(h,k))$ values on the different planes, where the color represents the value of the error. From the mean $\mu$ and variance $\sigma$ printed on top of the left plots we can see that the top case has mean much smaller than the bottom case, and the variance is half. This explains the much smaller error (downward triangle) compared to the other (upward triangle). }
    \label{fig:goodandbadhexample}
\end{figure}

As discussed in Section~\ref{sub:corr_3D}, our corrected trapezoidal rule is applied on every plane in $T_\varepsilon$ (see Figure~\ref{fig:visualization_normal_surface_grid}), and the error for the whole integral is a sum of the quadrature errors on each relevant plane. 

For a fixed target point the singular line intersects each plane at a different
location relative to the grid. The relative positions are given by the shift parameters $\alpha$ and $\beta$, which depend on both the plane's $z$-coordinate and the grid spacing $h$. 
More precisely, let $\mathbf{\bar y}_0(z)=(x_0(z),y_0(z))$. Then
$$
   \alpha = \left\{\frac{x_0(z)}{h}+\frac12\right\}-\frac12,
   \qquad
   \beta = \left\{\frac{y_0(z)}{h}+\frac12\right\}-\frac12,
$$
where $\{x\}$ denotes the fractional part of $x$. 
Defining the 1-periodic function $r(x) = \{x+1/2\}-1/2$
and using the expression for the singular line
$\mathbf{\bar y}_0(z)$ 
in \eqref{eq:y0(z)_expression} where $\mathbf{n}=(n_1,n_2,1)$, we can write
$$
\alpha = \alpha(h,z/h) = r\left(\frac{z\,n_1+x_0(0)}{h}\right),
\qquad
\beta = \beta(h,z/h) = r\left(\frac{z\,n_2+y_0(0)}{h}\right).
$$
This shows that the relative
location $(\alpha,\beta)$ may vary rapidly between
planes both in $h$ and $z$ when $h$ is small,
and since $r$ is discontinuous, the
variation is non-smooth. 

\begin{figure}[ht]
    \begin{center}
    \includegraphics{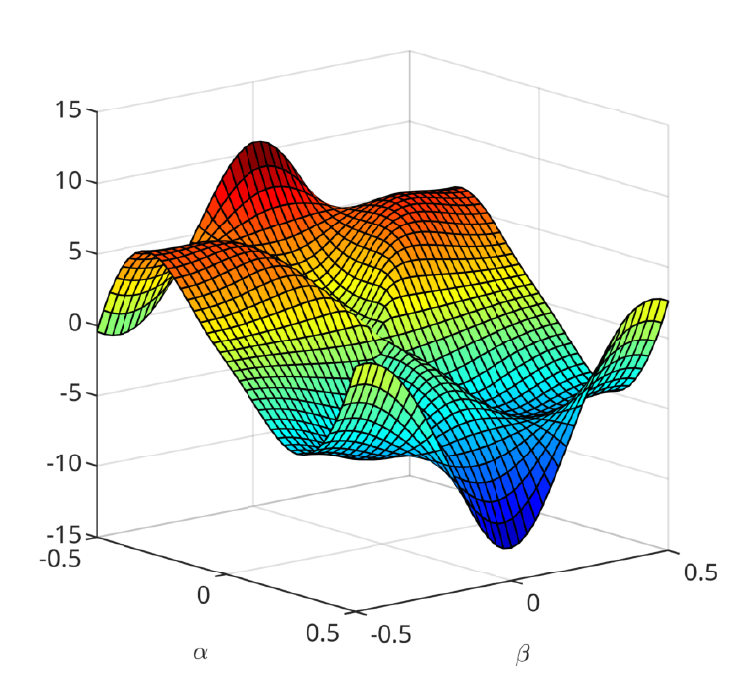}\includegraphics{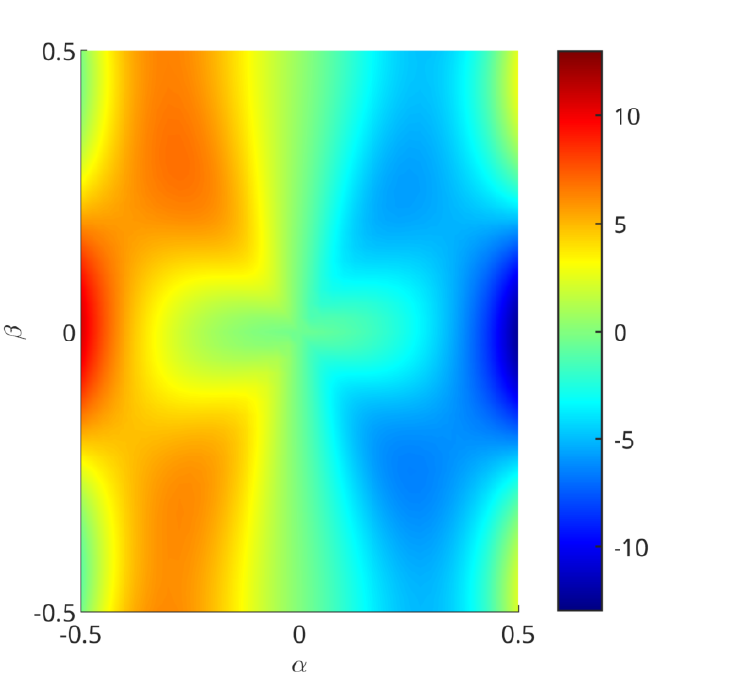}
    \end{center}
    \caption{Error as function of $\alpha$ and $\beta$}
    \footnotesize{Error $F_1(\alpha,\beta,\eta)$ seen for $\eta\approx-0.056$; the singularity line has direction defined by the spherical coordinates $(\theta,\phi)\approx(2.298,3.154)$. The mean over $\alpha$ and $\beta$ is $0.01625$.}
    \label{fig:errorsurface}
\end{figure}

Let $E(z,h)$ be the quadrature error for one plane.
Based on \eqref{eq:error_2D} we can express it as
\begin{equation}
    E(h,z) = F_1(\alpha(h,z/h),\beta(h,z/h),z)\,h^2  +\mathcal{O}(h^3),
\end{equation}
where $F_1$ is a 
smooth function of $(\alpha,\beta,z)$.
In Figure \ref{fig:errorsurface} we can see the function $F_1(\alpha,\beta,z)$ for a specific value of $z$.
The total error for the three-dimensional integral, is then
$$
  E_{\rm tot}(h) 
  =\sum_k h E(h,z_k)
  =D(h)h^2  +\mathcal{O}(h^3),
$$
where, noting that $z_k/h=k$,
$$
D(h) = \sum_k h
  F_1(\alpha(h,k),\beta(h,k),z_k).
$$
The coefficient $D(h)$
is thus the mean of the error
coefficients on the
different planes. Since $F_1$ is smooth, and evaluated
in a compact set, $D$ is therefore
bounded in $h$. 
However, $\alpha$
and $\beta$ are non-smooth in
$h$ and underresolved in the second argument in the sum.
Therefore,
$D$ is not
a smooth function of $h$.
This accounts for the irregular
convergence plots. 
See for instance the left subplot in
Figure~\ref{fig:goodandbadhexample} or the right one in Figure~\ref{fig:tiltedtorusDL}.

The analysis above would
predict second order accuracy
for the method, 
when $\varepsilon$ is independent of $h$. 
However, in
practical computations we typically
observe the higher order convergence
rate $2.5$. 
We believe this can be explained by 
a further property of $F_1$.
Looking at
Figure \ref{fig:errorsurface}, we may notice a skew-symmetry in $F_1(\alpha,\beta, z)$ for fixed $z$.
%a certain symmetry 
%is noticeable 
%in the shape of the function, and 
It is reasonable to expect that the average value of $F_1$ over $\alpha$ and $\beta$ is much smaller in module compared to the maximum error.
In fact, we conjecture that
the average value of $F_1$
for fixed $z$ is zero,
$$
  \int_{-1/2}^{1/2}   \int_{-1/2}^{1/2} F_1(\alpha,\beta,z)\text{d}\alpha \text{d}\beta = 0.
$$
In the sum of $F_1(\alpha(h,k),\beta(h,k),z_k)$,
defining $D(h)$ 
the first two arguments ($\alpha,\beta$) vary much faster than the third ($z_k$). If the sequence
$k\mapsto (\alpha(h,k),\beta(h,k))$ has some ergodic property
the sum will behave similar to the full integral,
which would be zero
$$
D(h)\approx    \int
\int_{-1/2}^{1/2}   \int_{-1/2}^{1/2} 
  F_1(\alpha,\beta,z)\text{d}\alpha \text{d}\beta \text{d}z = 0.
$$
A precise analysis of this effect is beyond the scope of
this article. Here we just show in Figure \ref{fig:goodandbadhexample}
an example of how
$F_1(\alpha(h,k),\beta(h,k),z_k)$ 
and $(\alpha(h,k),\beta(h,k))$
may vary for two different $h$ that are very close
to each other.

%\newpage
\section*{Acknowledgment}
Tsai's research is supported partially by National Science Foundation Grants DMS-1720171 and DMS-1913209.
Part of this research was performed while the third author was visiting the Institute for Pure and Applied Mathematics (IPAM), which is supported by the National Science Foundation (Grant No. DMS-1440415). This work was partially supported by a grant from the Simons Foundation.

\clearpage{}

\appendix

\section{Appendix}

\subsection{Relating curvatures and the principal directions on parallel surfaces}

Let $\Omega\subset\mathbb{R}^{3}$ be a bounded domain, and let $d_\Gamma$
and $P_{\Gamma}$ be the signed distance function and the closest point projection defined in \eqref{eq:signed-distance} and \eqref{eq:closest-point-projection} in Section~\ref{subsec:ImplicitParametrization}. 
We assume that the distance function is twice continuously differentiable in the tubular neighborhood $T_{\varepsilon}:=\left\{ \mathbf{z}\in\mathbb{R}^{3}\,:\  |d_\Gamma(\mathbf{z})|<\varepsilon\right\}.$

The derivation of the proposed quadratures relies heavily on the knowledge of geometrical information of the surface $\Gamma$, through that of the level sets of $d_\Gamma$.
In this Section, we relate the principal curvatures and the corresponding directions on different parallel surfaces $\Gamma_{\eta}:=\left\{ \mathbf{z}\in T_{\varepsilon}\,:\ {d}_\Gamma(\mathbf{z})=\eta\right\},$ for $\eta\in[-\varepsilon,\varepsilon].$

Let $\mathbf{z}$ be an arbitrary point in $T_\varepsilon,$ and $\eta=d_\Gamma(\mathbf{z}).$
The curvature information of $\Gamma_\eta$ at $\mathbf{z}$
can be retrieved from the Hessian of $d_\Gamma$. Through eigenvalue decomposition, we have
\[
H_{{d}_\Gamma}(\mathbf{z})=\nabla^{2}{d}_\Gamma(\mathbf{z})=\left[\begin{array}{ccc}
\mathbf{n} & \bm{\tau}_{1} & \bm{\tau}_{2}\end{array}\right]\begin{bmatrix}0\\
 & -\bar{\kappa}_{1}\\
 &  & -\bar{\kappa}_{2}
\end{bmatrix}\left[\begin{array}{ccc}
\mathbf{n} & \bm{\tau}_{1} & \bm{\tau}_{2}\end{array}\right]^{T}
\]
where $\bar{\kappa}_{1}$ and $\bar{\kappa}_{2}$ are the principal
curvatures of $\Gamma_{\eta}$ at $\mathbf{z}$ and $\bm{\tau}_1$, $\bm{\tau}_2$ the corresponding principal directions.  

One can derive easily that the following formula, relating the principle curvatures $\kappa_i$ of $\Gamma$ at $P_\Gamma\mathbf{z}$ and $\bar\kappa_i$ of $\Gamma_\eta$ at $\mathbf{z}$:
\[
-\kappa_{i}=\frac{-\bar{\kappa}_{i}}{1+{d}_\Gamma(\mathbf{z})\bar{\kappa}_{i}},~~~i=1,2.
\]
See for example \cite{gilbarg2015elliptic} (\emph{
$\S$14.6 Appendix: Boundary Curvatures and Distance Function}).
The principal directions will remain the same: 
\begin{lem}
Let $\Gamma$ be a $C^{2}$ surface, $\mathbf{\bar{z}}\in\Gamma$; let $\Gamma_{\eta}$
be a parallel surface, and $\mathbf{z}\in\Gamma_{\eta}$ such that $\mathbf{\bar{z}}=P_{\Gamma}\mathbf{z}$.
The principal directions at $\mathbf{z}$ coincide with the principal directions
at $\mathbf{\bar{z}}$.\label{principaldirectionsLemma}
\end{lem}

\begin{proof}
The tangent plane $TM_\eta(\mathbf{z})$ for $\Gamma_\eta$ at $\mathbf{z}$ is parallel
to the tangent plane $TM_0(\mathbf{\bar{z}})$ for $\Gamma$ at $\mathbf{\bar{z}}$. 

Let $(\mathbf{b}_{1},\mathbf{b}_{2})$
be an orthonormal basis for the plane $TM_0(\mathbf{\bar{z}})$, and $\mathbf{v}=\mathbf{b}_{1}\cos\theta+\mathbf{b}_{2}\sin\theta$
a unit vector.% on $TM_0({\mathbf{\bar{z}}})$. 
We can consider the plane $H_{\mathbf{v}}$
passing though $\mathbf{\bar{z}}$ and parallel to the normal vector, and $H_{\mathbf{v}}\cap\Gamma$
will locally be the support of the regular curve $\gamma_\theta(s)$,
which is the \emph{normal section} of $\Gamma$ at $\mathbf{\bar{z}}$. Corresponding
to the normal section we can calculate the \emph{normal curvature}
$\kappa(\theta)$ of $\Gamma$ at $\mathbf{\bar{z}}$ along $\mathbf{v}$. 

Then $\kappa(\theta)$ is a periodic function in $[0,\pi]$. The minimum
and maximum attained by the curvature are the two \emph{principal
curvatures} $\kappa_{1}:=\min_{\theta}\kappa(\theta)=\kappa(\theta_{1})$
and $\kappa_{2}:=\max_{\theta}\kappa(\theta)=\kappa(\theta_{2})$.
Consequently, $\kappa^{\prime}(\theta_{i})=0$,
$i=1,2$, and $\kappa^{\prime\prime}(\theta_{2})<0<\kappa^{\prime\prime}(\theta_{1})$. Corresponding to these values are two unit vectors on $TM_0({\mathbf{\bar{z}}})$,
which form an orthonormal basis, called \emph{principal directions}.

On $TM_\eta(\mathbf{z})$, we can use the same exact setup, the same basis
$(\mathbf{b}_{1},\mathbf{b}_{2})$, and same unit vectors $\mathbf{v}=\mathbf{b}_{1}\cos\theta+\mathbf{b}_{2}\sin\theta$.
We know that the curvatures will be transformed via the relation
\[
\bar{\kappa}(\theta)=\frac{\kappa(\theta)}{1-\eta\kappa(\theta)}\,.
\]
The maximum and minimum values of this function are going to be again
$\theta_{1}$ and $\theta_{2}$, as
\begin{align*}
\bar{\kappa}^{\prime}(\theta_{i}) & =\frac{\kappa^{\prime}(\theta_{i})}{\left[1-\eta\kappa(\theta_{i})\right]^{2}}=0\ \,,\ \ i=1,2\\
\bar{\kappa}^{\prime\prime}(\theta_{i}) & =\frac{\kappa^{\prime\prime}(\theta_{i})}{\left[1-\eta\kappa(\theta_{i})\right]^{2}}\,;
\end{align*}
then the values $\theta_{i}$, $i=1,2$ are extrema also for this case, and the second derivatives have the same sign as the ones on $\Gamma$. Consequently the angles at which maximum and minimum are attained
are the same, and the principal directions on the two parallel surfaces coincide.
\end{proof}

\subsection{Calculation of the regularizations of the double-layer kernels}
\label{sub:AppendixRegularization}
Given a target point $\mathbf{x}\in\Gamma$, and $r_0>0$, let $\mathcal{M}(\mathbf{x},r_0)\subset\Gamma$ be a neighborhood of $x$ dependent on the parameter $r_0$; we will define it more clearly later. We want to find $\Psi_{\Gamma, r_0}$ such that
\[
\int_{\mathcal{M}(x,r_{0})}\frac{1}{4\pi}\frac{(\mathbf{x}-\mathbf{y})^T \mathbf{n}_{y}}{\|\mathbf{x}-\mathbf{y}\|^{3}}\text{d}\sigma_{\mathbf{y}}=\int_{\mathcal{M}(\mathbf{x},r_{0})}\frac{\partial G_{0}}{\partial \mathbf{n}_{y}}(\mathbf{x},\mathbf{y})\text{d}\sigma_{\mathbf{y}}=\int_{\mathcal{M}(\mathbf{x},r_{0})}\Psi_{\Gamma, r_0}(\mathbf{x},\mathbf{y})\text{d}\sigma_\mathbf{y}\,.
\]

We approximate the surface using a paraboloid, and we assume
the surface is positioned with the target point in the origin $\mathbf{x}=\mathbf{0}$,
and the normal in the target placed along the $z$-axis: $\mathbf{n}_{\mathbf{x}}=(0,0,1)$. Given
the principal curvatures of the surface in the target point $\mathbf{x}$,
$\kappa_{1}$ and $\kappa_{2}$, and assuming the corresponding principal directions lie along the $x$-axis and $y$-axis respectively, the surface is described around the
origin as as $\mathcal{P}(x,y):=(x,y,z(x,y))$, with $z$ such that
\begin{align*}
z(0,0) & = \frac{\partial z}{\partial x}(0,0)=\frac{\partial z}{\partial y}(0,0)=\frac{\partial^2 z}{\partial x\partial y}(0,0)=0,
\end{align*}
and
\begin{equation*}
    \ \frac{\partial^2 z}{\partial x^2}(0,0)=\kappa_{1},\ \frac{\partial^2 z}{\partial y^2}(0,0)=\kappa_{2}\,.
\end{equation*}

The paraboloid is the surface $\tilde \Gamma$ defined for points
close to the origin with coordinates $\left(x,y,z(x,y)\right)$ with
$z(x,y):=\frac{1}{2}(\kappa_{1}x^{2}+\kappa_{2}y^{2})$. The paraboloid
$\tilde \Gamma$ approximates the surface $\Gamma$ with errors of the
third order: $\mathcal{O}(x^{3},y^{3},x^{2}y,xy^{2})$. The Jacobian
$J(x,y)$ is going to be the norm of the normal vector to the surface
$$
J(x,y)=\sqrt{1+\left(\frac{\partial z}{\partial x}(x,y)\right)^2+\left(\frac{\partial z}{\partial y}(x,y)\right)^2}=\sqrt{1+\kappa_1^2x^2+\kappa_2^2y^2}\,.
$$
By using this approximation of the surface and considering as the neighborhood $\mathcal{M}(\mathbf{x},r_0)$ the
set $$\mathcal{M}_{r_0}:=\left\{ \mathcal{P}(x,y)\ :\ \sqrt{x^{2}+y^{2}}\leq r_{0} \right\}$$
we can rewrite the integral as
\[
\int_{\mathcal{M}_{r_{0}}}\frac{\partial G_{0}}{\partial \mathbf{n}_{y}}(\mathbf{0},\mathcal{P}(x,y))\text{d}\sigma_{x,y}=\int_{\mathcal{M}_{r_{0}}}F(x,y)J(x,y)\text{d}x\text{d}y, 
\]
where 
\[
F(x,y):=\frac{1}{8\pi}\frac{\kappa_{1}x^{2}+\kappa_{2}y^{2}}{\left[x^{2}+y^{2}+\frac{1}{4}\left(\kappa_{1}x^{2}+\kappa_{2}y^{2}\right)^{2}\right]^{\frac{3}{2}}\sqrt{1+\kappa_{1}^{2}x^{2}+\kappa_{2}^{2}y^{2}}}\ .
\]

In the article \cite{kublik2013implicit} the function $\Psi_{\Gamma, r_0}(x,y)=C_{\Gamma,r_0}$ is defined as a constant with respect to $x$ and $y$:
\[
\int_{\mathcal{M}_{r_{0}}}\frac{\partial G_{0}}{\partial \mathbf{n}_{y}}(\mathbf{0},\mathcal{P}(x,y))\text{d}\sigma_{x,y}\approx\int_{\mathcal{M}_{r_{0}}}C_{\Gamma,r_{0}}\,\text{d}\sigma_{x,y}\,.
\]
The constant $C_{\Gamma,r_{0}}$ represents the average of the integrand
over $\mathcal{M}_{r_0}$. From elementary calculation, we have
\begin{align*}
    C_F^{DL}(r_0) &= \int_{\mathcal{M}_{r_0}}F(x,y)J(x,y)\text{d}x\text{d}y\\
    &=\int_0^{2\pi}\text{d}\theta\int_{0}^{r_0}\text{d}r\left\{ r\,F(r\cos\theta,r\sin\theta)J(r\cos\theta,r\sin\theta) \right\}\\
    &= \frac{\kappa_1+\kappa_2}{8}r_0+\frac{\kappa_1+\kappa_2}{512}(-5\kappa_1^2+2\kappa_1\kappa_2-5\kappa_2^2)r_0^3+\mathcal{O}(r_0^5)\,,\\
    C_{\Gamma}(r_0) &= \int_{\mathcal{M}_{r_0}}J(x,y)\text{d}x\text{d}y=\int_0^{2\pi}\text{d}\theta\int_{0}^{r_0}\text{d}r\left\{ r\,J(r\cos\theta,r\sin\theta) \right\}\\
    &= \pi r_0^2 +\frac{\pi}{8}(\kappa_1^2+\kappa_2^2)r_0^4+\mathcal{O}(r_0^6)\,.
\end{align*}
Then 
\[
C_{r_{0}} = \frac{C_F^{DL}(r_0)}{C_\Gamma(r_0)} = \frac{\kappa_{1}+\kappa_{2}}{8\pi r_{0}}-\frac{\kappa_{1}+\kappa_{2}}{512\pi}\left(13\kappa_{1}^{2}-2\kappa_{1}\kappa_{2}+13\kappa_{2}^{2}\right)r_{0}+\mathcal{O}(r_0^3).
\]
Finally, $\frac{\partial G_{0}}{\partial \mathbf{n}_{y}}(\mathbf{x}-\mathbf{y})$,
$\mathbf{x},\mathbf{y}\in\Gamma$ can then be regularized as:
\begin{equation}\label{eq:IBIM_C_DL}
K^{reg,DL}_{r_0,C}(\mathbf{x},\mathbf{y}):=
\begin{cases}
\dfrac{\partial G_{0}}{\partial \mathbf{n}_{y}}(\mathbf{x},\mathbf{y})\,, & \|\mathbf{x}-\mathbf{y}\|\geq r_0\,,\\[0.4cm]
C_{\Gamma,r_0}\,, & \|\mathbf{x}-\mathbf{y}\|<r_0\,.
\end{cases}
\end{equation}

The same reasoning can be applied to the double-layer conjugate kernel, where in the previous calculations the expression of $F$ is
\[
F(x,y):=\frac{1}{8\pi}\frac{\kappa_{1}x^{2}+\kappa_{2}y^{2}}{\left[x^{2}+y^{2}+\frac{1}{4}\left(\kappa_{1}x^{2}+\kappa_{2}y^{2}\right)^{2}\right]^{\frac{3}{2}}}\,,
\]
and the result is the following regularization:
\begin{equation}\label{eq:IBIM_C_DLC}
K^{reg,DLC}_{r_0,C}(\mathbf{x},\mathbf{y}):=
\begin{cases}
\dfrac{\partial G_{0}}{\partial \mathbf{n}_{x}}(\mathbf{x},\mathbf{y})\,, & \|\mathbf{x}-\mathbf{y}\|\geq r_0\,,\\[0.4cm]
C^{DLC}_{\Gamma,r_0}\,, & \|\mathbf{x}-\mathbf{y}\|<r_0\,,
\end{cases}
\end{equation}
where
\[
C^{DLC}_{\Gamma,r_0}=\frac{\kappa_1+\kappa_2}{8\pi r_0}-\frac{5}{1536}\frac{\kappa_1+\kappa_2}{\pi}(3\kappa_1^2+2\kappa_1\kappa_2+3\kappa_2^2)r_0+\mathcal{O}(r_0^3)\,.
\]

For the case of the secondary kernel of the Helmholtz equation, the function $F$ is
\[
F(x,y):=\frac{1}{2}\frac{\kappa_{1}x^{2}+\kappa_{2}y^{2}}{{x^{2}+y^{2}+\frac{1}{4}\left(\kappa_{1}x^{2}+\kappa_{2}y^{2}\right)^{2}}}\,,
\]
and the regularization becomes
\begin{equation}\label{eq:IBIM_C_Helmholtz_secondary}
K^{reg,HL}_{r_0,C}(\mathbf{x},\mathbf{y}):=
\begin{cases}
\dfrac{(\mathbf{x}-\mathbf{y})^T\mathbf{n}_y}{\|\mathbf{x}-\mathbf{y}\|^2}\,, & \|\mathbf{x}-\mathbf{y}\|\geq r_0\,,\\[0.4cm]
C^{HL}_{\Gamma,r_0}\,, & \|\mathbf{x}-\mathbf{y}\|<r_0\,,
\end{cases}
\end{equation}
where
\[
C_{\Gamma,r_0}^{HL}=\frac{\kappa_1+\kappa_2}{4}-\frac{\kappa_1+\kappa_2}{256}(13\kappa_1^2-2\kappa_1\kappa_2+13\kappa_2^2)r_0^2+\mathcal{O}(r_0^4)\,.
\]

\subsubsection*{New regularization with linear function (\emph{cappuccio})}
An potential improvement on \eqref{eq:KregIBIM-constant} can be made by building $\Psi_{\Gamma, r_0}$ as linear
with respect to the distance from the singularity $\Psi_{\Gamma, r_0}(\mathbf{x},\mathbf{y})=\Psi_{\Gamma,r_0^L}(\mathbf{x},\mathbf{y}):=a_0\,\frac{\|\mathbf{x}-\mathbf{y}\|}{r_0}+a_1$:
\[
\int_{\mathcal{M}_{r_{0}}}F(x,y)J(x,y)\text{d}x\text{d}y=\int_{\mathcal{M}_{r_{0}}}\left(a_0\,\frac{\|\mathbf{0}-\mathcal{P}(x,y)\|}{r_0}+a_1\right)J(x,y)\,\text{d}x\text{d}y\,.
\]
The second property we impose is the following: we express $\mathcal{P}(x,y)=\mathcal{P}(r\cos\theta,r\sin\theta)$ in polar coordinates, and impose
\[ 
\Psi_{\Gamma,r_0}^L(r_0\cos\theta,r_{0}\sin\theta)=\frac{1}{2\pi}\int_{0}^{2\pi}F(r_0\cos\theta,r_0\sin\theta)\text{d}\theta\,.
\]
We call:
\begin{align*}
    C_{r^2} &= \int_{\mathcal{M}_{r_0}}\|\mathbf{0}-\mathcal{P}(x,y)\|J(x,y)\text{d}x\text{d}y \\
    C_{F}^{DL} &= \int_{\mathcal{M}_{r_0}}F(x,y)J(x,y)\text{d}x\text{d}y\\
    C_\Gamma &= \int_{\mathcal{M}_{r_0}}J(x,y)\text{d}x\text{d}y\\
    \phi_0 &= \frac{1}{2\pi}\int_{0}^{2\pi}F(r_0\cos\theta,r_0\sin\theta)\text{d}\theta
\end{align*}
then the conditions imposed form the following linear system:
\begin{align*}
    &a_0+a_1=\phi_0 \\
    &\frac{a_0}{r_0}C_{r^2}+a_1C_\Gamma=C_F^{DL}
\end{align*}
from which we find:
\begin{align*}
    a_0 &= \frac{C_F^{DL}-\phi_0C_\Gamma}{C_{r^2}-r_0C_\Gamma}r_0  \\
    &= -\frac{3}{16}\frac{\kappa_1+\kappa_2}{\pi r_0}+\frac{3}{5120}\frac{\kappa_1+\kappa_2}{\pi}(21\kappa_1^2-2\kappa_1\kappa_2+21\kappa_2^2)r_0+\mathcal{O}(r_0^3)\\
    a_1 &= \frac{\phi_0C_{r^2}-C_F^{DL}r_0}{C_{r^2}-r_0C_\Gamma}\\
    &= \frac{\kappa_1+\kappa_2}{4 \pi r_0}-\frac{3}{2560}\frac{\kappa_1+\kappa_2}{\pi}(23\kappa_1^2-6\kappa_1\kappa_2+23\kappa_2^2)r_0+\mathcal{O}(r_0^3)
\end{align*}
This regularization is then:
\begin{equation}\label{eq:IBIM_L_DL}
K^{reg,DL}_{r_0,L}(\mathbf{x},\mathbf{y}):=
\begin{cases}
\dfrac{\partial G_{0}}{\partial \mathbf{n}_{y}}(\mathbf{x},\mathbf{y})\,, & \|\mathbf{x}-\mathbf{y}\|\geq r_0\,,\\[0.4cm]
a_0\dfrac{\|\mathbf{x}-\mathbf{y}\|}{r_0}+a_1\,, & \|\mathbf{x}-\mathbf{y}\|<r_0\,.
\end{cases}
\end{equation}

\printbibliography
\end{document}